\documentclass[11pt]{article}

\usepackage{amsmath, graphicx, latexsym, amssymb, amsthm, amsfonts}

\usepackage[mathscr]{eucal}
\usepackage{graphics}
\usepackage{color}
\usepackage{epsfig}
\usepackage{bbm}
\usepackage{hyperref}
\usepackage{tikz}

\usepackage{mathtools}
\usetikzlibrary{shapes,arrows}
\usetikzlibrary{graphs,graphs.standard}
\usepackage{pgfplots}
\usepackage{dsfont}
\usepackage{tcolorbox}
\usepackage{caption}
\usepackage{subfig}

\newtheorem{theorem}{Theorem}[section]
\newtheorem{lemma}[theorem]{Lemma}

\newtheorem{definition}[theorem]{Definition}

\newtheorem{proposition}[theorem]{Proposition}
\newtheorem{assumption}[theorem]{Assumption}

\newtheorem{remark}[theorem]{Remark}

\def\thelemma{\arabic{section}.\arabic{lemma}}
\def\thetheorem{\arabic{section}.\arabic{theorem}}
\def\thecorollary{\arabic{section}.\arabic{corollary}}
\def\thedefinition{\arabic{section}.\arabic{definition}}
\def\theexample{\arabic{section}.\arabic{example}}
\def\theproposition{\arabic{section}.\arabic{proposition}}

\def\theassumption{\arabic{section}.\arabic{assumption}}

\def\theremark{\arabic{section}.\arabic{remark}}

\def\theequation{\arabic{section}.\arabic{equation}}

\newcommand{\manualnames}[1]{
	\def\theequation{#1.\arabic{equation}}
	\def\thelemma{#1.\arabic{lemma}}
	\def\thetheorem{#1.\arabic{theorem}}
	\def\thecorollary{#1.\arabic{corollary}}
	\def\thedefinition{#1.\arabic{definition}}
	\def\theexample{#1.\arabic{example}}
	\def\theproposition{#1.\arabic{proposition}}
	\def\theassumption{#1.\arabic{assumption}}
	\def\theremark{#1.\arabic{remark}}
}

\newcommand{\beginsec}{
	\setcounter{equation}{0}
}

\newcommand{\la}{\lambda}
\newcommand{\eps}{\varepsilon}
\newcommand{\ph}{\varphi}

\newcommand{\al}{\alpha}

\newcommand{\s}{\sigma}
\newcommand{\sig}{\sigma}
\newcommand{\del}{\delta}
\newcommand{\om}{\omega}
\newcommand{\Gam}{\mathnormal{\Gamma}}
\newcommand{\Del}{\mathnormal{\Delta}}
\newcommand{\Th}{\mathnormal{\Theta}}
\newcommand{\La}{\mathnormal{\Lambda}}
\newcommand{\X}{\mathnormal{\Xi}}

\newcommand{\Ups}{\mathnormal{\Upsilon}}

\newcommand{\Om}{\mathnormal{\Omega}}

\newcommand{\N}{{\mathbb N}}

\newcommand{\R}{{\mathbb R}}

\newcommand{\Z}{{\mathbb Z}}

\newcommand{\E}{{\mathbb E}}

\newcommand{\PP}{{\mathbb P}}
\newcommand{\BH}{{\mathbb H}}

\newcommand{\calA}{{\cal A}}

\newcommand{\calD}{{\cal D}}

\newcommand{\calF}{{\cal F}}

\newcommand{\calH}{{\cal H}}

\newcommand{\calR}{{\cal R}}
\newcommand{\calS}{{\cal S}}

\newcommand{\scrA}{\mathscr{A}}

\newcommand{\frS}{\mathfrak{S}}
\newcommand{\frX}{\mathfrak{X}}

\renewcommand{\proof}{\noindent{\bf Proof.\ }}

\newcommand{\lan}{\langle}
\newcommand{\ran}{\rangle}

\newcommand{\skp}{\vspace{\baselineskip}}

\newcommand{\w}{\wedge}

\newcommand{\To}{\Rightarrow}

\newcommand{\wt}{\widetilde}

\newcommand{\iy}{\infty}

\newcommand{\be}{\begin{equation}}
\newcommand{\ee}{\end{equation}}

\newcommand{\noi}{\noindent}
\newcommand{\ds}{\displaystyle}
\newcommand{\one}{\mathds{1}}

\newcommand{\slp}{\calS_{\text{\rm LP}}}

\newcommand{\les}{\leqslant}
\newcommand{\ges}{\geqslant}
\newcommand{\p}{${\bf P}$}
\newcommand{\teeone}{{\bf T}${}_1$}
\newcommand{\teetwo}{{\bf T}${}_2$}

\newcommand{\m}{\mathbf{m}}

\begin{document}

\title{Asymptotic optimality of switched control policies
in a simple parallel server system  
under an extended heavy traffic condition}

\author{Rami Atar\thanks{Viterbi Faculty of Electrical and Computer Engineering, Technion, Haifa, Israel}
\and
Eyal Castiel${}^*$
\and Marty Reiman\thanks{Department of Industrial Engineering and Operations Research, Columbia University, New York City, NY}}

\maketitle

\begin{abstract}
	
	This paper studies a 2-class, 2-server parallel server system under the
	recently introduced extended heavy traffic condition \cite{ACR1},
	which states that the underlying 'static allocation' linear program (LP) is critical,
	but does not require that it has a unique solution. The main result is the construction
	of policies that asymptotically achieve a lower bound, proved in \cite{ACR1},
	on an expected discounted linear combination of diffusion-scaled queue lengths,
	and are therefore asymptotically optimal (AO).
	Each extreme point solution to the LP determines a control mode, i.e.,
	a set of activities (class--server pairs) that are operational.
	When there are multiple solutions, these modes can be selected dynamically.
	It is shown that the number of modes required
	for AO is either one or two. In the latter case
	there is a switching point in the (normalized) workload domain,
	characterized in terms of a free boundary problem.
	Our policies are defined by identifying pairs of elementary policies and switching between them
	at this switching point.
	They provide the first example in the heavy traffic literature
	where weak limits under an AO policy are given by
	a diffusion process where both the drift and diffusion coefficients are discontinuous.

	\skp
	
	\noi{\bf MSC 2020 classification:}
	60K25 ; 68M20 ; 93E20 ; 60F17 ; 90B36
	
	\skp
	
	\noi{\bf Keywords:}
	parallel server systems;
	product form service rates;
	extended heavy traffic condition;
	dynamic graph of basic activities;
	switched control systems;
	diffusion with discontinuous coefficients.
\end{abstract}

\section{Introduction}\label{sec1}

\subsection{Background}

Parallel server systems (PSS) are queueing control problems in which
a number of servers offer service to customers
of different classes, and choices
as to which customer class each server is dedicated to are made dynamically.
Since its introduction in \cite{har-lop}, its study in heavy traffic
has attracted much attention due to its simple structure, its practical significance,
and the theoretical challenges it poses.
The problem formulation in \cite{har-lop} includes a key assumption, referred to as
the {\it heavy traffic condition} (HTC), which states that an underlying
'static allocation' linear program (LP) satisfies a critical load condition
and that it has a unique solution. Whereas critical load is
universally considered a defining condition of any notion of heavy traffic,
uniqueness of solutions has been assumed mainly because it simplifies
the mathematical treatment. The {\it extended HTC} (EHTC), which merely
states that the LP is at criticality
but does not require uniqueness, has recently been introduced in \cite{ACR1}
in order to address a considerably broader notion of
heavy traffic.
The main result of \cite{ACR1} is a lower bound on the asymptotically achievable cost in
a general PSS under the EHTC.
%The reader is referred to \cite{ACR1} for further background on this model.
This paper focuses on the $2$-class, $2$-server PSS
referred to in this introduction as the $2\times2$ PSS,
which is the simplest case in which the EHTC is strictly broader than
the HTC. The goal is to complement the results of \cite{ACR1} in this case
by constructing policies that asymptotically achieve the lower bound,
which hence are asymptotically optimal (AO) in heavy traffic.

The structure of the $2\times 2$ PSS is as follows.
Each of the two servers is capable of
serving each of the two classes. The classes (respectively, servers)
are usually indexed using the symbol $i$ (respectively, $k$), and activities,
namely class-server pairs, by $j=(i,k)$.
Arriving customers await service in class-based queues, and upon receiving
a single service, leave the system. The control decisions consist of routing
(determining which server serves each customer) and sequencing
(determining the order in which they are served).
The rates of arrivals of customers of the two classes
are denoted by $\la_i^n$, $i=1,2$, and the rates of service
at each of the four activities are denoted by $\mu^n_{ik}$,
where $n$ denotes the usual heavy traffic parameter.
These rates are assumed to be asymptotic to
$\la_in+\hat\la_in^{1/2}$ and $\mu_{ik}n+\hat\mu_{ik}n^{1/2}$,
for some given $\la_i$, $\hat\la_i$, $\mu_{ik}$, $\hat\mu_{ik}$.
The cost consists of
an expected infinite horizon discounted linear combination
of the two queue lengths, and is rescaled at the diffusion scale.

Whereas the cost, and consequently the notion of AO, are set up at the
diffusion scale, the underlying LP alluded to above addresses
the behavior of the PSS at the fluid, or law-of-large-numbers (LLN) scale.
Posed in terms of the first order parameters,
$\la_i$, $\mu_{ik}$, it is concerned with the mean fraction
of time devoted by each server to each class. When the LP has
a unique solution, at least one activity is {\it non-basic}, in the sense that
the fraction allocated to it is zero. The so called {\it graph of basic activities}
(GBA), formed by the activities with positive allocation fraction, is
static.
In this case, the critical load condition dictates that any policy not adhering to
this solution, in the sense of effort allocation,
causes the total queue length to blow up, and in particular
cannot be AO. Under policies that adhere to this solution,
the LLN assures that the aforementioned fractions of effort converge
to those given by the LP solution (a necessary, but certainly not sufficient
condition for AO).
When there are multiple LP solutions,
a result from \cite{ACR1} states that for the $2\times2$ PSS, the space of
solutions, denoted $\slp$,
forms a line segment ${\rm ch}(\xi^{*,1},\xi^{*,2})$ in the space
of $2\times2$ matrices (where ch denotes the convex hull).
In each of the two extremal solutions, $\xi^{*,1},\xi^{*,2}$
there is again at least one non-basic activity.
We refer to these two extreme points as
{\it control modes}, or simply {\it modes}.
For similar reasons, any policy that does not lead to an unbounded cost
should keep the system critically loaded at all times,
and thus, asymptotically, the fractions of effort will vary
dynamically within $\slp$.
In the control literature, a control process that takes values only
at the vertices of the action space is called a bang-bang control.
The analogue of this notion in our setting is a policy for which
the limiting fractions of effort take values
only in $\{\xi^{*,1},\xi^{*,2}\}$, switching between the two
extremal solutions. Some of the policies introduced in this paper
are designed to act that way.

Contrary to the setting where the HTC holds,
it is impossible to construct an AO policy based only
on the first order data under the EHTC.
A second order approximation of the PSS, which is often referred to as
a {\it Brownian control problem} (BCP), is required.
The BCP represents a diffusion limit of the PSS,
in which Brownian motion (BM) replaces stochastic fluctuations
associated with cumulative arrival and service processes.
Closely related to the BCP is another diffusion control problem,
called a {\it workload control problem} (WCP).
Obtained by a certain projection of the BCP, it is a control problem
in which the process is one-dimensional, representing the total workload
asymptotics. The structure of the WCP obtained is quite simple to describe.
The state process is a reflected diffusion on $\R_+$, with
controlled drift and diffusion coefficients, $b=b(\xi)$, $\sig=\sig(\xi)$,
where the control process, $\xi=\xi_t$ takes values in $\slp$
and $\xi\mapsto(b(\xi),\sig(\xi)^2)$ is an affine map.
The cost is given as an expected discounted version of the
state process itself. By a standard argument based on the HJB equation, there
exists an optimal bang-bang control for the WCP.
There can therefore be two possibilities for the
WCP solution: the single mode case,
where one of the modes is always used, and the dual mode case,
where both modes are used by the optimal control in different parts
of the state space. Note that in this case the GBA can be changed dynamically.
The HJB equation also reveals the structure of the
feedback function from state to control. This particular
HJB equation was solved in \cite{Sheng78}. It was shown that
in the dual mode case there is a switching point
$z^*\in(0,\iy)$, such that one of the modes is used when the state is below
$z^*$ and the other otherwise. The HJB equation can be viewed,
in this case, as an equation involving a free boundary, in which the solution
is a pair, where one component is the value function and the other is
$z^*$. The results of
\cite{Sheng78} also characterize $z^*$ as the unique solution to
explicit equation, as well as a solution to the HJB equation.

Our policies are obtained by 'translating' the WCP solution.
In the case of a single mode,
the prescribed policy corresponds either to a threshold policy similar
to that of \cite{bw1} (see below)
or a simple priority policy, depending on the mode used and the cost.
In the case of dual mode, pairs of elementary policies are identified,
which are combined together to form switched control policies,
so that one is active when the normalized workload process
is below the switching point and the other
above it. In each case, the policy is designed to meet the target
allocation efforts determined by the corresponding
mode, and the set of operational activities is restricted by
the corresponding GBA.

The paper closest to ours is the aforementioned \cite{bw1},
that studies a 2-server, 2-class PSS with 3 activities.
This PSS is known as an 'N' network, because upon relabeling,
the activities are given by $(1,1)$, $(1,2)$ and $(2,2)$, forming
the symbol N. In this network the number of solutions to the LP
cannot exceed 1, and thus the requirement of a unique solution
does not pose a restriction. In an earlier work, \cite{har98},
it had been observed that when the larger '$c\mu$' value is in class 1,
the BCP solution suggests that the queue length of class 1 customers
and the idleness process at server 1 should both converge to zero
at the diffusion scale, and that a simple priority policy does not achieve this. 
In \cite{bw1} this was addressed by putting
a threshold on class 1 queue length, that when exceeded,
server 2 prioritizes class 1, and otherwise it prioritizes class 2.
The size of the threshold must converge to zero at the diffusion scale
so as to achieve the first goal. To achieve AO of a threshold policy
with logarithmic (in $n$) size threshold, as used in \cite{bw1}, the interarrival and service times are
assumed there to possess exponential moments.
(More on the history of the problem and the works that contributed
to its development can be read in \cite{ACR1}.)

As already mentioned, one of the policies we implement is a threshold
policy similar to the one used by \cite{bw1}.
However, our assumptions are positioned differently with respect to the threshold--moment
tradeoff, assuming a larger (still $o(\sqrt{n})$), polynomial size threshold,
but requiring only a polynomial moment assumption.
We assume $2+\eps$ moment assumptions
for all of our policies except the single-mode threshold policy
and the dual-mode policies that employ the threshold policy
when the workload is above $z^*$. For these, a finite $\m_0$-th moment
is assumed, where the number $4<\m_0<5$ is indicated explicitly.
Another difference between our results and those of \cite{bw1} is that our policies do not use preemption.
Although the policy introduced in \cite{bw1} uses preemption, it is plausible that an analogous non-preemptive policy is also AO under similar conditions. In this paper, our choice not to use preemption leads to non-trivial issues in the dual mode case. Instead of a simple switching between elementary policies when the workload crosses $z^*$, it is sometimes the case that one must wait for a particular server to become available before switching. This is described in \S \ref{sec:disc}.

It is also worth mentioning that we have argued
in \cite{ACR1} that the AO of the threshold policy
from \cite{bw1} extends
beyond the HTC to the case of multiple solutions and a single mode
(under some assumptions which
include the existence of exponential moments).

Beside the objective to break the uniqueness barrier,
an additional source of motivation for this work stems from the relation
between non-uniqueness and service rate decomposability.
As stated in Lemma \ref{lem1}, for the $2\times2$ PSS,
the LP exhibits multiple solutions if and only if the service rates decompose
as $\mu_{ij}=\al_i\beta_j$.
Service rates decompose this way when
the mean size of a job is characteristic to the class (and then
$\al_i$ is the reciprocal mean), and each server has its own processing speed
(here given by $\beta_j$).
As the HTC does not hold under decomposability,
this important class of service rates has been left out by earlier work.

\subsection{Results}

The description of the policies given above is only a sketch.
There are nontrivial issues
% The main reason
%is that, depending on the nature of the extreme point arising in the control, the prescribed policy corresponds either to the threshold policy similar to that of \cite{bw1} or a simple priority policy, and they may differ above and below $z^*$.
that arise
regarding the need to 'patch' 2 policy types, requiring us to slightly modify the policies, where the details differ from one pairing to another.

The main result states that, under the prescribed policies,
the rescaled workload process converges in law to the diffusion process
that solves the WCP, and these policies are AO.
As far as convergence is concerned,
in addition to the 'standard' issues involved in proving state-space collapse, we need to deal with issues related to switching control modes at $z^*$.
Moreover, to obtain AO from weak convergence, uniform integrability
needs to be established, and it is here where the $2+\eps$
and $\m_0$ moment assumptions are used.

An approach to proving convergence to a diffusion
with discontinuous coefficients, addressing especially the technicalities
involved with the discontinuity of the diffusion coefficient,
was developed in \cite{kry02}, going beyond
the general framework for convergence of semimartingales
such as that from \cite{lip-shi}.
Whereas the tools from \cite{kry02} are not directly applicable in our setting,
an argument which, as in \cite{kry02}, shows
that the time spent near the discontinuity set is negligible, is also at the basis of our proof.
The paper \cite{kry02} also gives an example of a queueing model
whose scaling limit  yields a diffusion process with discontinuities
in both drift and diffusion coefficients.
Our dual mode case provides what seems to be the first example
where this occurs under an AO policy of a queueing control problem
(for AO in heavy traffic leading to discontinuity in the drift only, see \cite{ata-lev}).

\subsection{Organization of the paper}

In \S 2.1 we describe our model 
and the control problem associated with it in more detail.
The LP and the extended heavy traffic condition are introduced in \S 2.2
and preliminary results about the LP from \cite{ACR1} are stated. In \S 2.3, the WCP and the associated
HJB equation are introduced, and in
Proposition \ref{prop0}, it is stated that there exists a unique
classical solution to the HJB equation.
This proposition also provides a condition which determines
whether an optimal solution to the WCP must employ a single mode or two modes
(not to be confused with the number of modes in the space
of LP solutions,
which is always two under multiplicity),
and asserts that in the dual mode case
there exists a single switching point $z^*$ in workload space. We also present in this section the lower bound from \cite{ACR1} stated in
Theorem \ref{thm:lowerbound}.
The main result is stated in \S \ref{sec:policy}.  The definitions of the proposed policies appear first,
and then, in Theorem \ref{th-ao-s}, the weak convergence and
AO results are stated.
 A discussion about the sampling times for switching policies 
 is contained in \S \ref{sec:disc}

In \S \ref{sec3} we state and prove some results related to the static allocation LP, providing, in particular, explicit expressions for the extreme points of the set of optimal solutions.
Development of the WCP is carried out in \S \ref{sec4}. Preliminary results proved in \cite{ACR1} in a general case are included in this section. 
This section also contains proofs of results related to the HJB equation, some of which rely on \cite{Sheng78}.

The proof of our main result, Theorem \ref{th-ao-s}, is the subject of
\S \ref{sec5}.
 In \S \ref{sec51}, we present the general scheme for proving Theorem \ref{th-ao-s}: the weak convergence result is stated in Theorem \ref{th3}. We then present four propositions that are used for the proof.
Each proposition corresponds to a specific section and step of the proof. 
Proposition \ref{lem:m1}, in   \S  \ref{sec:ui}, proves uniform integrability;   
state space collapse is proved in Proposition \ref{lem:ssc} in  \S\ref{sec:ssc}; 
a key non idling property is proved in Proposition \ref{lem:reflection} in \S  \ref{sec:boundarybeh}; 
and a  'fast switching'  property, showing the aforementioned property that the process spends asymptotically negligible time near the discontinuity, is proved in Propostion \ref{lem:correctmode} in \S \ref{sec:fs}.
Finally, the appendix contains proofs of several lemmas stated earlier.

\subsection{Notation}

$\N$, $\R$ and $\R_+$ are the sets of natural, real and, respectively,
nonnegative real numbers.
For $a, b \in \R$, $a \vee b$ and $a \wedge b$ denote the maximum and
minimum of $a$ and $b$, respectively, and
$a^+=a \vee 0$.
For a set $A$, $\one_A$ denotes its indicator function.
For $f:\R_+\to\R$ and $t,\del>0$,
$\|f\|_t=\sup_{s\in[0,t]}|f(s)|$ and
\[
w_t(f,\del)=\sup\{|f(s_1)-f(s_2)|:0\le s_1\le s_2\le(s_1+\del)\w t\}.
\]
For $0\le s\le t$, the notation $f[s,t]$ stands for $f(t)-f(s)$.
For real-valued functions and processes, the notation $X(t)$ is
used interchangeably with $X_t$.
Given a Polish space $E$, $C_E[0,\iy)$ and $D_E[0,\infty)$ denote
the spaces of $E$-valued, continuous and, respectively, c\`{a}dl\`{a}g functions on $[0,\infty)$,
equipped with the topology of convergence u.o.c.\ and, respectively,
the $J_1$ topology.
Denote by $C^+_\R[0,\iy)$ and $D^+_\R[0,\infty)$ the subset of
$C_\R[0,\iy)$ and, respectively, $D_\R[0,\iy)$,
of non-negative, non-decreasing functions, and by $C^{0,+}_\R[0,\iy)$
the subset of $C^+_\R[0,\iy)$ of functions that are null at zero.
Write $X_n\To X$ for convergence in law.
A sequence of processes with sample paths in $D_E[0,\iy)$ is said to be $C$-tight
if it is tight and the limit of every
weakly convergent subsequence has sample paths in $C_E[0,\iy)$ a.s.
The letter $c$ denotes a deterministic constant whose value may change from one
appearance to another.

\section{Model and main results}\label{sec2}
\beginsec

\subsection{Queueing model, scaling and queueing control problem}

The model under consideration is as in \cite{ACR1}, specialized to
the case of two job classes, two servers and four activities.
We will refer to is as the $2\times2$ PSS when there is need to distinguish it
from the general PSS treated in \cite{ACR1}.
The symbol $i\in\{1,2\}$ is used as a generic index to a class, and
$k\in\{1,2\}$ to a server. For a general PSS,
an {\it activity} is a class-server pair $(i,k)$ where
server $k$ is capable of serving class $i$. In this paper it is
assumed that each server is capable of serving each class,
hence there are four activities. They are labeled by $(i,k)$
or sometimes by $j=(i,k)$.

The model consists of a sequence of systems, indexed by $n\in\N$, that are
all defined on one probability space $(\Om,\calF,\PP)$.
For the $n$th system, one considers the following processes.
The processes denoted $A^n=(A^n_i)$  and $S^n=(S^n_{ik})$
represent arrival and potential service counting processes.
That is, $A^n_i(t)$ is the number of arrivals of class $i$ jobs until time $t$, $i=1,2$,
and $S^n_{ik}(t)$ is the number of service completions of class $i$ jobs
by server $k$, by the time server $k$ has devoted $t$ units of time to class $i$,
$i=1,2$, $k=1,2$. Next, $X^n=(X^n_i)$, $I^n=(I^n_k)$, $D^n=(D^n_{ik})$ and $T^n=(T^n_{ik})$
denote queue length, cumulative idleness, departure, and cumulative busyness
processes. In other words, $X^n_i(t)$ is the number of class $i$ customers in the system at time $t$,
$I^n_k(t)$ is the cumulative time server $k$ has been idle by time $t$, $D^n_{ik}(t)$
is the number of class $i$ departures from server $k$, and $T^n_{ik}(t)$
is the cumulative time devoted by server $k$ to class $i$.
The process $T^n_{ik}$ takes the form $T_{ik}(t)=\int_0^t\X^n_{ik}(s)ds$,
where $\X^n_{ik}(t)$ is the fraction of effort devoted by server $k$ to class-$i$ jobs at $t$. In particular, $\sum_i\X^n_{ik}(t)\le1$ for every $k$.
Thus $\X^n$ is referred to as the allocation process.

The aforementioned arrival and potential service processes
are constructed as follows. Arrival rates $\la^n_i$ and service rates
$\mu^n_{ik}$ are given, satisfying, for some constants
$\la_i\in(0,\iy)$, $\mu_{ik}\in(0,\iy)$, $\hat\la_i\in\R$,
$\hat\mu_{ik}\in\R$,
\begin{align*}
	&
	\hat\la^n_i:=n^{-1/2}(\la^n_i-n\la_i)\to\hat\la_i,\\
	&
	\hat\mu^n_{ik}:=n^{-1/2}(\mu^n_{ik}-n\mu_{ik})\to\hat\mu_{ik},
\end{align*}
as $n\to\iy$.
For each $i$ a renewal process $\check A_i$ is given,
with interarrival distribution that has mean 1 and squared
coefficient of variation $0<C^2_{A_i}<\iy$.
Similarly, for each $(i,k)$, a renewal process $\check S_{ik}$
is given with mean 1 interarrival and squared coefficient of
variation $0<C^2_{S_{ik}}<\iy$. It is assumed that
$A^n$ and $S^n$ are given by
\[
A^n_i(t)=\check A_i(\la^n_it),
\qquad
S^n_{ik}(t)=\check S_{ik}(\mu^n_{ik}t).
\]
It is assumed moreover that the six processes $\check A_i$, $\check S_{ik}$
are mutually independent, have strictly positive inter-arrival distributions
and right-continuous sample paths.
The (IID) interarrivals of $\check A_i$ and $\check S_{ik}$ are denoted by
$\check{a}_{i}(l)$ and $\check{u}_{ik}(l)$, $l\in\N$, respectively, and those of the accelerated
processes $A^n_i$ and $S^n_{ik}$ are given by
\begin{equation}\label{e01}
	a^n_{i}(l)=\dfrac{1}{\lambda^n_{i}}\check{a}_{i}(l),
	\qquad
	u^n_{ik}(l)=\dfrac{1}{\mu^n_{ik}}\check{u}_{ik}(l ).
\end{equation}
The system is assumed to start empty, that is, $X^n(0)=0$ for all $n$.
Simple relations between the processes are
\begin{equation}\label{40-}
	D^n_{ik}(t)=S^n_{ik}(T^n_{ik}(t)),
\end{equation}
\begin{equation}\label{40}
	X^n_i(t)=A^n_i(t)-\sum_kD^n_{ik}(t),
\end{equation}
\begin{equation}\label{41}
	I^n_k(t)=t-\sum_iT^n_{ik}(t),
\end{equation}
\begin{equation}\label{41+}
	\text{the sample paths of $X^n_i$ are nonnegative, and
		those of $I^n_k$ are in $C_\R^{0,+}[0,\iy)$.}
\end{equation}

The tuple $(\check A,\check S)$ is referred to as the {\it stochastic primitives}.
In our formulation we will consider $T^n$ as the control process
(equivalently, the allocation process $\X^n$ may be regarded the control).
In view of equations \eqref{40-}, \eqref{40}, \eqref{41},
given the stochastic primitives, the control uniquely determines the processes $D^n$, $X^n$, $I^n$.
Let an additional process be defined on the probability space
denoted by $\Ups=(\Ups(l),l\in\N)$,
taking values in a Polish space $\calS_\text{rand}$
and assumed to be independent of the stochastic primitives, for each $n$
(there is no need to let $\Ups$ vary with $n$, as the primitives
are all defined on the same probability space).
It is included in the model in order to allow the construction of randomized controls; for more details about its potential use see \cite[Remark 2.1.ii]{ACR1}.

The process $T^n$ is said to be an {\it admissible control for the queueing control problem
	(QCP) for the $n$-th system}
if for each $(i,k)$, $T^n_{ik}$ has sample paths in $C_\R^{0,+}[0,\iy)$
that are $1$-Lipschitz, and the associated processes $D^n$,
$X^n$ and $I^n$ given by \eqref{40-}, \eqref{40} and \eqref{41} satisfy \eqref{41+};
furthermore, $T^n$ is adapted to the filtration $\{\calF^n_t\}$
defined by $\calF^n_t=\s\{(A^n(s),D^n(s),s\in[0,t]), \Ups\}$.
Denote by $\calA^n$ the collection of all admissible controls for the QCP
for the $n$-th system. As argued in
\cite[Remark 2.1.i]{ACR1}, this definition allows for the control to
depend on the history of all processes involved in the model
(in addition to the auxiliary randomness $\Ups$).

The queue length process normalized at the diffusion scale is
defined by $\hat X^n_i(t)=n^{-1/2}X^n_i(t)$.
The cost of interest for the $n$-th system is given by
\begin{equation}\label{42}
	\hat J^n(T^n)=\E\int_0^\iy e^{-\gamma t}h(\hat X^n(t))dt,
	\qquad
	T^n\in\calA^n,
\end{equation}
where $\gamma>0$ and $h(x)=h_1x_1+h_2x_2$, with constants $h_1,h_2>0$,
and $\hat X^n$ is the rescaled queue length process associated with
the admissible control $T^n$.
The value for the $n$-th system is defined by
\[
\hat V^n=\inf\{\hat J^n(T^n):T^n\in\calA^n\}.
\]
This completes the description of the queueing models and QCP.
The complete set of problem data consists of the stochastic primitives
mentioned above and the collection of parameters
\[
(\la_i),(\mu_{ik}),(\hat\la_i),(\hat\mu_{ik}),(C_i),(C_{ik}),\gamma,(h_i).
\]
We sometimes refer to $(\la_i), (\mu_{ik})$ as the first order data
and to $(\hat\la_i),(\hat\mu_{ik}), (C_i), (C_{ik})$ as the second order data.

\subsection{The linear program and extended heavy traffic condition}\label{sec21}

Given the first order data $\la_i>0$,
$\mu_{ik}>0$, $i,k=1,2$,
consider the following linear program (LP) for the unknowns
$(\xi_{ik})\in\R_+^{2\times 2}$ and $\rho\in\R$.

\noi{\it Linear Program.}
Minimize $\rho$ subject to
\begin{equation}\label{01}
	\begin{cases}
		\ds
		\sum_{k=1}^2\xi_{ik}\mu_{ik}=\la_i & i=1,2,
		\\ \\
		\ds
		\sum_{i=1}^2\xi_{ik}\les\rho & k=1,2,
		\\ \\
		\xi_{ik}\geqslant0 & i,k=1,2.
	\end{cases}
\end{equation}
Denote the optimal objective value of \eqref{01} by $\rho^*$.

\noi{\it Extended heavy traffic condition.}
$\rho^*=1$.

The extended heavy traffic condition (EHTC)
is broader than the {\it heavy traffic condition}
that has been extensively used in the literature,
which requires, in addition to $\rho^*=1$, that there be
a unique corresponding $\xi$.

Under the EHTC, any solution is of the form $(\xi,1)$.
Let $\slp$ denote the subset of $\R^{2\times 2}$ for which
the set of all solutions is given by $\slp\times\{1\}$.
We say that the {\it EHTC with multiplicity}
(EHTCM) holds if the EHTC holds and the LP has multiple solutions
(that is, there exist two distinct pairs $(\xi^{(1)},1)$
and $(\xi^{(2)},1)$ satisfying \eqref{01}).
We say that the service rates $\mu_{ik}$ are of {\it product form}
if $\mu_{ik}=\al_i\beta_k$ for all $i,k$, for some constants
$\al_i$ and $\beta_k$.

A matrix $\xi\in\R_+^{2\times2}$ is called {\it column-stochastic} if
$\sum_i\xi_{ik}=1$ for both $k=1,2$. A column-stochastic matrix is called
a {\it mode} if (at least) one of its columns is either $(0,1)^T$ or $(1,0)^T$.
A mode is said to be {\it degenerate} if it has more than one zero
entry; otherwise it is said to be {\it nondegenerate}.
A pair of nondegenerate modes is said to be a {\it class-switched}
({\it server-switched}) pair of modes if the zero entries in
the two modes are in distinct rows but the same column
(respectively, distinct columns but the same row).
The following condition will be referred to as the {\it nondegeneracy condition},
namely
\begin{equation}\label{a2}
	\la_i\ne\mu_{ik} \text{ for all } (i,k)\in\{1,2\}^2.
\end{equation}
The following is proved in \S \ref{sec3}.

\begin{lemma}\label{lem1}
	Let the EHTC hold.
	\begin{enumerate}
		\item For any solution $(\xi,1)$, $\xi$ is column-stochastic.
		\item  The LP \eqref{01} has multiple solutions
		if and only if $(\mu_{ik})$ are of product form.
		\item If the LP has multiple solutions then there exists a pair
		of modes $(\xi^{*,1}, \xi^{*,2})$ such that
		\begin{equation}\label{a1}
			\slp={\rm ch}(\{\xi^{*,1},\xi^{*,2}\}).
		\end{equation}
		\item
		If the LP has multiple solutions and the nondegeneracy condition
		\eqref{a2} holds
		then both $\xi^{*,1}$ and $\xi^{*,2}$ of \eqref{a1} are nondegenerate. Moreover,
		they form either a class-switched or a server-switched pair.
	\end{enumerate}
\end{lemma}

The main result will be proved under the following.
\begin{assumption}\label{assn1}
	\begin{enumerate}
		\item
		The EHTCM holds.
		\item
		The nondegeneracy condition \eqref{a2} holds.
	\end{enumerate}
\end{assumption}
Note that the case where the EHTC holds but EHTCM does not hold is already
covered in the work \cite{bw1}.

In view of Lemma \ref{lem1}, under Assumption \ref{assn1},
the rates $(\mu_{ik})$ are of product form. Thus $\mu_{ik}=\al_i\beta_k$,
and clearly there is a degree of freedom in choosing $(\al_i)$ and $(\beta_k)$.
In this paper we will always assume that they are chosen so that
$\sum_k \beta_k=1$, and it is easy to see that, given $(\mu_{ik})$, this normalization
uniquely determines these parameters.

It is also guaranteed by the lemma that, under Assumption \ref{assn1},
the extreme points of $\slp$ are two nondegenerate modes $\xi^{*,1}, \xi^{*,2}$
forming a class- or a server-switched pair.
Once a labeling of these modes has been fixed,
we will sometimes slightly abuse the terminology by referring to them as modes $1$ and $2$
rather than modes $\xi^{*,1}$ and $\xi^{*,2}$.

In earlier work on PSS, under the assumption that the LP has a unique solution $(\xi^*,1)$, activities are categorized as
{\it basic} or {\it nonbasic} according to the positivity of the fraction
allocated to them by $\xi^*$, that is, an activity
$(i,k)$ is {\it basic} if $\xi^*_{ik}>0$ and {\it nonbasic} if $\xi^*_{ik}=0$.
We extend this terminology to the case of multiple solutions
as follows. For $m\in\{1,2\}$,
an activity $(i,k)$ is said to be {\it basic in mode $m$} if the allocation associated to it by this mode
does not vanish, namely $\xi^{*,m}_{ik}>0$. If $\xi^{*,m}_{ik}=0$ (respectively, $\xi^{*,m}_{ik}=1$)
it is said to be {\it non-basic (respectively, full) in mode $m$}.

\begin{figure}[h]
	\centering
	\begin{tikzpicture}[scale=1]
		%% nodes left
		\node[shape=circle,draw=black] (A) at (-1,1) {1};
		\node[shape=circle,draw=black] (B) at (1,1) {2};
		\node[shape=circle,draw=black] (C) at (1,-1) {2};
		\node[shape=circle,draw=black] (D) at (-1,-1) {1};
		%%edges left;
		\path [thick](B) edge node[right] {} (C);
		\path [dashed] (B) edge node[left] {} (D);
		\path [line width=2.5pt](D) edge node[left] {\quad\quad} (A);
		\path[thick](A)edge node[left]{}(C);
		%%nodes right
		%	\node[shape=circle,draw=black] (E) at (3,1) {1};
		%	\node[shape=circle,draw=black] (F) at (5,1) {2};
		%	\node[shape=circle,draw=black] (G) at (5,-1) {2};
		%	\node[shape=circle,draw=black] (H) at (3,-1) {1};
		%%edges right
		%	\path [thick](E) edge node[right] {$\dfrac{\lambda_1}{\alpha_1\beta_1}$} (H);
		%	\path [dashed] (E) edge node[left] {} (G);
		%	\path [color=red, line width=2.5pt](F) edge node[left] {$\xi^{*,2}_{22}=1$} (G);
		%	\path[thick](F)edge node[left]{}(H);
		
		\node[shape=circle,draw=black] (E) at(7,1) {1};
		\node[shape=circle,draw=black] (F) at (9,1) {2};
		\node[shape=circle,draw=black] (G) at (9,-1) {2};
		\node[shape=circle,draw=black] (H) at (7,-1) {1};
		%%edges left;
		\path  [dashed](F) edge node[right] {} (G);
		\path [thick] (F) edge node[left] {} (H);
		\path [thick](H) edge node[left] {\quad\quad} (E);
		\path[line width=2.5pt](E)edge node[left]{}(G);
		
		\node[shape=circle,draw=black] (I) at (3,1) {1};
		\node[shape=circle,draw=black] (J) at (5,1) {2};
		\node[shape=circle,draw=black] (K) at (5,-1) {2};
		\node[shape=circle,draw=black] (L) at (3,-1) {1};
		%%edges left;
		\path [thick](J) edge node[right] {} (K);
		\path  [line width=2.5pt](J) edge node[left] {} (L);
		\path [dashed](L) edge node[left] {\quad\quad} (I);
		\path[thick](I)edge node[left]{}(K);
		
	\end{tikzpicture}
	\\
	\flushleft
	\hspace{10.5em} (a) \hspace{8.5em} (b) \hspace{8.5em} (c)
	\caption{\label{fig1}\sl
		The full and non-basic activities are shown in
		thick and dashed lines, respectively.
		Graph (a) corresponds to a mode in canonical form.
		Graphs (a,b) correspond to a pair of modes where the non-basic activity
		switches a class, whereas in (a,c) it switches a server.
		The pair (b,c), in which the non-basic activity switches both a class
		and a server is neither class- nor server-switched.
	}
\end{figure}
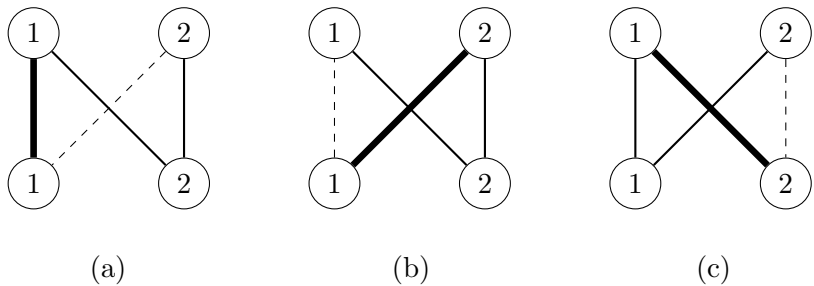

A mode is said to be in {\it canonical form}
if its first column is $(1,0)^T$. It is clear by the definition
of a mode that it is always
possible to relabel the classes and the servers so that a given mode
is in canonical form, and that if the mode is nondegenerate there is only one
such relabeling.
The graph of a mode in canonical form is shown in Figure \ref{fig1}(a).
Because of its resemblance to the symbol {\it N}, this form is sometimes called an {\it N-system.}

The graphs in Figure \ref{fig1}(a) and (b) correspond to a class-switched pair of modes,
whereas those in Figure \ref{fig1}(a) and (c) correspond to a server-switched pair.
Hence it is seen that when switching between (nondegenerate)
modes, the non-basic activity may change either a class or a server, but not both;
under a class-switched pair, the non-basic
activity switches a class, whereas under a server-switched pair it switches a server.
The terms class-switched and server-switched will sometimes be abbreviated as
\textbf{CS} and \textbf{SS}.

If the EHTCM holds but \eqref{a2} does not, that is,
there exist $i,k$ such that $\lambda_i=\mu_{ik}$, the situation is different: at least one of the modes will be degenerate, i.e., have two non-basic activities (see Lemma \ref{lem:solvelp} below). In the terminology of linear programming
this corresponds to a case where one of the basic solutions of the LP \eqref{01}
is degenerate (cf.\ \cite[Definition 3.1]{gol-tod}). The degenerate case is not
covered in this paper.

\subsection{Workload control problem}\label{sec:wcpp}

The WCP was derived and studied in \cite{ACR1} under the EHTC.
We describe this problem in the special case needed here, namely
under the setting of a $2\times2$ PSS and assuming that the EHTCM holds.
In particular, as mentioned above, the parameters $(\al_i)$, $(\beta_k)$ are
uniquely determined by the problem data.
Define the workload process and its scaled version as
\begin{equation}\label{241}
	W^n(t)=\sum_{i}\dfrac{X^n_i(t)}{\alpha_i},
	\qquad
	\hat W^n(t)=\sum_{i}\dfrac{\hat{X}^n_i(t)}{\alpha_i}.
\end{equation}
(That this definition agrees with that of \cite{ACR1} follows from \cite[Lemma 2.4.2]{ACR1}).
Let the process that appears in the definition of the cost \eqref{42} be denoted by
\begin{equation}\label{143}
	\hat H^n_t=h(\hat X^n(t))=h_1\hat X^n_1(t)+h_2\hat X^n_2(t).
\end{equation}
Throughout, denote by $p,q\in\{1,2\}$ the two distinct indices for which
\begin{equation}\label{145}
	h_p\al_p\ge h_q\al_q,
\end{equation}
where in the special case $h_1\al_1=h_2\al_2$, set $p=1$ and $q=2$.
The policies constructed in this paper aim at keeping $\hat X^n_p$ close to zero.
Hence we call $p$ the {\it high priority class} (HPC) and $q$
the {\it low priority class} (LPC).
Next, let $\sig_{A,i}=\la_i^{1/2}C_{A_i}$, $\sig_{S,ik}=\mu_{ik}^{1/2}C_{S_{ik}}$, and 
\begin{equation}\label{50}
	b(\xi)=\sum_i\frac{\hat\la_i-\sum_k\hat\mu_{ik}\xi_{ik}}{\al_i},
	\qquad
	\sig(\xi)^2=
	\sum_i\frac{\sig_{A,i}^2+\sum_k\sig_{S,ik}^2\xi_{ik}}{\al_i^2},
	\qquad
	\xi=(\xi_{ik})\in\slp.
\end{equation}
Let also
\begin{equation}
	b_m=b(\xi^{*,m}),\qquad \sig_m=\sig(\xi^{*,m}),\qquad m=1,2.\label{eq:param}
\end{equation}
Consider a one-dimensional controlled diffusion
with controlled drift and diffusion coefficients and reflection at the origin, given by
\begin{equation}\label{15}
	Z_t=z+\int_0^tb(\X_s)ds+\int_0^t\sig(\X_s)dB_s+L_t,
\end{equation}
where $\X$ is a control process, $B$ is a standard BM (SBM), $L$ is a reflection term
at zero, and $z\geqslant0$. A precise definition is as follows.

Given a filtration $(\calF'_t)$, let $\frX((\calF'_t))$ denote
the collection of $(\calF'_t)$-progressively measurable
processes $\{\X_t\}$ taking values in $\R^{2\times 2}$,
for which $\PP(\text{for a.e.\ $t$, } \X_t\in\slp)=1$.
A tuple $\frS=(\Om',\calF',(\calF'_t),\PP', B,\X,Z,L)$ is said to be
an {\it admissible control system for the WCP with initial condition $z$}
if $(\Om',\calF',(\calF'_t),\PP')$ is a filtered probability space,
$B$, $\X$, $Z$ and $L$ are processes defined on it,
$B$ is a SBM and an $(\calF'_t)$-martingale,
$\X\in\frX((\calF'_t))$, $Z$ is continuous nonnegative and $(\calF'_t)$-adapted,
$L$ has sample paths in $C_\R^{0,+}[0,\iy)$ and is $(\calF'_t)$-adapted,
and equation \eqref{15} and
the identity $\int_{[0,\iy)}Z_tdL_t=0$ are satisfied $\PP'$-a.s.

Denoting by $\calA_\text{WCP}(z)$ the collection of all control systems for the WCP
with initial condition $z$, the cost and value are defined as
\begin{equation}\label{51}
	J_\text{WCP}(z,\frS)=\E_{\frS}\int_0^\iy e^{-\gamma t}Z_tdt,
	\qquad
	V_\text{WCP}(z)=\inf\{J_\text{WCP}(z,\frS):\frS\in\calA_\text{WCP}(z)\}.
\end{equation}
The main result of \cite{ACR1}, when specialized to the $2\times2$ PSS,
states the following.

\begin{theorem}[\cite{ACR1}]
	\label{thm:lowerbound}
	Let Assumption \ref{assn1} (multiplicity and nondegeneracy) hold. Then
	\[
	\liminf_n\hat V^n\ges V_0:=h_q\alpha_qV_{\rm WCP}(0).
	\]
\end{theorem}
To prove this result one only needs to verify that the assumptions
from \cite{ACR1} hold under Assumption \ref{assn1} of this paper.
This is done in \S\ref{sec3}.
\begin{remark}
	The result from \cite{ACR1} is not limited to the EHTCM, only to the EHTC.
	However, without this condition, the parameters $(\al_i)$, $(\beta_k)$
	are not defined, and the WCP requires the use of the dual to the LP \eqref{01},
	which in this paper need not be introduced.
\end{remark}

In view of Theorem \ref{thm:lowerbound},
a sequence $T^n\in\calA^n$ of admissible controls for the QCP,
also referred to as a {\it sequence of policies}, is said to be
{\it asymptotically optimal} (AO)
if $\limsup_{n\to\iy}\hat J^n(T^n)\les V_0$.

The value function can be characterized in terms of an HJB equation.
To present this equation, for $(v_1,v_2,\xi)\in\R^2\times\slp$, let
\[
\bar\BH(v_1,v_2,\xi)=b(\xi)v_1+\frac{\sig(\xi)^2}{2}v_2,
\qquad
\BH(v_1,v_2)=\inf_{\xi\in\slp}\bar\BH(v_1,v_2,\xi)
=\min_{\xi\in\{\xi^{*,1},\xi^{*,2}\}}\bar\BH(v_1,v_2,\xi),
\]
where the identity on the RHS follows from the fact that both $b$ and $\sig^2$
are affine as a function of $\xi$, and $\slp$ is the convex hull of $\{\xi^{*,1},\xi^{*,2}\}$.
A classical solution to the HJB equation is a $C^2(\R_+:\R)$ function $u$ satisfying
\begin{equation}
	\label{14}
	\BH(u'(z),u''(z))+z-\gamma u(z)=0,
	\qquad z\in(0,\iy),
\end{equation}
and the boundary conditions at $0$ and $\iy$,
\begin{equation}\label{14+}
	u'(0)=0,
	\qquad
	u(z)<c(1+z), z\in\R_+, \text{ for some constant $c$}.
\end{equation}
Given a $C^2$ function $u$, denote
$\BH^u_m(z)=\bar\BH(u'(z),u''(z),\xi^{*,m})$, $m=1,2$,
and $\BH^u(z)=\min_{m}\BH^u_m(z)$.

The following two conditions, which play an important role in what follows,
are complementary:
\begin{equation}\label{90}
	\text{there exist distinct $m,m'\in\{1,2\}$ such that
		$b_m\les b_{m'}$ and $\sig_m\les \sig_{m'}$,}
\end{equation}
\begin{equation}\label{91}
	\text{there exist distinct $m,m'\in\{1,2\}$ such that
		$b_m<b_{m'}$ and $\sig_m>\sig_{m'}$.}
\end{equation}
The $C^2$ smoothness of the value function is tied to the question of existence
of a classical solution to the HJB equation.
Owing to the uniform ellipticity ($\sig(\xi)^2>0$ at both $\xi=\xi^{*,1}$ and
$\xi^{*,2}$), one can show that a classical solution uniquely exists \cite[Proposition 2.5]{ACR1}, a type of result that, in the context of optimal switching of a diffusion process,
has been called {\it the principle of smooth fit}.
In this paper we rely, in addition, on the results of \cite{Sheng78} which address
the specific control problem treated here,
and provide further structural properties (parts 2 and 3 of the result below)
that are based on the construction of a solution to the HJB equation
and are harder to obtain via the general approach.

\begin{proposition}\label{prop0}
	\begin{enumerate}
		\item There exists a unique classical solution $u$ to \eqref{14}--\eqref{14+}.
		Moreover, $u=V_\text{\rm WCP}$.
		\item  If \eqref{90} holds then, with $m$ as in \eqref{90},
		\begin{equation}\label{58-}
			\BH^u(z)=\BH^u_m(z) \quad z\in\R_+.
		\end{equation}
		\item Alternatively, if \eqref{91} holds then there exists $z^*\in(0,\iy)$
		such that, with the pair $(m,m')$ of \eqref{91},
		\begin{equation}\label{58}
			\BH^u(z)=
			\begin{cases}
				\BH^u_{m'}(z) & z<z^*,\\
				\BH^u_{m}(z) & z>z^*.
			\end{cases}
		\end{equation}
	\end{enumerate}
\end{proposition}

The proof of this result appears in \S\ref{sec4}.
Some details on the construction from \cite{Sheng78} appear in Appendix \ref{app:a}.

In view of this result, the case \eqref{90} will be referred to as the
{\it single mode} case, and \eqref{91} as the {\it dual mode} case
(not to be confused with uniqueness and multiple solutions to the LP).
In the former, the mode $\xi^{*,m}$ for  which $m$ satisfies \eqref{90} will be referred to
as the {\it active mode} and denoted $\xi^A$,
because the above result indicates that it is optimal to always
select $\X_t=\xi^{*,m}$.
In the latter, the modes $\xi^{*,m'}$ and $\xi^{*,m}$ for which $m'$ and $m$ satisfy
\eqref{91} will be referred to as $\xi^L$, the {\it lower} and, respectively,
$\xi^H$ the {\it higher workload mode}, and $z^*$ as the {\it switching point}.
These terms refer to the fact that the result suggests that it is optimal to select $\X_t=\xi^L$
(respectively, $\X_t=\xi^H$) when $Z_t< z^*$ (respectively, $Z_t>z^*$).

Going back to \eqref{15}, selecting a control process
according to the above description results in two different diffusion processes.
In the single mode case, the optimally controlled process $Z$ is given by
\begin{equation}\label{142-}
	Z^{(1)}_t=z+b^At+\sig^A B_t+L^{(1)}_t,
\end{equation}
with $b^A=b(\xi^A)$, $\sig^A=\sig(\xi^A)$ and $L^{(1)}$ a reflection term at $0$,
which is nothing but a reflecting BM with drift
$b^A$ and diffusivity $\sig^A$.
In the dual mode case, one is led to consider the SDE
\begin{equation}\label{142}
	Z^{(2)}_t=z+\int_0^tb^*(Z^{(2)}_s)ds+\int_0^t\sig^*(Z^{(2)}_s)dB_s+L^{(2)}_t,
\end{equation}
where, throughout, we denote
\begin{equation}\label{142+}
	b^*=b\circ\ph^*,\qquad \sig^*=\sig\circ\ph^*,
	\qquad
	\ph^*(z)=\xi^L\one_{[0,z^*]}(z)+\xi^H\one_{(z^*,\iy)}(z), \qquad z\in\R_+,
\end{equation}
and again $L^{(2)}$ is a reflection term.
For this equation, weak existence and uniqueness of solutions hold,
as we shall argue in Lemma \ref{lem3}.
As a result, there exists a control system for the WCP that behaves
exactly as described above, with
$\X_t=\xi^L$ (respectively, $\X^H$) when $Z^{(2)}_t\le z^*$ ($>z^*$),
and moreover, this description uniquely determines the law of the process $Z^{(2)}$.

As for the asymptotics of the QCP, the preceding discussion,
and the fact that the system starts empty,
suggest that in order to achieve the lower bound, the convergence
\begin{equation}\label{146}
	(\hat X^n_p,\hat W^n)\To(0,Z^{(1)}),
	\qquad \text{with } z=0,
\end{equation}
should hold in the single mode case, and
\begin{equation}\label{147}
	(\hat X^n_p,\hat W^n)\To(0,Z^{(2)}),
	\qquad \text{with } z=0,
\end{equation}
in the dual mode case, where $Z^{(1)}$ is given by \eqref{142-},
$(Z^{(2)},\X,L,B)$ is a weak solution to \eqref{142},
and in both cases the initial condition is $z=0$.

\subsection{Asymptotic optimality results}\label{sec:policy}

This section is devoted to the description of several policies
that are shown to be AO under different conditions.
We have already assumed that the interarrival times of the primitive processes
possess finite second moments. Our main results require a stronger assumption.

\begin{assumption}\label{as:polymo}
	There exists $\m>2$ such that
	\[
	\max_{i,k}\mathbb{E}[\check{a}_{i}(1)^\m]\vee\mathbb{E}[\check{u}_{ik}(1)^\m] <\infty.
	\]
\end{assumption} 
Whereas the assumption $\mathbf{m}>2$ is required for all our results,
some of them will require yet a stronger moment assumption, namely
$\mathbf{m}>\mathbf{m}_0$, where, throughout, we denote
$$\mathbf{m}_0=\frac{1}{2}(5+\sqrt{17}).$$

As already mentioned, under Assumption \ref{assn1},
both modes have a single non-basic activity.
Then, given a mode, the graph of basic activities has exactly three edges,
and one can speak of the single-activity class (the one associated
with only one nonbasic activity), the dual-activity class, and similarly,
the single- and dual-activity server. These terms allow us to refer to the roles
of classes and servers in the graph without considering a particular labeling.
It is also useful to accompany these terms with matching notation.
For a mode $\xi$, let the single- (respectively,
dual-) activity class be denoted by $i_1(\xi)$ (respectively, $i_2(\xi)$), and similarly,
let the single- (respectively, dual-) activity server be denoted by
$k_1(\xi)$ (respectively, $k_2(\xi)$).

Different policies are proposed in different cases.
The distinction between the various cases is based on whether the single-mode condition \eqref{90} or the dual-mode condition \eqref{91}
holds, and further, for each of the relevant modes
($\xi^A$ in the former case and both $\xi^L$ and $\xi^H$ in the latter),
whether the HPC is the single- or dual-activity class.

As a rule, all policies we describe are non-preemptive, that is,
the processing of a job is not interrupted once started.
A job is said to be {\it in the queue} if it is waiting to be served,
whereas it is {\it in the system} if it is either in the queue or being processed.
(In what is a bit of an abuse of terminology we use the term {\it queue length} to refer to the number in the system.)
A server is said to be {\it available} at a time $t$ if
either it has just completed a job or has already been idle at that time.

Some of the policies to be described are defined in terms of a sequence of thresholds,
$\Th^n$, put on the queue length at one of the two buffers.
Under Assumption \ref{as:polymo},
$\m>2$. Fix $\bar a$ satisfying
\begin{equation}\label{140}
	\frac{1}{2}-\bar\zeta(\m)<\bar a<\frac{1}{2}
	\qquad \text{where} \qquad
	\bar\zeta(\m)=\begin{cases}
		\frac{\m-2}{4\m}, & \m\in(2,\m_0],\\
		\frac{\m-2}{4\m}\w\frac{\m^2-5\m+2}{2\m(3\m-2)}
		=\frac{\m^2-5\m+2}{2\m(3\m-2)}, & \m\in(\m_0,\iy).
	\end{cases}
\end{equation}
Set the sequence of threshold levels $\Th^n$ and their normalized version
$\hat \Th^n$ to
\begin{equation}\label{141}
	\Th^n=\lceil n^{\bar a} \rceil, \qquad \hat \Th^n=n^{-1/2}\Th^n.
\end{equation}

\begin{definition}\label{def:0}(Server dedicated to / prioritizes a class).
	\begin{enumerate}
		\item
		A server is said to be
		{\em dedicated} to class $i$ at a given time if it acts as follows:
		if available at that time, it admits a job from class $i$
		provided there is one in the queue,
		or there is a new class-$i$ arrival, but does not admit a job from the other class.
		\item
		A server is said to {\em prioritize} class $i$ at a given time
		if it acts as follows:
		if available at that time,
		it admits a job from class $i$ provided there is one in the queue;
		otherwise it admits a job from the other class provided there is one in the queue.
		If the server is idle at that time and there is a new arrival of any class,
		it admits this arrival unless the other server is dedicated to that class and is free at that time
		(in which case the other server admits it).
	\end{enumerate}
\end{definition}

\begin{definition}\label{def:1}(\p{}, \teeone{} and \teetwo{} rules).
	Let a mode $\xi$ be given. At any moment in time the single-activity server is dedicated to the dual-activity class.
	\begin{enumerate}
		\item
		The servers are said to obey the {\em priority rule},
		abbreviated \p{} rule, at a given time if
		the dual-activity server prioritizes the single-activity class at that time.
		\item
		The servers are said to obey the {\em single-activity class threshold rule},
		abbreviated \teeone{} rule, at a given time if
		in the $n$-th system, the dual-activity server prioritizes the single-activity class
		when the queue length of the single-activity class equals or exceeds $\Th^n$ at that time,
		and otherwise  prioritizes
		the dual-activity class.
		\item
		The servers are said to obey the {\em dual-activity class threshold rule},
		abbreviated \teetwo{} rule, at a given time if
		in the $n$-th system, the dual-activity server prioritizes the dual-activity class
		when the queue length of the dual-activity class equals or exceeds $\Th^n$ at that time,
		and otherwise  prioritizes
		the single-activity class.
	\end{enumerate}
\end{definition}

Recall the notation $\xi^A$, $\xi^L$, $\xi^H$ and $p$ from \S\ref{sec:wcpp}.
In the case of a single mode, the following two policies are proposed.
(In Definitions \ref{def:2} and \ref{def:3} below, the text in square brackets is not
a part of the definition, but serves to indicate when each policy is to be applied).

\begin{definition}\label{def:2} (Single mode policies \p{} and \teetwo{}).
	\label{def:PP}\label{def:thresh}
	Let the single mode condition \eqref{90} hold.
	\begin{enumerate}
		\item
		The \p{} policy [to be applied when $i_1(\xi^A)=p$]
		is as follows: The servers obey the \p{} rule corresponding to $\xi^A$ at all times.
		\item
		The \teetwo{} policy [to be applied when $i_2(\xi^A)=p$] is as follows: The servers obey
		the \teetwo{} rule corresponding to $\xi^A$ at all times.
	\end{enumerate}
\end{definition}

Whereas in the case of a single mode the servers obey one rule at all times,
in the dual mode case they switch between two rules, depending, roughly speaking,
on whether the workload process is below or above the switching level $n^{1/2}z^*$.
The precise definition requires the use of a state variable called {\it current mode},
that determines which rule is applicable. This variable is not updated continuously
in time about whether $W^n>n^{1/2}z^*$
but only at certain sampling times, that differ from one case to another,
as detailed below.
The four policies used in the case of a dual mode are as follows.

\begin{definition} (Dual mode policies \p{}\p{}, \teetwo{}\teetwo{}, \teeone{}\teetwo{}
	and \teetwo{}\teeone{}).
	\label{def:3}
	\label{def:switchpolpppp}
	Let the dual-mode condition \eqref{91} hold.
	\begin{enumerate}
		\item The \p{}\p{} policy [applied when $i_1(\xi^L)=i_1(\xi^H)=p$] is as follows.
		\begin{enumerate}
			\item
			The workload is sampled at each service completion of the single activity server.
			If the workload is below $n^{1/2}z^*$, the current mode is set to $\xi=\xi^L$; otherwise
			it is set to $\xi=\xi^H$.
			\item
			The servers always obey the \p{} rule w.r.t.\ the current mode $\xi$.
			
		\end{enumerate}
		\item
		The \teetwo{}\teetwo{} policy [applied when $i_2(\xi^L)=i_2(\xi^H)=p$] is as follows.
		\label{def:switchpoltptp}
		\begin{enumerate}
			\item
			The workload is sampled at each service completion of the single activity server.
			If the workload is below $n^{1/2}z^*$, the current mode is set to $\xi=\xi^L$; otherwise
			it is set to $\xi=\xi^H$.
			\item
			The servers always obey the \teetwo{} rule w.r.t.\ the current mode $\xi$.
		\end{enumerate}
		\item
		The \teeone{}\teetwo{} policy [applied when $i_1(\xi^L)=i_2(\xi^H)=p$] is as follows.
		\label{def:switchpol}
		\begin{enumerate}
			\item
			The workload is sampled at each arrival and service completion.
			If the workload is below $n^{1/2}z^*$, the current mode is set to $\xi=\xi^L$; otherwise
			it is set to $\xi=\xi^H$.
			\item
			Whenever $\xi=\xi^L$, the servers obey the \teeone{} rule w.r.t.\ $\xi$;
			whenever $\xi=\xi^H$, they obey the \teetwo{} rule w.r.t.\ $\xi$.
		\end{enumerate}
		\item
		The \teetwo{}\teeone{} policy [applied when $i_2(\xi^L)=i_1(\xi^H)=p$]
		is as \teeone{}\teetwo{}, except that the roles of \teeone{} and \teetwo{}
		are interchanged.
	\end{enumerate}
\end{definition}

Our main result states conditions under which each of the six policies just introduced are AO.

\begin{theorem}\label{th-ao-s}
	Let Assumptions \ref{assn1} and \ref{as:polymo} hold.
	In parts 1(b), 2(b) and 2(c) below, assume moreover that $\mathbf{m}>\mathbf{m}_0$.
	\begin{enumerate}
		\item
		Assume that the single-mode condition \eqref{90} holds.
		\begin{enumerate}
			\item
			If $i_1(\xi^A)=p$ then under the \p{} policy \eqref{146} holds and this policy is AO.
			\item
			If $i_2(\xi^A)=p$ then under the \teetwo{} policy \eqref{146} holds and this policy is AO.
		\end{enumerate}
		\item
		Assume that the dual-mode condition \eqref{91} holds.
		\begin{enumerate}
			\item
			If $i_1(\xi^L)=i_1(\xi^H)=p$ then under the \p{}\p{} policy \eqref{147} holds and this policy is AO.
			\item
			If $i_2(\xi^L)=i_2(\xi^H)=p$ then under the \teetwo{}\teetwo{} policy \eqref{147} holds
			and this policy is AO.
			\item
			If $i_1(\xi^L)=i_2(\xi^H)=p$ then under the \teeone{}\teetwo{} policy \eqref{147} holds
			and this policy is AO.
			\item
			If $i_2(\xi^L)=i_1(\xi^H)=p$ then under the \teetwo{}\teeone{} policy \eqref{147} holds
			and this policy is AO.
		\end{enumerate}
	\end{enumerate}
\end{theorem}

\begin{remark}
	In \S\ref{sec51}, Theorem \ref{th3} provides more detailed information
	than \eqref{146}--\eqref{147} on the weak limit of the processes involved.
\end{remark}

\subsection{Considerations for the construction of dual-mode policies}
\label{sec:disc}

Two of the key results that we need to prove Theorem \ref{th-ao-s} are state-space collapse of the HPC (Proposition \ref{lem:ssc}), and a boundary property stating that there is asymptotically no idleness of either server when there is work in the system (Proposition \ref{lem:reflection}).
The proofs of these results differ by case/policy, and, for dual-mode policies, rely on the rules used as well as the way that workload is sampled.
Here we provide a brief description of the reasons underlying the policy definitions.

\begin{itemize}
	\item  A difficulty arises in the proof of  Proposition \ref{lem:ssc} in case 2(a), where the policy switches between two \p{} rules. Only the dual activity server in the current mode processes the HPC and the identity of the dual activity server changes when switching modes. At each switching time, if the server processing the HPC in the new mode is busy with low priority jobs and the service of the high priority job at the other server ends, no server will process high priority jobs for a time $O(n^{-1})$. In principle, this time could accumulate to let the number of high priority jobs increase to a non-negligible value. In order to prevent this, switching between \p{} rules only occurs at the time of service completion at the single activity server. When switching, this server becomes dual activity and gives priority to the HPC (which is the single activity class in \p{}). \vspace{0.2cm}
	
	\textbf{In case 2(a), there is always at least one server processing the HPC}, regardless of switching between modes.\vspace{0.2cm}
	
	\item Similarly, in order to prove Proposition \ref{lem:reflection} in case 2(b), we need to show that the number of HPC is not zero when there are low priority jobs in the system. To make sure that the high priority class does not receive too much service, switching between two \teetwo{} rules only occurs at the time of service completion at the single activity server. When switching, the single activity server becomes dual activity and now gives priority to the low priority jobs if the number of HPC jobs is low.\vspace{0.2cm}
	
	\textbf{	 In case 2(b), as long as the number of HPC jobs is below the threshold there is at most one server working on HPC jobs}, regardless of switching between modes.\vspace{0.2cm}
	
	\item In case 2(c) it should be noted that in the lower mode the system looks similar to case 1(a), so one could consider using a \p{} rule.
	Doing so, however, causes difficulty in proving Proposition \ref{lem:reflection} for this case.
	Using a \p{} rule keeps the  number of HPC jobs  $O(1)$. At the moment of switching into the \teetwo{} rule this can lead to the new single activity server, which is dedicated to HPC, incurring idle time when there is work in the system. To avoid this situation the \teeone{} rule is used instead of the \p{} rule.
	This guarantees that, at a switching time, there will not be too few HPC jobs in the system. A similar argument applies to case 2(d).
\end{itemize}

\section{The LP under the EHTC}\label{sec3}
\beginsec

In this section we prove Lemma \ref{lem1}. In addition, in a sequence of three
lemmas, we provide an explicit solution to the LP and a criterion for
determining whether the two modes are class-switched or server-switched.
Finally, Theorem \ref{thm:lowerbound} is proved.
The section is structured as follows. In \S\ref{sec31} we first prove
Lemma \ref{lem1}(1--3) based mostly on results from
\cite{ACR1}. Then we state Lemmas \ref{lem:solvelp} and \ref{lem:reordering},
which provide the LP solution, and Lemma \ref{rem:alg}
which is concerned with how the modes are paired.
Lemma \ref{lem1}(4) is then proved based on Lemma \ref{rem:alg}.
In \S \ref{sec32} we prove Lemmas \ref{lem:solvelp}--\ref{rem:alg}
and Theorem~\ref{thm:lowerbound}.

\subsection{LP-related lemmas}\label{sec31}

\noi{\bf Proof of Parts 1--3 of Lemma \ref{lem1}.}

1. Let $\xi\in\slp$. It is impossible that $\sum_i\xi_{ik}<1$ for both
$k=1$ and $2$ as this contradicts the EHTC $\rho^*=1$.
Assume then that, say, $\sum_i\xi_{i1}<1$.
Then $\xi_{11}<1$.
Define
\[
\tilde\xi=\xi+
\begin{pmatrix}
	\eps & -c\eps \\
	0 & 0
\end{pmatrix},
\]
where $c=\mu_{12}^{-1}\mu_{11}$. Then for $\eps>0$ small,
$\tilde\xi$ satisfies \eqref{01} with $\rho<1$, which again contradicts
the EHTC. This proves Part 1.

2. This follows from \cite[Lemma 2.4(4)]{ACR1}.
Note that uniqueness of the dual problem, which is a standing assumption
in \cite{ACR1}, is not used in the proof of this statement.

3. This statement follows from \cite[Lemma 2.3(1)]{ACR1}
and \cite[Lemma 2.4(4)]{ACR1}, where again the uniqueness of the dual
is not used.
\qed

\skp

The following lemma computes the two modes.

\begin{lemma}\label{lem:solvelp}
	Let the EHTCM hold.
	Then
	\begin{equation}\label{eq:ht}\sum_i\dfrac{\lambda_i}{\alpha_i}=\sum_k \beta_k=1,
	\end{equation}
	where the last equality merely expresses
	the normalization convention mentioned earlier in \S\ref{sec21}.
	Moreover, any $\xi\in\slp$ is determined by its entry $\xi_{11}$ via
	\begin{equation}\label{b1}
		\xi=
		\begin{pmatrix}
			\xi_{11} &
			\dfrac{\lambda_1}{\alpha_1\beta_2}-\dfrac{\beta_1}{\beta_2}\xi_{11}
			\\
			1-\xi_{11} & 1-\dfrac{\lambda_1}{\alpha_1\beta_2}+\dfrac{\beta_1}{\beta_2}\xi_{11}
		\end{pmatrix}.
	\end{equation}
	The two modes $\xi^{*,1}$ and $\xi^{*,2}$ can be expressed by \eqref{b1}
	with $\xi_{11}$ given by
	\begin{equation}\label{b2}
		\xi_{11}^{*,1}=\max\Big(0, \dfrac{\lambda_1}{\alpha_1\beta_1}-\dfrac{\beta_2}{\beta_1}\Big),
		\qquad
		\xi_{11}^{*,2}=\min\Big(\dfrac{\lambda_1}{\alpha_1\beta_1},1\Big).
	\end{equation}
\end{lemma}

Recall that under the nondegeneracy condition, for any mode
there is a unique relabelling of classes and servers which transforms it to
a canonical form. The following result shows that both modes, once
put in canonical form, are given by the same formula.

\begin{lemma}\label{lem:reordering}
	Let Assumption \ref{assn1} (EHTCM and nondegeneracy) hold.
	Fix $m\in\{1,2\}$.
	Relabel classes and servers so that $\xi^{*,m}$ is in canonical form. Then
	$\lambda_1> \alpha_1 \beta_1$ and
	\begin{equation}
		\label{b3}
		\xi^{*,m}=\left(\begin{array}{ll}
			1\quad\dfrac{\lambda_1}{\alpha_1\beta_2}-\dfrac{\beta_1}{\beta_2}\\
			0\quad\dfrac{\lambda_{2}}{\alpha_2 \beta_2}	\end{array}\right).
	\end{equation}
	In particular, if $\xi,\xi'\in \slp$ and there are $i,k$ such that $\xi_{ik}=\xi'_{ik}=0$
	then $\xi=\xi'$.
\end{lemma}

Lemma \ref{lem1}(4), which is yet to be proved, states that under the nondegeneracy
condition, the two modes must be either class- or server-switched.
The following lemma contains this result, and in addition provides
a criterion for distinguishing between these cases.
We will say that the {\it class-switching condition} holds if
\begin{equation}\label{93}
	\max_i\dfrac{\lambda_{i}}{\alpha_{i}}<\max_k \beta_k,
\end{equation}
and the {\it server-switching condition} holds if
\begin{equation}\label{92}
	\max_i \dfrac{\lambda_i}{\al_i}>\max_k \beta_k.
\end{equation}

\begin{lemma}
	\label{rem:alg}
	Let Assumption \ref{assn1} hold. Then both modes are nondegenerate.
	Moreover, under the class-switching condition \eqref{93}, the modes are
	class-switched, and under the server-switching condition \eqref{92},
	they are server-switched.
\end{lemma}

\noi{\bf Proof of Part 4 of Lemma \ref{lem1}.}
The statement is contained in Lemma \ref{rem:alg}.
\qed

\begin{remark}
	Note that
	cases 2(c) and 2(d) of Theorem \ref{th-ao-s} correspond to class-switched modes
	(for example, under case 2(c) one has $i_1(\xi^L)=i_2(\xi^H)=p$
	hence the non-basic activity must have switched from class $p$
	to class $q$ when moving from $\xi^L$ to $\xi^H$).
	By Lemma \ref{rem:alg}, this occurs under \eqref{93}. Also note that
	in both cases, the proposed policies apply a different rule
	for the lower and upper workload modes.
	On the other hand, cases 2(a) and 2(b) of Theorem \ref{th-ao-s}
	correspond to server-switching, and hold under \eqref{92},
	and our policies are such that
	the same rule is used for the lower and upper workload modes.
\end{remark}

\subsection{Proof of Lemmas \ref{lem:solvelp}--\ref{rem:alg}
	and Theorem \ref{thm:lowerbound}}\label{sec32}

\noi{\bf Proof of Lemma \ref{lem:solvelp}.}
In view of Lemma \ref{lem1}(1), every solution $\xi$ is column-stochastic, and,
recalling $\mu_{ik}=\al_i\beta_k$, must satisfy
\begin{align}\label{eq:lpext}
	&
	\xi_{11}\beta_{1}+\xi_{12}\beta_{2}=\frac{\lambda_1}{\alpha_1},
	\\& \notag
	\xi_{21}\beta_{1}+\xi_{22}\beta_{2}=\frac{\lambda_2}{\alpha_2},
	\\& \notag
	\xi_{11}+\xi_{21}=1,
	\\& \notag
	\xi_{12}+\xi_{22}=1,
	\\&\notag 
	\xi_{i,k}\geqslant 0,\, \ i,k\in \{1,2\}.
\end{align}
Identity \eqref{eq:ht} follows.

Next, because the expression \eqref{b1} is also column-stochastic,
proving that any solution $\xi$ is determined by $\xi_{11}$  as in \eqref{b1}
amounts to proving that $\xi_{12}$ is given as in \eqref{b1}.
This follows from the first line in \eqref{eq:lpext}.

It remains to prove \eqref{b2}. By the expression just obtained for
$\xi_{12}$ it follows that
as long as $$\xi_{11}\geqslant \dfrac{\lambda_1}{\alpha_1\beta_1}-\dfrac{\beta_2}{\beta_1},$$ we obtain $\xi_{12}\leqslant 1$. Clearly, in addition, $\xi_{11}\ges0$
must hold. Similarly, as long as 
$$\xi_{11}\leqslant \dfrac{\lambda_1}{\alpha_1\beta_1},$$ 
we obtain $\xi_{12}\geqslant 0$, and in addition, $\xi_{11}\les1$ must hold.
As a result, it is necessary that
\begin{equation}\label{b4}
	\xi_{11}\in \Big[\max\Big(0, \dfrac{\lambda_1}{\alpha_1\beta_1}-\dfrac{\beta_2}{\beta_1}\Big),\min\Big(\dfrac{\lambda_1}{\alpha_1\beta_1},1\Big)\Big].
\end{equation}
Moreover, setting $\xi_{11}$ to each of the two endpoints of the interval
indicated in \eqref{b4}
and letting $\xi$ be the corresponding expression from \eqref{b1}
gives rise to a solution satisfying all of \eqref{eq:lpext}, as can be checked
directly. Because by \eqref{b1} a solution $\xi$ is an affine function
of its entry $\xi_{11}$,
these two endpoints correspond to the two extreme points of
$\slp$, that is, to the two modes $\xi^{*,1}$, $\xi^{*,2}$.
This proves the lemma.
\qed

\noi{\bf Proof of Lemma \ref{lem:reordering}.}
Note that relations \eqref{eq:lpext} are invariant to relabeling of classes and servers.
Hence so is relation \eqref{b1}, which was derived solely from \eqref{eq:lpext}.
Let $m$ be given and assume a relabeling has been performed
to put $\xi^{*,m}$ in canonical form. Then $\xi^{*,m}$
satisfies \eqref{b1} with its first column given by $(1,0)^T$. Consequently
$\xi^{*,m}_{11}=1$. Substituting $\xi^{*,m}_{11}=1$ into \eqref{b1} proves \eqref{b3}.
Because under the nondegeneracy assumption there can be only one
zero entry,  in \eqref{b3} we have $\xi^{*,m}_{12}>0$. Hence
$\la_1>\alpha_1 \beta_1$.
The final assertion follows from \eqref{b3} using again the fact that there
can be at most one zero entry.
\qed

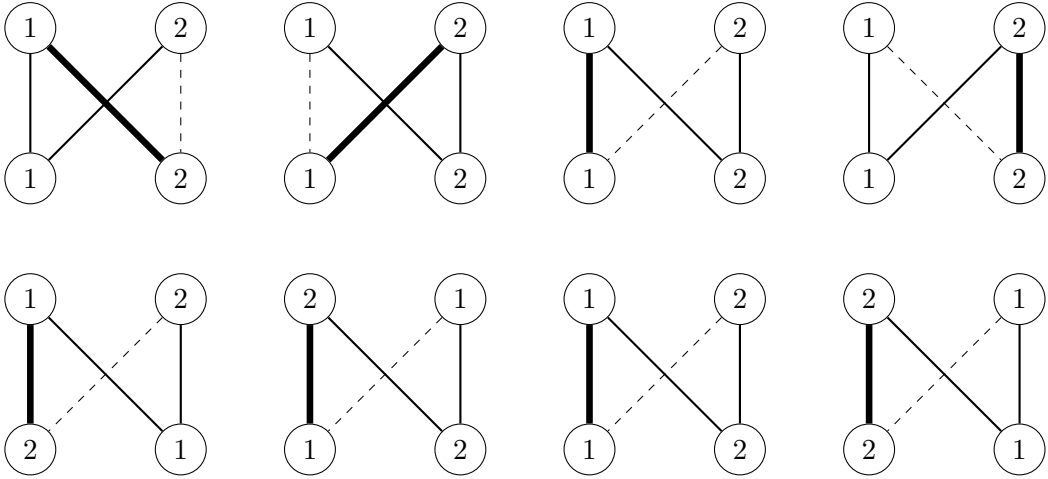
\begin{figure}[h]
	\centering
	\captionsetup[subfigure]{labelformat=empty}
	\subfloat[]{
		\begin{tikzpicture}[scale=1]
			%% nodes left
			\node[shape=circle,draw=black] (A) at (-2,1) {1};
			\node[shape=circle,draw=black] (B) at (0,1) {2};
			\node[shape=circle,draw=black] (C) at (0,-1) {2};
			\node[shape=circle,draw=black] (D) at (-2,-1) {1};
			%%edges left
			\path [thick](A) edge node[left] {} (D);
			\path [dashed] (B) edge node[right] {} (C);
			\path [color=black, line width=2.5pt](A) edge node[right] {} (C);
			\path[thick](B)edge node[right]{}(D);
			%% nodes right
		\end{tikzpicture}
	}
	\qquad
	\subfloat[]{
		\begin{tikzpicture}[scale=1]
			\node[shape=circle,draw=black] (E) at (-2,1) {1};
			\node[shape=circle,draw=black] (F) at (0,1) {2};
			\node[shape=circle,draw=black] (G) at (0,-1) {2};
			\node[shape=circle,draw=black] (H) at (-2,-1) {1};
			%% edges right
			\path [thick](E) edge node[sloped, anchor=south east] {} (G);
			\path [dashed] (E) edge node[left] {} (H);
			\path [color=black, line width=2.5pt](F) edge node[left] {} (H);
			\path[thick](F)edge node[right]{}(G);
		\end{tikzpicture}
	}
	\qquad
	\subfloat[]{
		\begin{tikzpicture}[scale=1]
			\node[shape=circle,draw=black] (A) at (-2,1) {1};
			\node[shape=circle,draw=black] (B) at (0,1) {2};
			\node[shape=circle,draw=black] (C) at (0,-1) {2};
			\node[shape=circle,draw=black] (D) at (-2,-1) {1};
			%%edges left;
			\path [thick](B) edge node[right] {} (C);
			\path [dashed] (B) edge node[left] {} (D);
			\path [color=black, line width=2.5pt](D) edge node[left] {} (A);
			\path[thick](A)edge node[sloped, anchor=south east]{}(C);
			%%nodes right
		\end{tikzpicture}
	}
	\qquad
	\subfloat[]{
		\begin{tikzpicture}[scale=1]
			\node[shape=circle,draw=black] (E) at (-2,1) {1};
			\node[shape=circle,draw=black] (F) at (0,1) {2};
			\node[shape=circle,draw=black] (G) at (0,-1) {2};
			\node[shape=circle,draw=black] (H) at (-2,-1) {1};
			%%edges right
			\path [thick](E) edge node[right] {} (H);
			\path [dashed] (E) edge node[left] {} (G);
			\path [color=black, line width=2.5pt](F) edge node[left] {} (G);
			\path[thick](F)edge node[]{}(H);
		\end{tikzpicture}	
	}
	
	%	\caption{Four possible resource allocations}
	\subfloat[Subfigure 1 list of figures text][]{
		\begin{tikzpicture}[scale=1]
			%% nodes left
			\node[shape=circle,draw=black] (A) at (-2,1) {1};
			\node[shape=circle,draw=black] (B) at (0,1) {2};
			\node[shape=circle,draw=black] (D) at (0,-1) {1};
			\node[shape=circle,draw=black] (C) at (-2,-1) {2};
			%%left edges
			\path [thick](A) edge node[left] {} (D);
			\path [dashed] (B) edge node[right] {} (C);
			\path [color=black, line width=2.5pt](A) edge node[left] {} (C);
			\path[thick](B)edge node[right]{}(D);
			%% nodes right
		\end{tikzpicture}
	}
	\qquad
	\subfloat[Subfigure 2 list of figures text][]{
		\begin{tikzpicture}[scale=1]
			\node[shape=circle,draw=black] (F) at (3,1) {2};
			\node[shape=circle,draw=black] (E) at (5,1) {1};
			\node[shape=circle,draw=black] (G) at (5,-1) {2};
			\node[shape=circle,draw=black] (H) at (3,-1) {1};
			%right edges
			\path [thick](E) edge node[right] {} (G);
			\path [dashed] (E) edge node[right] {} (H);
			\path [color=black, line width=2.5pt](F) edge node[right] {} (H);
			\path[thick](F)edge node[left]{}(G);
		\end{tikzpicture}
	}
	\qquad
	\subfloat[Subfigure 3 list of figures text][]{
		\begin{tikzpicture}[scale=1]
			\node[shape=circle,draw=black] (A) at (-2,-3) {1};
			\node[shape=circle,draw=black] (B) at (0,-3) {2};
			\node[shape=circle,draw=black] (C) at (0,-5) {2};
			\node[shape=circle,draw=black] (D) at (-2,-5) {1};
			%%left edges
			\path [thick](B) edge node[right] {} (C);
			\path [dashed] (B) edge node[left] {} (D);
			\path [color=black, line width=2.5pt](D) edge node[left] {} (A);
			\path[thick](A)edge node[left]{}(C);
		\end{tikzpicture}
	}
	\qquad
	\subfloat[Subfigure 4 list of figures text][]{
		\begin{tikzpicture}[scale=1]
			\node[shape=circle,draw=black] (F) at (3,-3) {2};
			\node[shape=circle,draw=black] (E) at (5,-3) {1};
			\node[shape=circle,draw=black] (H) at (5,-5) {1};
			\node[shape=circle,draw=black] (G) at (3,-5) {2};
			%%right edges
			\path [thick](E) edge node[right] {} (H);
			\path [dashed] (E) edge node[left] {} (G);
			\path [color=black, line width=2.5pt](F) edge node[left] {} (G);
			\path[thick](F)edge node[left]{}(H);
		\end{tikzpicture}	
	}
	\caption{Top: Four possible graphs of activities. Bottom:
		Corresponding relabelings for canonical form.\label{fig:reorder}}
\end{figure}

The four possible graphs and their relabelings are described in Figure \ref{fig:reorder}.
Namely, if $(i',k)$ is the non-basic activity in $\xi^{*,m}$, then defining
\begin{align*}
	\widetilde{\xi}^{*,m}_{11}&=\xi^{*,m}_{ik}=1, \\
	\widetilde{\xi}^{*,m}_{22}&=\xi^{*,m}_{i'k'},\\
	\widetilde{\xi}^{*,m}_{21}&=\xi^{*,m}_{i'k}=0,\text{ and }\\
	\widetilde{\xi}^{*,m}_{12}&=\xi^{*,m}_{ik'},
\end{align*}
$\widetilde{\xi}^{*,m}$ is obtained from $\xi^{*,m}$ upon relabeling
in the form of an "N".

\noi{\bf Proof of Lemma \ref{rem:alg}.}
The nondegeneracy of both modes follows from Lemma \ref{lem:reordering}.

Next, let the class switching condition \eqref{93} hold.
Because of \eqref{eq:ht},
\begin{equation}
	\label{eq:detailedcs}\max_k \beta_k>\max_i \dfrac{\lambda_i}{\al_i}\geqslant\min_i \dfrac{\lambda_i}{\al_i} >\min_k \beta_k.
\end{equation}
Recall from the proof of Lemma \ref{lem:solvelp} that the two endpoints
of the interval defined in \eqref{b4} correspond to the two modes.
Consider the right endpoint.
If the minimum in expression in \eqref{b4} is 1 then by \eqref{b1},
the non-basic activity in that mode is $(2,1)$, and moreover,
$\frac{\lambda_1}{\al_1}>\beta_1$. By \eqref{eq:ht}, this gives
$\frac{\lambda_2}{\al_2}<\beta_2.$
In view of \eqref{eq:detailedcs} this gives
$$\beta_2> \max_i \dfrac{\lambda_i}{\al_i}.$$
Hence the maximum in \eqref{b4} is 0.
By \eqref{b1}, this shows that the non-basic activity in the other mode is $(1,1)$.
If, on the other hand, the minimum in \eqref{b4} is $\frac{\lambda_1}{\al_1\beta_1}$ then $\frac{\lambda_1}{\al_1}<\beta_1$ and the non-basic activity in
the corresponding mode is $(1,2)$. Similarly, by \eqref{eq:detailedcs}
$$\min_i \dfrac{\lambda_i}{\al_i}>\beta_2.$$
This means the maximum in \eqref{b4} is not 0. The non-basic activity in the other mode is then $(2,2)$. In both cases, the two modes form a class-switched pair
as claimed.

Consider now the server switching condition \eqref{92}. because of \eqref{eq:ht},
\begin{equation}
	\label{eq:detailedss}\max_i \dfrac{\lambda_i}{\al_i} > \max_k \beta_k \geqslant\min_k \beta_k>\min_i \dfrac{\lambda_i}{\al_i}.
\end{equation}
If the minimum in \eqref{b4} is 1 then $\frac{\lambda_1}{\al_1}>\beta_1$ and the non-basic activity in one mode is $(2,1)$. By \eqref{eq:detailedss},
$\frac{\lambda_1}{\alpha_1}>\beta_2$.
Hence the maximum in \eqref{b4} is not zero and the non-basic activity in the other mode is $(2,2)$.
Finally, if the minimum in \eqref{b4} is not 1 then $\frac{\lambda_1}{\al_1}<\beta_1$ and the non-basic activity in one mode is $(1,2)$. By \eqref{eq:detailedss},
$\frac{\lambda_1}{\alpha_1}<\beta_2$.
Hence the maximum in \eqref{b4} is zero and the non-basic activity in the other mode is $(1,1)$. In both cases, the two of modes forms a server-switched pair.
\qed

\skp

\noi{\bf Proof of Theorem \ref{thm:lowerbound}.}
This lower bound is precisely the one stated in \cite[Theorem 2.6]{ACR1},
when specialized to the $2\times2$ PSS.
To prove that it is valid we must verify that the standing assumption of
\cite{ACR1}, namely \cite[Assumption 2.2]{ACR1}, holds.

First, \cite[Assumption 2.2.1]{ACR1}, which states that the EHTC holds,
is valid because of our Assumption \ref{assn1}.1.
Next, \cite[Assumption 2.2.2]{ACR1}, when translated to the notation
of this paper, states that every $\xi\in\slp$ is column-stochastic.
This holds by our Lemma \ref{lem1}.1.
It remains to show that \cite[Assumption 2.2.3]{ACR1}, which
states that the dual of \eqref{01} has a unique solution, is satisfied.

To this end we shall adopt in the remainder of this proof some notation and terminology from \cite{ACR1}.
By Lemma \ref{lem:reordering}, the non basic activity is different in both modes
$\xi^{*,1}$, $\xi^{*,2}$, for else one would have $\xi^{*,1}=\xi^{*,2}$,
contradicting the EHTCM.
As a result, with the terminology introduced in \cite[Section 2]{ACR1},
all the activities are potentially basic. Using strict complementary slackness ((36) in  \cite[Chapter 7]{schrijver}) in the same way as in
\cite[Lemma 2.3.4]{ACR1}, any $(y,z)$ solution of the dual satisfies
$$
y_i=\mu_j z_k,
\qquad i,k\in\{1,2\}, \ j=(i,k).
$$
It follows that $y_1y_2=\mu_{11}\mu_{21}z_1^2=\mu_{12}\mu_{22}z_2^2$.
By positivity of both $z_k$, this gives $z_1=cz_2$ for some constant $c>0$.
Moreover, one of the constraints of the dual problem \cite[eq.\ (2.8)]{ACR1}
is $z_1+z_2=1$. Therefore $(z_k)$ are uniquely determined.
As a consequence so are $(y_i)$, which shows that
the dual problem has at most 1 solution. The existence of a dual solution
follows from the EHTC, as shown in \cite[Lemma 2.4.2]{ACR1}.
This completes the verification of \cite[Assumption 2.2]{ACR1}.
\qed

\section{The WCP and HJB equation}\label{sec4}
\beginsec

In this section, Proposition \ref{prop0} is proved.
Lemma \ref{lem3}, that is used to prove it, contains two additional results:
An identification of optimal control systems for the WCP,
and weak uniqueness of solutions to the SDE \eqref{142},
both needed for the weak convergence proofs in \S\ref{sec5}.

\noi{\bf Proof of Part 1 of Proposition \ref{prop0}.}
This is a special case of \cite[Proposition 2.5]{ACR1}.
We comment that
the fact that $V_{\rm WCP}$ is a classical solution to \eqref{14}--\eqref{14+}
has been established already in \cite{Sheng78} (except when $b_1-b_2=\gamma$).
However, uniqueness of solutions is not covered there.
\qed

In what follows, $u$ always denotes the unique solution to \eqref{14}--\eqref{14+}.

In the following lemma, Parts 1 and 2 are largely
based on results from \cite{Sheng78}.
For completeness, we have included details on the construction from \cite{Sheng78}
in Appendix \ref{app:a}.

\begin{lemma}\label{lem3}
	\begin{enumerate}
		\item (Optimality in the single mode case).
		Assume \eqref{90} and recall that in this case $\xi^A=\xi^{*,m}$ for $m$ as in \eqref{90}.
		Then, with $u$ as above, equation \eqref{58-} of Proposition~\ref{prop0} holds. Moreover,
		$J_{\rm WCP}(z,\frS^{(1)})=V_{\rm WCP}(z)$ for the admissible control system
		$\frS^{(1)}=(\Om,\calF,\{\calF_t\},\PP,B,\X,Z^{(1)},L^{(1)})$
		where $(Z^{(1)},L^{(1)},B)$ is the RBM from \eqref{142-} (assumed to be constructed on
		the original probability space), $\calF_t=\sig\{B_s:s\in[0,t]\}$
		and $\X_t=\xi^A$.
		\item (Optimality in the dual mode case).
		Assume \eqref{91} and recall that in this case $(\xi^L,\xi^H)=(\xi^{*,m'},\xi^{*,m})$.
		Then there exists a switching point $z^*\in(0,\iy)$
		such that \eqref{58} of Proposition~\ref{prop0} holds.
		Moreover, SDE \eqref{142} possesses a weak solution
		$(\Om',\calF',\{\calF'_t\},\PP',B,Z^{(2)},L^{(2)})$.
		Furthermore, one has $J_{\rm WCP}(z,\frS^{(2)})=V_{\rm WCP}(z)$
		for the admissible control system defined by
		$\frS^{(2)}=(\Om',\calF',\{\calF'_t\},\PP',B,\X,Z^{(2)},L^{(2)})$,
		where $\X_t=\ph^*(Z^{(2)}_t)$.
		\item
		Weak uniqueness holds for solutions to SDE \eqref{142}.
	\end{enumerate}
\end{lemma}

\noi{\bf Proof of Parts 2 and 3 of Proposition \ref{prop0}.}
These results are contained in Parts 1 and 2 of Lemma~\ref{lem3}.
\qed

For Markov control problems, a map from the state space to
the control action space is often called a {\it stationary (feedback) control policy},
or a {\it policy} for short. For our WCP, a policy is thus
a measurable map $\bar\xi:\R_+\to\slp$. This term is used in \cite{Sheng78} and we adopt it
in the next proof. Parts 1 and 2 take full advantage of several results from \cite{Sheng78},
where a classical solution to equation \eqref{14}--\eqref{14+} is constructed,
and a description of an optimal policy is provided.

\noi{\bf Proof of Parts 1 and 2 of Lemma \ref{lem3}.}
The single mode condition \eqref{90} corresponds to eq.\ [41]--[42] on p.\ 105
of \cite[Section 5.3]{Sheng78}.
The dual mode condition \eqref{91} corresponds to the complementary case.
It is stated in \cite[Theorem 1, Chapter 4]{Sheng78} that a policy is optimal if and only if
the cost associated to it is $C^2$ and satisfies \cite[eq.\ (14), Chapter 4]{Sheng78},
which is the HJB equation \eqref{14}--\eqref{14+} in our notation. However,
the equation studied there is more general, and in order to reduce it to
our \eqref{14}--\eqref{14+}
one must take the switching costs to vanish (by setting the expression $K=0$),
and the reflection-absorption parameter to correspond to reflection only
(by setting $\la=1$).

Under the single mode condition, the policy constructed has the simple form
$$\bar{\xi}(z)=\xi^{*,m},\qquad z\in\R_+,$$
with the same $m$ as in \eqref{90}. It is shown that it is an optimal policy
by computing the cost associated it and showing that
it solves the HJB equation with its boundary conditions.
The fact that \eqref{58-} holds in this case
is a direct consequence of the fact that the cost associated to the single mode
policy indeed solves the HJB equation.
This completes the proof of part 1.

As for part 2, the claim that the SDE \eqref{142}
possesses a weak solution is proved in \cite[Lemma~4.1]{ACR1}.
Under the dual mode case, the policy constructed in \cite{Sheng78} is
\begin{equation}\label{p01}
	\bar{\xi}_{z^*}(z)=\xi^{*,m'}\one_{z\le z^*}+\xi^{*,m}\one_{z> z^*},
	\qquad z\in\R_+,
\end{equation}
with the same $m$ and $m'$ as in \eqref{91}.
By \cite[Theorem 4, Chapter 3]{Sheng78}, any policy of this form gives rise to
a cost function that is $C^1$ in all of $(0,\iy)$, and $C^2$ in $(0,\iy)\setminus\{z^*\}$.
The value of $z^*$ which leads to an optimal
policy is found by the principle of smooth fit, namely by requiring that
the second derivative is continuous also at $z^*$. 
The equation that states
this smoothness condition, \cite[(61), Chapter 5]{Sheng78}, turns out to have a unique solution
$z^*\in(0,\iy)$, provided that $0<b_1-b_2\ne\gamma$.
(It is assumed in \cite{Sheng78}, without loss of generality, that $b_1 \geqslant b_2$.) 
The solution to the HJB equation is then the cost associated to $\bar\xi_{z^*}$
with that value of $z^*$. The statement of part 2, by which the corresponding
control system is optimal for the WCP, therefore follows from
\cite[Theorem 1, Chapter 4]{Sheng78}. Once again, the very fact that
this function solves the HJB equation implies that equation \eqref{58} holds in this case.

The case $0<b_1-b_2=\gamma$ must be treated separately.
The argument we give is based on a scaling property.
Let $\al>0$ be a constant and multiply the SDE \eqref{15} by $\al$.
Then $Y_t:=\al Z_t$ satisfies the same dynamics \eqref{15}
with initial condition $y=\al z$ and coefficients $\tilde b(\cdot)=\al b(\cdot)$
and $\tilde\sigma(\cdot)=\al\sigma(\cdot)$.
For the function $V_{\rm WCP}$ of \eqref{51}, denote the dependence on these
two functions by writing $V_{\rm WCP}(z;b,\sig)$
(the parameter $\gamma$ is fixed throughout).
Then we immediately get the scaling property
$V_{\rm WCP}(\al z;\al b,\al\sig)=\al V_{\rm WCP}(z;b,\sig)$.
We want to use this property to show that a classical solution
to \eqref{14}--\eqref{14+} exists for $b_1-b_2=\gamma$ based
on existence for other values of $b_1,b_2$. Thus consider parameters
$b_1,b_2,\sig_2,\sig_2,\gamma$ such that $b_1-b_2=\gamma$.
Fix $\al\ne1$. Let $\tilde b_m=\al b_m$, $\tilde\sig_m=\al\sig_m$.
Let $\tilde b(\cdot)$ and $\tilde\sig(\cdot)$ be the corresponding functions,
and note by \eqref{50} that they satisfy the same scaling relation
$\tilde b(\cdot)=\al b(\cdot)$ and $\tilde\sig(\cdot)=\al\sig(\cdot)$.
Now, the new set of parameters satisfies $\tilde b_1-\tilde b_2\ne\gamma$ hence
the corresponding $\tilde u:=V_{\rm WCP}(\cdot\,;\tilde b,\tilde\sig)$ satisfies the HJB equation.
Define $u$ by scaling, that is, $u(z)=\al^{-1}\tilde u(\al z)$.
We will verify that $u$ thus defined, which is clearly smooth and satisfies
the boundary conditions \eqref{14+}, satisfies the HJB equation \eqref{14}.
Now, $u'(z)=\tilde u'(\al z)$, $u''(z)=\al\tilde u''(\al z)$.
Hence, for $z\in(0,\iy)$,
\begin{align*}
	\min_{m=1,2}[b_mu'(z)+\frac{\sig_m^2}{2}u''(z)]+z-\gamma u(z)
	&=
	\min_{m=1,2}[\al^{-1}\tilde b_m\tilde u'(\al z)
	+\frac{\tilde\sig_m^2}{2}\al^{-1}\tilde u''(\al z)]
	+z-\al^{-1}\gamma\tilde u(\al z)
	\\&=
	\al^{-1}\Big\{\min_{m=1,2}[\tilde b_m\tilde u'(y)
	+\frac{\tilde\sig_m^2}{2}\tilde u''(y)]+y-\gamma\tilde u(y)\Big\}
	\\&=0,
\end{align*}
where on the second line we have denoted $y=\al z$,
and the last line expresses the fact that $\tilde u$ satisfies
\eqref{14} with the parameters $\tilde b_m,\tilde\sig_m,\gamma$.

Finally, note that the dual mode condition \eqref{91} is left invariant
under the type of scaling of parameters that we applied.
Therefore, if \eqref{91} holds then
the function $\tilde u$ has a switching point $z^*$, and
thus the existence of a switching point for the function $u$ follows by its construction.
\qed

\noi{\bf Proof of Part 3 of Lemma \ref{lem3}.}
We slightly simplify the notation by writing \eqref{142}--\eqref{142+} as
\[
Z_t=z+\int_0^tb^*(Z_s)ds+\int_0^t\sig^*(Z_s)dB_s+L_t,
\]
where, with $a^*(x)=\sig^*(x)^2$, there exist constants
$b_0,b_1\in\R$ and $a_0>0$, $a_1>0$, such that
\[
b^*(x)=\left\{\begin{array}{ll}
	b_0,&x\le z^{*},\\
	b_1,&x> z^{*},\end{array}
\right.\qquad a^*(x)=\left\{\begin{array}{ll}
	a_0,&x\le z^{*}\\
	a_1,&x> z^{*}\end{array}
\right..
\]
In the remainder of this proof, $a^*$ and $b^*$ are written as $a$ and $b$.
Let $\scrA$ be the operator
\[\scrA f(x):=b(x)f'(x)+\frac 12a(x)f^{\prime\prime}(x),\]
on the domain 
${\cal D}(\scrA):=\{f\in C^{\infty}_0[0,\infty ):f'(0)=0\}$,
where $C^{\infty}_0[0,\infty )$ denotes the set of compactly supported
members of $C^{\infty}[0,\infty )$.
If $f\in\calD(\scrA)$ and $(Z,L,B)$ is a weak solution to the SDE then by Ito's formula,
the boundary property of $L$ and the boundary condition $f'(0)=0$,
the process $f(Z_t)-\int_0^t\scrA f(Z_s)ds$ is a martingale.
Therefore it follows from the existence result in Part 2 of the lemma
that, for every probability distribution $\nu$ on $[0,\infty)$,
there exists a solution of the $D_{[0,\infty )}[0,\infty )$ martingale
problem for $(\scrA,\nu)$. We will use the definition of
the stopped martingale problem of
Ethier and Kurtz \cite{ethkur}, Chapter 4, Section 6.
Then, for every probability distribution $\nu$ on $[0,\infty)$, 
there exists a solution of the stopped martingale problem for 
$\big(\scrA,\nu ,[0,\frac 23\,z^{*})\big)$ and of the stopped martingale problem for 
$\big(\scrA,\nu ,(\frac 13\,z^{*},\infty )\big)$. 

Define the operators $\scrA_0$ and $\scrA_1$ in the following way:
\[\scrA_0f(x):=b_0f'(x)+\frac 12a_0f^{\prime\prime}(x),\qquad\mbox{\rm on the domain }
{\cal D}(\scrA_0):={\cal D}(\scrA)\]
\[\scrA_1f(x):=b(x)f'(x)+\frac 12a(x)f^{\prime\prime}(x),\qquad\mbox{\rm on the domain }
{\cal D}(\scrA_1):=C^{\infty}_0(\R).\]
The $D_{[0,\infty )}[0,\infty )$ martingale problem for $\scrA_
0$ is well 
posed (for instance by Corollary 8.1.2, Theorem 4.4.1 and 
Proposition 4.3.1 of \cite{ethkur}). 
Therefore, for every probability distribution $\nu$ on $[0,\infty 
)$, 
there exists a unique solution of the stopped martingale problem for 
$\big(\scrA_0,\nu ,[0,\frac 23\,z^{*})\big)$ by 
Theorem 4.6.1 of \cite{ethkur}. Since, for every $f\in {\cal D}
(\scrA_0)={\cal D}(\scrA)$,
\[\scrA_0f\big|_{[0,\frac 23\,z^{*}]}=\scrA f\big|_{[0,\frac 23\,z^{*}]},\]
every solution of the stopped martingale problem for 
$\big(\scrA,\nu ,[0,\frac 23\,z^{*})\big)$ is also a solution of the stopped martingale problem for 
$\big(\scrA_0,\nu ,[0,\frac 23\,z^{*})\big)$, therefore 
the solution of the stopped martingale problem for 
$\big(\scrA,\nu ,[0,\frac 23\,z^{*})\big)$ is unique. 

Next, the $D_{\R}[0,\infty )$ martingale problem for $\scrA_1$ is well posed by Exercise  
7.3.3 of \cite{str-var} (it is shown there
that the $C_{\R}[0,\infty )$ martingale problem for $\scrA_1$ is well 
posed, but every solution of the $D_{\R}[0,\infty )$ martingale 
problem for $\scrA_1$ has paths in $C_{\R}[0,\infty )$ almost surely).
Therefore, for every probability distribution $\nu$ on $[0,\infty)$,  
there exists a unique solution of the stopped martingale problem for 
$\big(\scrA_1,\nu ,(\frac 13\,z^{*},\infty)\big)$ by Theorem 4.6.1
of \cite{ethkur}.
Since for every $f_1={\cal D}(\scrA_1)$ there exists $f\in {\cal D}
(\scrA)$ such that 
\[f\big|_{[\frac 13\,z^{*},\infty )}=f_1\big|_{[\frac 13\,z^{*},\infty 
	)},\qquad \scrA f\big|_{[\frac 13\,z^{*},\infty )}=A_1f_1\big|_{[\frac 
	13\,z^{*},\infty )},\]
every solution of the stopped martingale problem for 
$\big(\scrA,\nu ,(\frac 13\,z^{*},\infty )\big)$ is also a solution of the stopped martingale problem for 
$\big(\scrA_1,\nu ,(\frac 13\,z^{*},\infty )\big)$, therefore the solution of 
the stopped martingale problem for $\big(\scrA,\nu ,(\frac 13\,z^{*},
\infty )\big)$ is 
unique.  

Now one can apply Theorem 4.6.2 in \cite{ethkur}
with $U_1:=[0,\frac 23\,z^{*})$, $U_k:=(\frac 13\,z^{*},\infty 
)$ for $k\geq 2$, to 
conclude that, for every probability distribution $\nu$ on $[0,\infty 
)$, 
the solution of the $D_{[0,\infty )}[0,\infty )$ martingale problem for 
$(\scrA,\nu)$ is unique.

Finally, because, as already mentioned,
every solution to the SDE is a solution to the martingale
problem, the weak uniqueness of solutions to the SDE follows.
\qed

\begin{remark}\label{prop:noswitch}
	(Symmetry conditions).
	As Lemma \ref{lem3} states, \eqref{90} is a necessary and sufficient condition for
	the optimal control of the WCP to not switch between modes.
	It is natural to ask whether it can be determined if the single or the dual-mode
	condition holds under various symmetry conditions. Specifically,
	consider
	\begin{equation}\label{sy1}
		\dfrac{\hat{\mu}_{1k}}{\alpha_1}=\dfrac{\hat{\mu}_{2k}}{\alpha_2},
		\qquad k=1,2,
	\end{equation}
	\begin{equation}\label{sy2}
		\dfrac{\hat{\mu}_{i1}}{\beta_1}=\dfrac{\hat{\mu}_{i2}}{\beta_2},
		\qquad i=1,2,
	\end{equation}
	\begin{equation}\label{sy3}
		C_{S_{i1}}=C_{S_{i2}},\qquad i=1,2.
	\end{equation}
	As it turns out, each of the above three conditions is sufficient for \eqref{90}.
	This is proved in Lemma~\ref{lem:sameb} in the appendix.
	As a result, each of these conditions is sufficient for non-switching.
	These conditions can arise naturally, and occur in certain cases in the literature.
\end{remark}

\section{Asymptotic optimality results}\label{sec5}
\beginsec

In this section we prove Theorem \ref{th-ao-s}.
Toward this, an important intermediate goal is to establish
a weak convergence result, stated in Theorem \ref{th3}.
The proofs of both theorems rely on four main steps
stated in Propositions \ref{lem:m1}--\ref{lem:correctmode}.

\subsection{Weak convergence}
\label{sec51}

\subsubsection{Statement of weak convergence result}
\label{sec511}

We will adopt the following convention regarding the six cases listed in
Theorem \ref{th-ao-s}. When a statement is said to hold in a certain case
of Theorem \ref{th-ao-s}, it is meant that the assumptions as well as
the policy specified in this case are in force. For example, saying that
a certain statement holds in case 1(a) of Theorem \ref{th-ao-s} means that it holds
when the single-mode condition \eqref{90} and
the condition $i_1(\xi^A)=p$ hold, and the \p{} policy is applied, and moreover,
the weaker moment assumption $\mathbf{m}>2$ is assumed
(but $\mathbf{m}>\mathbf{m}_0$ in case 1(b), for example).
When a claim is stated without specifying a case,
it is meant that it holds in each one of the six.

In addition to the rescaled processes already defined, the weak convergence
results will be concerned with additional rescaled processes, namely
\begin{align}\label{004}
	\hat A^n_i(t) &= n^{-1/2}(A^n_i(t)-\la^n_it),
	\qquad
	\hat S^n_{ik}(t) = n^{-1/2}(S^n_{ik}(t)-\mu^n_{ik}t),
	\\
	\hat I^n_k(t)&=n^{1/2}I^n_k(t),
	\hspace{6.3em}
	\hat L^n(t)=\beta_1\hat I^n_1(t)+\beta_2\hat I^n_2(t).
	\label{004+}
\end{align}
With this notation, the balance equation \eqref{40} for $X^n$ translates
under scaling to
\begin{align}\label{eq:balancex}
	\hat X^n_i(t) &=\hat A^n_i(t)+n^{-1/2}\la^n_it
	-\sum_k\hat S^n_{ik}(T^n_{ik}(t))
	-\sum_k n^{-1/2}\mu^n_{ik}T^n_{ik}(t).\\
	\label{eq:balancexx}
	&=\hat A^n_i(t)-\sum_k\hat S^n_{ik}(T^n_{ik}(t))
	+(n^{1/2}\lambda_i+\hat{\lambda}^n_i)t-\sum_{k}T^n_{ik}(t)(n^{1/2}\alpha_i\beta_k+\hat{\mu}^n_{ik}).
\end{align}
Thus if we denote
\begin{equation}\label{096}
	\hat F^n_t=\sum_i \dfrac{\hat{A}^n_i(t)}{\alpha_i}-\sum_{ik}\dfrac{\hat S^n_{ik}(T^n_{ik}(t))}{\alpha_i}+\sum_i\dfrac{\hat{\lambda}^n_i t}{\alpha_i}-\sum_{ik}\dfrac{\hat{\mu}^n_{ik}}{\alpha_i}T^n_{ik}(t),
\end{equation}
we obtain the following representation for the workload:
\begin{equation}
	\hat W^n_t=\hat F^n(t)+\hat L^n(t).
	\label{242}
\end{equation}

The convergence of the rescaled primitives is a direct consequence of
the central limit theorem for renewal processes \cite[\S 17]{bill}. Namely,
the tuple $(\hat A^n_i,\hat S^n_{ik})$
converges to what we denote by $(A_i,S_{ik})$, comprising 6 mutually independent BMs
with zero drift and diffusivity given by the constants $\la_i^{1/2}C_{A_i}$
and $\mu_{ik}^{1/2}C_{S_{ik}}$, which in \S\ref{sec:wcpp} we have denoted by $\sig_{A,i}$
and $\sig_{S,ik}$, respectively.

\begin{theorem}\label{th3}
	Let the assumptions of Theorem \ref{th-ao-s} hold.
	Then $(T^n,\hat W^n,\hat L^n,\hat X^n) \To (T,W,L,X)$,
	where the latter is defined as follows.
	\\
	1. In cases 1(a)--(b) of Theorem \ref{th-ao-s},
	$(W,L)$ is the RBM and boundary term given by \eqref{142-}
	and initial condition $z=0$, and $T_t=\xi^At$ for all $t$.
	\\
	2.
	In cases 2(a)--(d) of Theorem \ref{th-ao-s}, $(W,L)$ is the (unique in law) weak solution
	to the SDE \eqref{142} with initial condition $z=0$,
	and letting $\X_t=\ph^*(W_t)$, one has
	$T_t=\int_0^t\X_sds$ for $t\ge0$.
	\\
	3.
	In all cases, $X_p=0$ and $X_q=\al_q W$.
\end{theorem}

The significance of this result is that it implies that, under
each of the proposed
policies, limits of the processes exist and form admissible
control systems for the WCP, which are, in view of Lemma~\ref{lem3}.1--2,
optimal.

We next present a lemma from \cite{ACR1} concerning limits of
$(\hat A^n,\hat S^n,T^n,\hat F^n)$ under general sequences of policies.
The {\it Skorohod map},
$\Gam:D_\R[0,\iy)\to D_{\R_+}[0,\iy)\times D^+_\R[0,\iy)$,
takes a function $\psi$ to a pair $(\ph,\eta)$, where
\begin{equation}\label{e02}
	\ph(t)=\psi(t)+\eta(t),\qquad \eta(t)=\sup_{0\les s\les t}\psi(s)^-,
	\qquad t\ges0.
\end{equation}
The corresponding maps $\psi\mapsto\ph$ and $\psi\mapsto\eta$ are denoted
by $\Gam_1$ and $\Gam_2$, respectively.

\begin{lemma}\label{lem4}
	Let $\{T^n\}$, $T^n\in\calA^n$ be any sequence of admissible controls
	for the QCP for which $\limsup_n\hat J^n(T^n)<\iy$.
	Then the following conclusions hold.
	The sequence $(\hat A^n,\hat S^n,T^n)$ is $C$-tight.
	Along any convergent subsequence
	where $(\hat A^n,\hat S^n,T^n)\To(A,S,T)$, one has
	\[
	(\hat A^n,\hat S^n,T^n,\hat F^n)\To(A,S,T,F),
	\]
	where $F$ is defined below in \eqref{eq:wtilde}.
	There exist on $(\Om,\calF)$
	processes $(B,\X,Z',L')$ and a filtration $\{\calH_t\}$ such that the tuple
	$\frS=(\Om,\calF,\{\calH_t\},\PP,B,\X,Z',L')$
	forms an admissible control system for the WCP with initial condition $0$.
	These processes satisfy the relations
	\begin{align}\label{eq:wtilde}
		T&=\int_0^\cdot\X_sds,
		\qquad
		(Z',L')=\Gam[F],
		\notag
		\\
		F_t &=\sum_i\frac{A_i(t)+\hat\la_it}{\al_i}
		-\sum_{i,k}\frac{S_{ik}(T_{ik}(t))+\hat\mu_{ik}T_{ik}(t)}{\al_i}\\
		&= \int_0^tb(\X_s)ds+\int_0^t\sig(\X_s)dB_s.
		\notag
	\end{align}
\end{lemma}

\proof This is the content of \cite[Lemmas 5.1 and 5.5]{ACR1}.
\qed

This lemma is our starting point for proving Theorem \ref{th3}.
Whereas it relates limits of processes associated with the QCP
to an admissible control system for the WCP,
note that it does not make any claim regarding
$\hat X^n$ or $\hat W^n$, hence by itself is not sufficient to relate
the prelimit cost (defined in terms of $\hat X^n$)
to the WCP cost. In particular,
the pair $(Z',L')$ need not be the weak limit of
$(\hat W^n,\hat L^n)$.
To proceed one must show that under the proposed policies,
along the sequence specified in Lemma \ref{lem4}, one has
\begin{equation}\label{100}
	(\hat A^n,\hat S^n,T^n,\hat F^n,\hat W^n,\hat L^n)
	\To
	(A,S,T,F,W,L)
\end{equation}
where $(W,L)=(Z',L')$. Once this is achieved, the lemma guarantees
that the weak limit $(W,L)$ satisfies
\[
W_t=\int_0^tb(\X_s)ds+\int_0^t\sig(\X_s)dB_s+L_t,
\]
with $\int_{[0,\iy)} W_tdL_t=0$,
assuring that the limit tuple indeed forms an admissible system for the WCP.
To go from here to the statements made in Theorem \ref{th3},
one further needs to show that $(W,L)$ are as given in
\eqref{142-} or \eqref{142}, and that $X_p=0$.

The Propositions presented below address these issues as follows.
Propositions \ref{lem:m1} provides uniform integrability required
to eventually deduce Theorem \ref{th-ao-s} from Theorem \ref{th3},
and in addition ensures that the prelimit cost remains bounded
so Lemma \ref{lem4} may be applied.
Proposition \ref{lem:ssc} shows that $\hat X^n_p\to0$ in probability.
Proposition \ref{lem:reflection} states precisely what is needed
to attain \eqref{100}. Finally, Proposition \ref{lem:correctmode}
implies that \eqref{142} holds in the dual mode case.

\subsubsection{Main steps toward weak convergence}
\label{sec:512}

Throughout what follows, the assumptions of Theorem \ref{th-ao-s} are in force.
The four main steps required to achieve weak convergence, and later,
Theorem \ref{th-ao-s}, are as follows.

\begin{proposition}\label{lem:m1}
	There exists $\eps_0>0$ such that
	\[
	\limsup_n \E\int_0^\iy e^{-\gamma t}(\hat H^n_t)^{1+\eps_0}dt<\iy.
	\]
\end{proposition}
As already mentioned, this uniform integrability result
will allow us to deduce convergence of the costs, as stated
in Theorem~\ref{th-ao-s}, from the convergence stated
in Theorem \ref{th3}.
Moreover, because it also implies boundedness of the cost under the proposed policies, it enables us to use Lemma \ref{lem4}.
The proof is given in \S \ref{sec:ui}.

\begin{proposition}\label{lem:ssc}
	For every $t_0>0$, as $n\to\iy$,
	\[
	\mathbb{P}\left(\|\hat{X}^n_p\|_{t_0}\geqslant 2\hat\Th^n\right)\to 0.
	\]
\end{proposition}
The above type of result is often referred to as state-space collapse (SSC), as it asserts that asymptotically
all workload is kept in one buffer, a property crucially used in establishing
the one-dimensional state space description of the limiting dynamics, as well as asymptotic optimality, since all workload is held in the 'less expensive' class.
It is proved in \S \ref{sec:ssc}.

\begin{proposition}\label{lem:reflection}
	Consider a subsequence as in Lemma \ref{lem4},
	where $(\hat A^n,\hat S^n,T^n)\To(A,S,T)$.
	Then along this sequence one has
	$(\hat A^n,\hat S^n,T^n,\hat F^n,\hat W^n,\hat L^n)
	\To(A,S,T,F,W,L)$,
	where $F$ is given by \eqref{eq:wtilde}, and $(W,L)=\Gam(F)$.
	In particular, $(\hat{W}^n,\hat L^n)$ is a C-tight sequence,
	and the conclusions of Lemma \ref{lem4} hold with $(W,L)=(Z',L')$.
\end{proposition}

The proof appears in \S\S \ref{sec:boundarybeh}--\ref{sec525}.

Finally, the policies \p{} and \teetwo{} that are proposed under the single mode condition
do not use the non-basic activity of the active mode $\xi^A$. For the four policies
employed under the dual-mode condition, we need the following control over the
use of the non-basic activities corresponding to $\xi^L$ and $\xi^H$.
\begin{proposition}\label{lem:correctmode}
	Consider the same subsequence as in Proposition \ref{lem:reflection} and cases 2(a)--(d) of Theorem \ref{th-ao-s}.
	If $(i^l,k^l)$ (resp. $(i^h,k^h)$) denotes the non-basic activity in $\xi^L$ (resp. $\xi^H$),
	then for any $t_0>0$ and $\eps>0$,
	\begin{equation}\label{eq:modes}
		\int_{0}^{t_0}\mathds{1}_{\hat W^n_t \leqslant z^*-\eps}dT^n_{i^lk^l}(t)\to 0,
		\qquad
		\int_{0}^{t_0}\mathds{1}_{\hat W^n_t \geqslant z^*+\eps}dT^n_{i^hk^h}(t)\to 0,
	\end{equation}
	in probability, as $n\to\iy$.
\end{proposition}
To explain the role of this proposition, recall that
if $\xi\in\slp$ then the condition $\xi_{ik}=0$ for some activity $(i,k)$ not only
implies that $\xi$ is one is the modes $\xi^{*,1}$ or $\xi^{*,2}$
but also identifies which one by Lemma \ref{lem:reordering}.
Since any limit $T$ of $T^n$ is given as $\int_0^\cdot\X_sds$, $\X_t\in\slp$,
this proposition implies that when workload is either below or above the switching 
point, the resource allocation
asymptotically follows the respective mode of operation $\xi^L$ or $\xi^H$.
This is very close to stating that the pair $(W,L)$
follows the SDE \eqref{142}, and indeed is the basis for proving this fact.
The proof of this proposition appears in \S \ref{sec:fs}.

\subsubsection{Proof of weak convergence}
\label{sec513}

Here we prove Theorem \ref{th3} based on the four propositions.

\noi{\bf Proof of Theorem \ref{th3}.}
First, by Proposition \ref{lem:m1} one has that $\limsup_n\hat J^n(T^n)<\iy$
for each one of the relevant policies $T^n$. As a result, the assumptions
of Lemma \ref{lem4} and Proposition \ref{lem:reflection} are valid.
To summarize the conclusions from these results,
fix a subsequence along which $(\hat A^n,\hat S^n,T^n)\To(A,S,T)$.
Then there exists a tuple $(\{\calH_t\},F,W,L,\X,B)$
such that, along this sequence,
\[
(\hat A^n,\hat S^n,\hat F^n,T^n,\hat W^n,\hat L^n)\To (A,S,F,T,W,L),
\]
where $F$ is given in terms of $A$, $S$, $T$
by \eqref{eq:wtilde},
and $W$ and $L$ are given by $(W,L)=\Gam(F)$.
Moreover,
$\frS=(\Om,\calF,\{\calH_t\},\PP,B,\X,W,L)$ forms an admissible control
for initial condition $0$, and $T=\int_0^\cdot\X_sds$. In particular,
\begin{equation}\label{150}
	W=F+L=\int_0^\cdot b(\X_s)ds+\int_0^\cdot\sigma(\X_s)dB_s+L_t.
\end{equation}

Next, by Proposition \ref{lem:ssc}, $\hat X^n_p\to0$ in probability.
Also, by \eqref{241}, $\hat X^n_q=\al_q\hat W^n-\al_p\hat X^n_p$,
and therefore we now have, along the subsequence,
$(\hat A^n,\hat S^n,\hat F^n,T^n,\hat W^n,\hat L^n,\hat X^n)
\To (A,S,F,T,W,L,X)$, where we denote
\begin{equation}\label{151}
	X_p=0, \qquad X_q=\al_q W.
\end{equation}

Note that the control system $\frS$ thus constructed
may depend on the subsequence. However, consider the following.

\noi{\it Claim.
	In cases 1(a)--(b) of Theorem \ref{th-ao-s},
	$(W,L)$ is the RBM and its boundary term given by \eqref{142-},
	and $T_t=\xi^At$ for all $t$.
	In cases 2(a)--(d) of Theorem \ref{th-ao-s}, $(W,L)$ is a weak solution
	to the SDE \eqref{142}, and moreover,
	$\X_t=\ph^*(W_t)$ for a.e.\ $t$.
}

Suppose the above claim holds true. Then the law of $(W,L)$ is uniquely determined:
in case 1 as a RBM; in case 2 as a weak solution to \eqref{142}, for which weak 
uniqueness holds by Lemma \ref{lem3}.3.
In particular, this law does not depend on the subsequence.
Moreover, since by this claim and \eqref{151},
the pair of processes $(T,X)$ is uniquely determined by $W$ (away from a $\PP$-null set),
it follows that the law of $(T,W,L,X)$ does not depend on the subsequence.
This yields the convergence
$(T^n,\hat W^n,\hat L^n,\hat X^n)\To(T,W,L,X)$ along the full sequence
and completes the proof of the result.
In what follows, the claim is proved.

Consider first the single mode case, namely cases 1(a)--(b) of Theorem \ref{th-ao-s}.
The policies employed are \p{} and \teetwo{},
and both do not use the non-basic activity of the active mode $\xi^A$.
In other words, if we denote this activity by $(i^a,k^a)$ then under these
policies, $T^n_{i^ak^a}=0$ for all $n$. As a consequence, the limit process
$T$ must satisfy $T_{i^ak^a}(t)=0$ a.s., hence
$\X_{i^ak^a}(t)=0$ for a.e.\ $t$, a.s.
By the uniqueness statement made in Lemma \ref{lem:reordering},
whenever $\tilde\xi\in\slp$ and $\tilde\xi_{i^ak^a}=0$, one must have
$\tilde\xi=\xi^A$.
It follows that $\X_t=\xi^A$ for a.e.\ $t$, and hence $T_t=\xi^At$.
As a consequence, \eqref{150}
holds as $W_t=b(\xi^A)t+\sig(\xi^A)B_t+L_t$.
That is, the pair $(W,L)$ satisfies \eqref{142-}.
This proves the first part of the claim.

Next consider the dual mode, namely, cases 2(a)--(d) of Theorem \ref{th-ao-s}.
First, we will show based on Proposition \ref{lem:correctmode} that,
for every $\eps>0$ and $t_0>0$,
\begin{equation}\label{102}
	\int_0^{t_0}\one_{\{W_t<z^*-\eps\}}dT_{i^lk^l}(t)=0,
	\qquad
	\int_0^{t_0}\one_{\{W_t>z^*+\eps\}}dT_{i^hk^h}(t)=0
\end{equation}
holds a.s. Let $g$ be a continuous function such that
\begin{equation}
	\label{eq:g}\one_{w<z^*-\eps}\leqslant g(w)\leqslant \one_{w<z^*-\frac{\eps}{2}}, \qquad w\in\R_+.
\end{equation}
By the continuous mapping theorem, we have along the subsequence
$(g(\hat{W}^n),T^n)\Rightarrow (g(W),T)$.
In addition, $T^n$ is continuous with finite variation
over compacts. By Theorem 2.2 of \cite{kur-pro},
\begin{equation}\label{101}
	\int_0^{\cdot}g(\hat{W}^n_t)dT^n_{i^lk^l}(t)
	\Rightarrow \int_0^{\cdot}g(W_t)dT_{i^lk^l}(t).
\end{equation}
By Proposition \ref{lem:correctmode},
$\int_0^{t_0}\one_{\{\hat{W}^n_t<z^*-\frac{\eps}{2}\}}dT^n_{i^lk^l}(t)\to0$
in probability. Hence the LHS in \eqref{101} converges
to zero in probability. Thus the RHS in \eqref{101} equals zero a.s.,
and therefore by the first inequality in \eqref{eq:g},
the first part of \eqref{102} is proved.
The second part of \eqref{102} is proved analogously.

Next, clearly $\X_{i^lk^l}(t)\one_{\{\X_t=\xi^L\}}=0$. Hence by \eqref{102},
$\int_0^{t_0}\one_{\{W_t<z^*-\eps\}}\one_{\{\X_t\ne\xi^L\}}\X_{i^lk^l}(t)dt=0$.
Arguing again by the uniqueness statement in Lemma \ref{lem:reordering},
one has $\X_{i^lk^l}(t)>0$ whenever $\X_t\ne\xi^L$.
It follows that, a.s.,
\begin{equation}\label{152}
	\int_0^{t_0}\one_{\{W_t<z^*-\eps\}}\one_{\{\X_t\ne\xi^L\}}dt=0,
	\qquad
	\int_0^{t_0}\one_{\{W_t>z^*+\eps\}}\one_{\{\X_t\ne\xi^H\}}dt=0,
\end{equation}
where the second equality is proved analogously to the first one.
Going back to \eqref{150}, note that by \eqref{152} one has both
$\int_0^{t_0}b(\X_s)\one_{\{W_t<z^*-\eps\}}\one_{\{\X_t\ne\xi^L\}}ds=0$
and
$\int_0^{t_0}\sig(\X_s)\one_{\{W_t<z^*-\eps\}}\one_{\{\X_t\ne\xi^L\}}dB_s=0$.
A similar statement hold for $\{W_t>z^*+\eps\}$ and $\xi^H$.
Hence by \eqref{150} and the definition
\eqref{142+} of the functions $b^*$ and $\sigma^*$, it follows that
\begin{multline*} W_t=\int_{0}^{t}\mathds{1}_{\lvert W_s-z^*\rvert\geqslant\eps}b^*(W_s)ds+\int_{0}^{t}\mathds{1}_{\lvert W_s- z^*\rvert\geqslant \eps}\sigma^*(W_s)dB_s\\
	+\int_{0}^{t}\mathds{1}_{\lvert W_s-z^*\rvert<\eps}b(\X_s)ds+\int_{0}^{t}\mathds{1}_{\lvert W_s- z^*\rvert< \eps}\sigma(\X_s)dB_s+L_t,
\end{multline*}
for $t\le t_0$. Because $t_0$ is arbitrary, this is true for all $t$.
Hence, denoting $\del^\eps_t=\mathds{1}_{\lvert W_t-z^*\rvert<\eps}$
and using the boundedness of $b(\cdot)$,
\[
U_t:=\left\lvert W_t-\int_{0}^{t}b^*(W_s)ds-\int_{0}^{t}\sigma^*(W_s)dB_s-L_t\right\rvert\label{eq:diffw}\le \gamma_t^\eps+|\tilde\gamma_t^\eps|
+|\check\gamma_t^\eps|,
\]
where
\[
\gamma_t^\eps=c\int_{0}^{t}\del^\eps_sds,
\qquad
\tilde\gamma_t^\eps=
\int_{0}^{t}\del^\eps_s\sig(\X_s)dB_s,
\qquad
\check\gamma^\eps_t=
\int_{0}^{t}\del^\eps_s\sig^*(W_s)dB_s.
\]
Notice that $U_t$ does not depend on $\eps$,
so if we manage to prove that the RHS converges to zero in probability
as $\eps\downarrow0$, it follows that, for every $t$, $U_t=0$ a.s.
Hence by continuity of this process, $U_t=0$ for all $t$, a.s.
That is, the processes $(W,L,B)$ satisfy \eqref{142} a.s.

To this end,
apply the occupation times formula, \cite[Corollary VI.1.6]{rev-yor},
by which for any continuous semimartingale $Y$, one has, a.s.,
$$\int_0^t\one_{\{Y_s=0\}}d\lan Y,Y\ran_s
=\int_{-\infty}^{\infty}\one_{y=0}L^y_t(Y)dy,$$
with $L^y(Y)$ the local time of $Y$ at $y$.
Consider the above with $Y_t=W_t-z^*$.
Clearly the RHS in the above display is zero.
Moreover, $\langle Y,Y\rangle_t=\int_0^t\sigma(\X_s)^2ds$,
and since $\sigma$ is bounded away from zero, we obtain that
$\int_0^t\one_{\{W_s=z^*\}}ds=0$ a.s.
For fixed $\om$, one has for all $t$, $\one_{\{|W_t-z^*|<\eps\}}
\to\one_{\{W_t=z^*\}}$ as $\eps\downarrow0$.
Hence, by dominated convergence,
$\gamma^\eps_t\to0$ a.s.\ as $\eps\downarrow0$.

Next, for the stochastic integral $\tilde\gamma^\eps$ we have
\[
\E(\tilde\gamma^\eps)^2=\E\int_0^t(\del_s^\eps\sig(\X_s))^2ds
\le c\E\int_0^t\del_s^\eps ds=c\E(\gamma^\eps_t).
\]
By local boundedness of $\gamma^\eps$ and its a.s.\ convergence
to $0$ as $\eps\downarrow0$, it follows that $\tilde\gamma^\eps_t\to0$ in probability.
A similar argument holds for $\check\gamma^\eps_t$.
We conclude that $(W,L,B)$ satisfies \eqref{142}.
Finally, by \eqref{152} and the fact $\int_0^t\one_{\{W_s=z^*\}}ds=0$
a.s.\ just proved, one has $\X_t=\ph^*(W_t)$ for a.e.\ $t$, a.s.,
which completes the proof of the claim, and also of the result.
\qed

\subsubsection{Proof of main result}\label{sec:514}

As a consequence of Theorem \ref{th3} and Proposition \ref{lem:m1}
we have the following.

\noi{\bf Proof of Theorem \ref{th-ao-s}.}
The weak convergence statements asserted in the theorem are already established
in Theorem \ref{th3}. For the AO results we need to show that in each of the six cases
$\hat J^n(T^n)\to V_0=h_q\al_qV_{\rm WCP}(0)$.
Combining Theorem \ref{th3} with
the identification of an optimal control for the WCP
given in Lemma \ref{lem3} shows that
$\hat X^n\To X$ where $X_p=0$ and $X_q=\al_qW$,
$W$ is given by \eqref{142-} or \eqref{142} in the respective cases,
and moreover
\[
V_0=h_q\al_q\E\int_0^\iy e^{-\gamma t}W_tdt.
\]
The convergence stated above implies
\begin{equation}\label{149}
	\hat H^n_t=h\cdot\hat X^n_t\To h_q\al_q W_t.
\end{equation}
By \eqref{42}, $\hat J^n(T^n)= \E\int_0^\iy e^{-\gamma t}\hat H^n_tdt$.
Hence the convergence $\hat J^n(T^n)\to V_0$ will follow from \eqref{149}
once uniform integrability is established. Arguing
along the lines of \cite{bw1} (pp.\ 640--643), introduce the measure
$dm=\gamma e^{-\gamma t}dt$ on $(\R_+,\calR_+)$ and invoke
the Skorohod representation theorem to obtain from \eqref{149}
that $\hat H^n\to h_q\al_qW$ $(m\times \PP)$-a.e.
Accordingly, uniform integrability of $H^n$ w.r.t.\ $m\times\PP$
suffices to obtain $\hat J^n(T^n)\to V_0$.
However, this is ensured by Proposition \ref{lem:m1}, for
\[
\limsup_n\int_{[0,\iy)\times\Om}(\hat H^n)^{1+\eps_0}d(m\times\PP)
=\gamma\limsup_n\E\int_0^\iy e^{-\gamma t}(\hat H^n_t)^{1+\eps_0}dt<\iy,
\]
and the result is proved.
\qed

\subsection{Proof of Propositions \ref{lem:m1}--\ref{lem:correctmode}}\label{sec52}

In \S\ref{sec521}, we provide several useful estimates.
Then, Propositions \ref{lem:m1}, \ref{lem:ssc}, \ref{lem:reflection} and
\ref{lem:correctmode} are proved in
\S\ref{sec:ui}, \S\ref{sec:ssc}, \S\S\ref{sec:boundarybeh}--\ref{sec525},
and \S\ref{sec:fs}, respectively.

\subsubsection{Auxiliary lemmas}
\label{sec521}

Two estimates on the rescaled primitives,
$\hat A^n_i$ and $\hat S^n_{ik}$, are provided in
\eqref{250} and Lemma \ref{lem10}, and a certain estimate on
the maximum service duration is given in Lemma \ref{lem:1service}.

Because the assumptions on $A^n_i$ and $S^n_{ik}$ are similar, the estimates
are stated for $\hat A^n_i$ but apply also for $\hat S^n_{ik}$.
Recall that $\check a_i(l)$ are interarrival times of $A_i$
and thus $a^n_i(l):=(\la^n_i)^{-1}\check a_i(l)$
are the interarrival times of $A^n_i$. Recall $E[\check a_i(1)^\m]<\iy$
for a constant $\m>2$.
The first estimate is
\cite[Theorem 4]{krichagina} which states that
for any $2\le\kappa\le\m$,
\begin{equation}\label{250}
	E[\|\hat A^n_i\|_t^\kappa]\le c(1+t)^{\kappa/2}
\end{equation}
for a constant $c$ that does not depend on $n$ or $t$.
The next useful estimate is as follows.
\begin{lemma}\label{lem10}
	Let $\nu_1,\nu_2\in(0,1)$ be such that $\nu_1\les \nu_2+\frac{1}{2}$ and assume that $h_0:=(\frac{\m}{2}-1)\nu_1-\m\nu_2>0$
	(note, in particular, that for every $\nu_1\in(0,\frac{1}{2}]$ there exists $\nu_2>0$
	satisfying these conditions). Fix $c_1>0$. Then
	for any $h<h_0$,
	\[
	P(w_{t_0}(\hat A^n_i,n^{-\nu_1})\ges c_1n^{-\nu_2})
	\les cn^{-h}(1+t_0)^{\m/2}, \qquad n\in\N,\, t_0\ge1,
	\]
	where $c=c(c_1,\nu_1,\nu_2,\m)$ does not depend on $n$ or $t_0$.
\end{lemma}
\begin{remark}
	We sometimes use the balance equation in the following form,
	which follows from \eqref{eq:balancexx}, namely
	\begin{equation}\label{eq:balancexxx}
		\hat{X}^n_i(t)=\hat f^n_i(t)
		+n^{1/2}\int_{0}^{t}\Big(\lambda_i-\sum_k\mu_{ik}\X^n_{ik}(s)\Big)ds
	\end{equation}	
	where
	\[\hat{f}_i^n(t)\coloneqq \hat A^n_i(t)-\sum_k\hat S^n_{ik}(T^n_{ik}(t))+t\hat{\lambda}^n_i-\sum_{k}T^n_{ik}(t)\hat{\mu}^n_{ik}.\]
	Note that $\hat F^n$ of \eqref{096} is given by
	$\hat{F}^n(t)=\sum_{i}\al_i^{-1}\hat{f}^n_i(t).$
	Then using the $1$-Lipschitz property of the trajectories of $T^n_{ik}$,
	it is easy to see that both estimates above imply some estimates for $\hat f^n_i$. In particular, by \eqref{250},
	$\E\|\hat f^n_i\|_t^\kappa\le c(1+t)^\kappa$.
	Moreover, under the assumptions of Lemma \ref{lem10} and the
	additional assumption that $\nu_2<\nu_1$, the conclusion of the lemma holds
	for $\hat f^n_i$ and $\hat F^n$.
	\label{rem:stillok}
\end{remark}

\noi{\bf Proof.} In this proof, $\m$ is written as $m$.
Note first that  it suffices to prove
\begin{equation}\label{200}
	\PP(w_{t_0}(\hat A^n,n^{-\nu_1})\ges n^{-\nu_2})
	\les cn^{-h_0}(1+t_0)^{m/2}, \qquad n\in\N,
\end{equation}
where $c=c(\nu_1,\nu_2,m)$ does not depend on $n$ or $t_0$.
Indeed, if $c_1\ges 1$, the result follows directly from \eqref{200}.
If $c_1\in(0,1)$, let $\bar \nu_2>\nu_2$ be such that it
satisfies all hypotheses
of the lemma. Namely $\bar \nu_2\in(0,1)$, $\nu_1\le \bar\nu_2+\frac{1}{2}$,
$\bar h_0:=(\frac{m}{2}-1)\nu_1-m\bar \nu_2>0$. Then by \eqref{200},
for $n$ such that $c_1n^{-\nu_2}\ges n^{-\bar \nu_2}$,
\[
\PP(w_{t_0}(\hat A^n,n^{-\nu_1})\ges c_1n^{-\nu_2})
\les c(1+t_0)^{\frac{m}{2}}n^{-\bar h_0}.
\]
Moreover, $\bar h_0$ can be made arbitrarily close to $h_0$
by choosing $\bar \nu_2$ close to $\nu_2$. This shows that the desired
inequality holds for all large $n$, and by making $c=c(c_1,\nu_1,\nu_2,m)$
larger, for all $n\in\N$.

We now prove \eqref{200}.
As before, $c$ denotes a positive constant whose value may change
from line to line; here it may depend on $(\nu_1,\nu_2,m)$ but not on $n$, $t_0$.
By \eqref{004},
\[
\hat A^n(t)-\hat A^n(s)=n^{-1/2}(A^n(t)-A^n(s))
-n^{-1/2}\la^n(t-s).
\]
Thus
\begin{equation}\label{003}
	\PP(w_{t_0}(\hat A^n,n^{-\nu_1})\ges n^{-\nu_2})\les \PP(A^n(t_0)>2\la^nt_0)
	+\PP(\Om^n_+)+\PP(\Om^n_-),
\end{equation}
where
\[
\Om^n_+=\{A^n(t_0)\les 2\la^nt_0,\exists s,t\in[0,t_0],
0\les t-s\les n^{-\nu_1}, A^n(t)-A^n(s)\ges n^{\frac12-\nu_2}+\la^n(t-s)\},
\]
\[
\Om^n_-=\{A^n(t_0)\les 2\la^nt_0,\exists s,t\in[0,t_0],
0\les t-s\les n^{-\nu_1}, A^n(t)-A^n(s)\les -n^{\frac12-\nu_2}+\la^n(t-s)\}.
\]
Consider the event $\Om^n_+$.
Because $a^n(l)$ are the interarrival times of $A^n$, on this event
there must exist $l^0\les 2\la^nt_0$ and $R\les n^{-\nu_1}$
such that
\[
\sum_{l=l^0}^{l^0+n^{\frac12-\nu_2}+\la^n R}a^n(l)\les R.
\]
Recall that $\E\check a(l)=1$. Letting $\bar a^n(l)=
a^n(l)-(\la^n)^{-1}$, we have $\E\bar a^n(l)=0$.
Then taking $r=\la^n R$, using $\la^n\les c_1n$ where $c_1=\sup\frac{\la^n}{n}<\iy$, we have
\begin{align*}
	\PP(\Om^n_+) &\les
	\PP(\exists j\les 2\la^nt_0,\exists r\les n^{-\nu_1}\la^n,
	\sum_{l=l^0+1}^{l^0+n^{\frac12-\nu_2}+r}a^n(l)\les (\la^n)^{-1}r)\\
	&\les
	\PP(\exists j\les 2c_1nt_0,\exists r\les c_1n^{1-\nu_1},
	\sum_{l=l^0+1}^{l^0+n^{\frac12-\nu_2}+r}\bar a^n(l)\les -(\la^n)^{-1}n^{\frac12-\nu_2}).
\end{align*}
Let $M_{l^1}^n=\sum_{l=1}^{l^1}\bar a^n(l)$, $l^1=0,1,2,\ldots$,
and note that it is a martingale.
For a real-valued function $X$ on $\Z_+$ let
\[
{\rm osc}(X,l_1,l_2)=\max\{|X(l_3)-X(l_4)|:l_3,l_4\in[l_1,l_2]\},
\qquad 0\le l_1\le l_2.
\]
Then using $\frac{1}{2}-\nu_2\les 1-\nu_1$ and denoting $\rho=[c_1n^{1-\nu_1}]$,
\begin{align*}
	\PP(\Om^n_+) &\les
	\PP(\exists j\les 2c_1nt_0, {\rm osc}(M^n,l^0,l^0+2c_1n^{1-\nu_1})\ges c_1^{-1}n^{-\frac12-\nu_2})\\
	&\les
	\PP(\exists j\in[0,2c_1nt_0]\cap\{0,\rho,2\rho,\ldots\},\ {\rm osc}(M^n,l^0,l^0+\rho)\ges
	c_1^{-1}n^{-\frac12-\nu_2}/3)\\
	&\les
	1-(1-\PP(\|M^n\|_\rho\ges c_1^{-1}n^{-\frac12-\nu_2}/6))^{ct_0n^{\nu_1}}.
\end{align*}
Because $n^{-1}\la^n\to\la>0$, we have the lower bound $\la^n>cn$ for some
$c>0$ and all large $n$. Hence Burkholder's inequality shows that
\[
\E(\|M^n\|_\rho)^m\les c\E(\rho|\bar a^n(1)|^2)^{\frac{m}{2}}\les
cn^{(1-\nu_1)\frac{m}{2}}(\la^n)^{-m}\les cn^{-(1+\nu_1)\frac{m}{2}}.
\]
Hence for any $a>0$, $\PP(\|M^n\|_\rho>a)\les c a^{-m}n^{-(1+\nu_1)\frac{m}{2}}$,
and therefore
\[
\PP(\Om^n_+) \les 1-(1-cn^{-(\nu_1-2\nu_2)\frac{m}{2}})^{ct_0n^{\nu_1}}.
\]
Note that $\nu_1>2\nu_2$ by the assumption $h_0>0$. Now, if $\eps\in(0,\frac{1}{2})$
and $a>0$ then, with $c_2=2\log 2$,
$1-(1-\eps)^a\les 1-e^{-c_2a\eps}\les c_2a\eps$.
This gives
\[
\PP(\Om^n_+) \les cn^{-(\nu_1-2\nu_2)\frac{m}{2}}t_0n^{\nu_1}=ct_0n^{-h_0}.
\]
By a similar argument, the same estimate holds for $\PP(\Om^n_-)$.
An application of \eqref{250} gives
\[
\PP(A^n(t_0)>2\la^nt_0)
\les \PP(\|\hat A^n\|_{t_0}>cn^{1/2})\les c\frac{\E(\|\hat A^n\|_{t_0}^m)}{n^{m/2}}
\les
cn^{-\frac{m}{2}}(1+t_0)^{\frac{m}{2}}.
\]
Hence by \eqref{003},
\[
\PP(w_{t_0}(\hat A^n,n^{-\nu_1})\ges n^{-\nu_2})
\les ct_0n^{-h_0}+c(1+t_0)^{\frac{m}{2}}n^{-\frac{m}{2}}.
\]
It follows from $\nu_1,\nu_2\in(0,1)$ that $h_0<\frac{m}{2}$.
The result follows.
\qed
\vspace{0.2cm}

Next, we give an estimate on the maximal service duration and interarrival time
up to a given time.
The {\it time in service by $t$} of a given job is defined as the time that the job
has spent in service up to time $t$.
Let ${\rm TIS}(n,i,k,l,t)$ denote the time in service by $t$ of the $l$th job in activity $(i,k)$.
If service to job $l$ has completed by time $t$ then clearly ${\rm TIS}(n,i,k,l,t)=u^n_{ik}(l)$,
but if it is still in service,
${\rm TIS}(n,i,k,l,t)<u^n_{ik}(l)$. Of course, ${\rm TIS}(n,i,k,l,t)=0$ for jobs for which service has not
started by $t$.
For $t>0$ and a real-valued path $\ph$, denote
\[
\La(\ph,t)=\sup\{t_2-t_1: \ 0\le t_1\le t_2\le t,\ \ph_{t_2}=\ph_{t_1}\}.
\]
Then, for activity $(i,k)$, the maximal time in service by time $t$,
namely $\sup_l{\rm TIS}(n,i,k,l)$, is bounded above by $\La(S^n_{ik},t)$.
We will need an upper bound on the service time as well as the interarrival times,
and to this end define
\begin{equation}\label{eq:maxservice}
	e^n_{\rm max}(t)=\max_{i,k}\La(S^n_{ik},t)\vee\max_i\La(A^n_i,t).
\end{equation}
This process also bounds from above all service durations
completed by time $t$.
\begin{lemma}
	\label{lem:1service}
	One has $\E[(e^n_{\rm max})^2]\le cn^{-1}(1+t)$
	for a constant $c$ that does not depend on $n$ or $t$. Moreover,
	for any $t<\infty$ and $c_1>0$,
	$\mathbb{P}\left(e_{\max}^n(t)\geqslant c_1 n^{\bar a-1}\right)\to 0,$
	as $n\to +\infty$.
\end{lemma}

\proof
Fix $i,k$. Denote $Y^n_t=\La(S^n_{ik},t)$.
For any $u>0$, if $Y^n_t\ge u$ then there must exist $0\le t_1<t_2\le t$
with $t_2-t_1=u$ and $S^n_{ik}(t_2-)=S^n_{ik}(t_1)$, hence by the definition \eqref{004}
of $\hat S^n_{ik}$,
\[
\hat S^n_{ik}(t_1)-\hat S^n_{ik}(t_2-)=n^{-1/2}\mu^n_{ik}(t_2-t_1)\ge c_2n^{1/2}u,
\]
for some constant $c_2$ that depends only on the sequence $\{\mu^n_{ik}\}$.
Hence $Y^n_t\le 2c_2^{-1}n^{-1/2}\|\hat S^n_{ik}\|_t$.
The first result now follows from \eqref{250} with $\kappa=2$.

To prove the second statement we will proceed in two steps.
First, we will show that the number of interarrival times involved in the maximum is at most $cn$ with probability going to 1, then show that the maximum over $cn$ variables has the right order of magnitude.

To simplify the notation, fix $(i,k)$ and remove them from the notation
of $S^n_{ik}$, $u^n_{ik}$, $\hat\mu^n_{ik}$, etc. The claim will be proved for $e^n_{\rm max}(t)$
defined as in \eqref{eq:maxservice} but without maximizing over $(i,k)$;
clearly, this is sufficient.
Note that 
$$S^n(t)=\sup\left\lbrace s\geqslant 0: \, \sum_{l=1}^su^n(l)\leqslant t\right\rbrace.$$
Let 
$$K^n= \inf\left\lbrace s\geqslant 0 :  \sum_{l=1}^s u^n(l)\geqslant t_0\right\rbrace.$$
Then
\begin{align*}
	e^n_{\rm max}(t_0)&\leqslant  \sup\left\lbrace u^n(s): \sum_{l=1}^{s-1}u^n(l)\leqslant t_0\right\rbrace
	\leqslant \sup\left\lbrace u^n(s): \, s\leqslant K^n\right\rbrace.
\end{align*}
For any $c_2\geqslant 0$,
\begin{align*}
	\mathbb{P}\left( K^n\geqslant c_2n\right)&\leqslant
	\mathbb{P}\left(\sum_{p=1}^{c_2n} u^n(p)\leqslant t_0\right)\\
	&=\mathbb{P}\left(\sum_{p=1}^{c_2n} \check{u}(p)\leqslant \mu^nt_0\right)\\
	&=\mathbb{P}\left(\dfrac{1}{c_2n}\sum_{p=1}^{c_2n} \check{u}(p)\leqslant \dfrac{\mu t_0}{c_2}+\dfrac{\hat\mu^nt_0}{c_2}n^{-1/2}\right).
\end{align*}
If $c_2> \mu t_0$, by the law of large numbers, the RHS
converges to 0 as $n\to\iy$.
Next, by independence of service times, for any $c>0$,
\begin{align*}
	U_n:=
	\mathbb{P}\left(\max_{l\leqslant c_2n }u^n(l)\geqslant c n^{\bar a-1}\right)
	=1-\mathbb{P}\left(u^n(1)\leqslant  cn^{\bar a-1}\right)^{c_2n}.
\end{align*}
Using $1-(1-x)^n\leqslant c_3nx$, letting $c_4=c_2c_3$
and denoting $\bar\eps=\m-2>0$, we obtain
\begin{align*}
	U_n&\leqslant c_4n\mathbb{P}\left(u^n(1)\geqslant  cn^{\bar a-1}\right)\\
	&\leqslant c_4n\mathbb{P}\left(  \check{u}(1)\geqslant  c\mu^nn^{\bar a-1}\right)\\
	&= c_4n\mathbb{P}\left(  (\check{u}(1))^{2+\bar\eps}\geqslant  \left({c\mu^n}n^{\bar a-1}\right)^{2+\bar\eps}\right)\\
	&\leqslant c_4n\dfrac{\mathbb{E}\left[ (\check{u}(1))^{2+\bar\eps}\right]}{\left({c\mu^n} n^{\bar a-1}\right)^{2+\bar\eps}},
\end{align*}
where the last inequality uses Assumption \ref{as:polymo}. Next,
by the definition in \eqref{140} of $\bar a$ we have that
$\bar a >\frac{1}{2}-\frac{\bar\eps}{4(\bar\eps+2)}$, hence
\begin{align*}
	(2+\bar\eps)(1+\bar a-1)&=\bar a(2+\bar\eps)\\
	&\geqslant (2+\bar\eps)(\frac{1}{4}+\dfrac{1}{4+2\bar\eps})\\
	&=\frac{1}{2}+\frac{\bar\eps}{4}+\dfrac{2+\bar\eps}{4+2\bar\eps}=1+\frac{\bar\eps}{4}.
\end{align*}
Thus 
$$
n\dfrac{\mathbb{E}\left[ (\check{u}(1))^{2+\bar\eps}\right]}{\left({c\mu^n} n^{\bar a -1}\right)^{2+\bar\eps}}=O(n^{-\frac{\bar{\eps}}{4}}).
$$
Hence for any $c>0$ there exists $c_2$ such that
\begin{align*}
	\mathbb{P}\left(e_{\max}^n\geqslant cn^{\bar a-1}\right)&\leqslant \mathbb{P}\left(K^n\geqslant c_2n\right)+ \mathbb{P}\left(e_{\max}^n\geqslant cn^{\bar a-1},\,K^n\leqslant c_2n\right)\\
	&\leqslant \mathbb{P}\left(K^n\geqslant c_2n\right)+\mathbb{P}\left( \max_{l\leqslant c_2n }u^n(l)\geqslant  cn^{\bar a-1}\right),
\end{align*}
and both terms have been shown to converge to zero.
\qed

\subsubsection{Uniform integrability}\label{sec:ui}

Here we prove Proposition \ref{lem:m1}.
We sometimes need to refer to the variable keeping track of the current mode,
which in the various policies is defined in slightly different ways according
to different sampling times.
Denote
$$
\mathrm{MODE}^n(t)=\begin{cases}
	\xi^L \text{ if the current mode is lower workload,}&\\
	\xi^H\text{ if the current mode is higher workload.}
\end{cases}
$$

\noi{\bf Proof of Proposition \ref{lem:m1}.}
Fix any one of the sequences of policies $T^n\in\calA^n$ for which we
attempt to prove AO. The proof has three parts.

Part 1.
This part is concerned with the case where, whenever the workload in the system
is sufficiently large, policy \p{} is active.
This covers case 1(a) of Theorem \ref{th-ao-s},
where the system has a single mode and the fixed priority policy \p{} is applied,
as well as case 2(a) where the dual mode policy \p{}\p{} is applied.
In this case we will prove that the statement of the lemma holds with $\eps_0=1$,
and specifically that, under the given sequence,
$\E[(\hat H^n_t)^2]$ is bounded by a polynomial in $t$ for all $n$.
By \eqref{241}, this is equivalent to the same property
holding for $\E[(\hat W^n_t)^2]$.

To prove the result in this case,
assume that in the single (resp., dual) mode case,
the active mode $\xi^A$ (resp., the high workload mode $\xi^H$)
is in canonical form. Thus, provided that the workload in the system
exceeds $z^*$ (when applicable), class $1$ (resp., $2$)
is the dual (single) activity class, and server $1$ (resp., $2$) is the single (dual) activity server. 
In addition, $p=2$: class 2 is the HPC.
First we provide a bound on the second moment of $\hat X^n_2$. In the dual mode
case, fix $K$ large enough so that
$\hat X^n_2(t)\ge K$ implies $\hat W^n_t\ge z^*+1$; in the single mode case
let $K=1$.
Given $t>0$, consider the event $\hat X^n_2(t)>K$. Let $\tau=\tau_n(t)$ be defined by
\[
\tau=\sup\{s\in[0,t]:\hat X^n_2(s)\le K\}.
\]
Because the system starts empty,
$0\le\tau\le t$, and because the jumps of the normalized queue length are of size
$n^{-1/2}$, $\hat X^n_2(\tau)\le K+1$. Thus
\[
\hat X^n_2(t)\le \hat X^n_2(t)-\hat X^n_2(\tau)+K+1.
\]
By \eqref{eq:balancex}, denoting $C^n(t)=\sum_i\|\hat A^n_i\|_t+\sum_{ik}\|\hat S^n_{ik}\|_t$,
\[
\hat X^n_2(t)-\hat X^n_2(\tau)\le 2C^n(t)+n^{-1/2}\la^n_2(t-\tau)
-n^{-1/2}\mu^n_{22}(T^n_{22}(t)-T^n_{22}(\tau)).
\]
Because $\hat X^n_2\ge K$ in $[\tau,t]$, the priority policy corresponding to
a mode given in canonical form (either $\xi^A$ or $\xi^H$) is in force after an initial time before the current mode updates and there are class 2 jobs to serve.
We can bound the time until server $2$
prioritizes class 2 and starts serving it at full rate in $[\tau,t]$ by $e^n_{\rm max}(t)$. If there currently is a class 2 job being served by server 2, server 2 only serves class 2 jobs for the whole period and 
$$T^n_{22}(t)-T^n_{22}(\tau)= t-\tau.$$
If there currently is a class 1 job being served by server 2, because service is non-preemptive, server 2 has to finish this job.
If at that time the current mode is $\xi^H$ (resp. $\xi^A$ in the single mode case), server 2 prioritizes class 2 from that point on. If at that time the current mode  is $\xi^L$ the workload is then sampled because a service finished at the single activity server. The current mode  changes to $\xi^H$ and stays until $t$. In both cases we get
$$T^n_{22}(t)-T^n_{22}(\tau)\ge t-\tau-e^n_{\max}.$$

Also, in the case under consideration one has $\mu_{22}>\la_2$.
Hence for all sufficiently large $n$,
$\mu^n_{22}>\la^n_{22}$. Moreover,
$c_1:=\sup_nn^{-1}\mu^n_{22}<\iy$. Hence
\begin{align*}
	\hat X^n_2(t)
	\le 2C^n(t)+c_1n^{1/2}e^n_{\rm max}(t)+K+1.
\end{align*}
The above inequality holds also on the complementary event, namely when
$\hat X^n_2(t)\le K$. Thus using \eqref{250} and Lemma \ref{lem:1service} we obtain
\begin{equation}\label{251}
	\E[\|\hat X^n_2\|_t^2]\le c(1+t).
\end{equation}

In the next step we bound $\hat W^n_t$.
Let $\tilde K$ be a constant that is sufficiently large
to ensure that $\hat X^n_1(t)\ge\tilde K$ implies $\hat W^n(t)\ge z^*+1$ (where again,
$\tilde K=1$ in the single mode case).
Given $t>0$, consider the event $\hat X^n_1(t)> \tilde K$ and let $\sig=\sig_n(t)$ be
\[
\sig=\sup\{s\in[0,t]:\hat X^n_1(s)\le \tilde K\}.
\]
Then clearly $\sig\in[0,t]$.
Because during $[\sig,t]$ there are at least
two jobs of class $1$ in the system, both servers are never idle during $[\sig,t]$.
Since one  has $\hat X^n_1(\sig)\le \tilde K+n^{-1/2}$,
it follows that $\hat W^n_\sig\le c(1+\tilde K+\hat X^n_2(\sig))$.
Hence
\[
\hat W^n_t\le\hat W^n_t-\hat W^n_\sig+c(1+\tilde K+\|\hat X^n_2\|_t).
\]
By \eqref{242} and the nonidling of both servers during $[\sig,t]$,
by which $\hat L^n$ remains flat during this interval,
we have $\hat W^n_t-\hat W^n_\sig=\hat F^n_t-\hat F^n_\sig$.
Thus by \eqref{096}, for a constant $c$ (that may depend on $\tilde K$),
\[
\hat W^n_t\le c(C^n(t)+t+1+\|\hat X^n_2\|_t).
\]
Clearly, the above bound is valid also in the complementary event,
$\hat X^n_1(t)\le \tilde K$. We can therefore apply
\eqref{250} and \eqref{251},
to obtain $\E[(\hat W^n_t)^2]\le c(1+t)$ for some constant
$c$, for all $n$ and $t$. Consequently the same holds for the second moment
of $\hat H^n_t$, and the result follows.

Part 2.
Next consider the case where \teeone{} is applicable when
the workload is sufficiently large.
This covers case 2(d) of Theorem \ref{th-ao-s}, in which
the policy \teetwo{}\teeone{} applies (none of our proposed
policies implement \teeone{} as a single mode policy).
The proof given in Part 1 is applicable, for the following reasons.
In the first step, during the analyzed time interval $[\tau,t]$, one has
$\hat X^2\ge K>\hat\Theta^n$, and therefore there is no difference
between how \p{} and \teeone{} behave during this interval. In addition, the workload is sampled at each arrival/service, which means that 
$$T^n_{22}(t)-T^n_{22}(\tau)\ge t-\tau-e^n_{\rm max}.$$
In the second step, the argument given for nonidling of both
servers during $[\sig,t]$ again holds here similarly to Part 1, upon noticing
that, since $\sigma$ corresponds to an arrival, it is a sampling time,
and so the current mode is either
already the high mode or switches to the high mode at that time.
(It is possible that, if $\sigma$ is a mode switching time, then server 1
was idle just before $\sigma$. But, if so, it starts to serve class 1 at $\sigma$.)
The remaining details need no adaptation.

Part 3.
Finally, consider the policies which employ \teetwo{} for high workload levels,
namely the policies \teetwo{}, \teetwo{}\teetwo{} and \teeone{}\teetwo{},
covering all remaining cases of Theorem \ref{th-ao-s}.
In these cases the stronger moment assumption is in force,
and the goal is to prove that there exists $\eps_0>0$ such that
$\E[(\hat W^n_t)^{1+\eps_0}]\le {\rm pol}(t)$ for all $n$ and $t$, for some
polynomial ${\rm pol}$.
As before, let $\xi^A$ or $\xi^H$ be in canonical form,
in the single and, respectively, dual mode case.
Fix a constant $K>z^*$ in the dual mode case
and $K=1$ in the single mode case. Given $t$ consider the event
$\hat W^n_t>K$. Let
\[
\tau_1=\tau_1^n=\sup\{s\in[0,t]:\hat W^n_s\le K\}.
\]
Then $\hat W^n_{\tau_1}\le K+1$, and moreover $\hat W^n\ge K$ during
$[\tau_1,t]$.
Although this lower bound on the workload is sufficiently large to guarantee that
$\hat W^n$ corresponds to the higher workload mode, it is possible that at time
$\tau_1$ the current mode variable still equals the lower workload mode. We argue that the
time it takes to switch to the upper workload mode is bounded by $2e^n_{\rm max}$
in both \teeone{}\teetwo{} and \teetwo{}\teetwo{}. In the former case,
the current mode switches as soon as a there is a new arrival or departure.
In the latter case, one possibly has to complete the service of a job at server 2,
which is server $k_1(\xi^L)$. It is possible that there are no jobs allowed to be routed to server 2 at time $\tau_1$. Wait for an arrival of class 1 job, and service completion
of this job at server 2. At this time it is guaranteed that mode has switched to $\xi^H$
if it was not $\xi^H$ earlier.
Thus if we let
\[
\tau_2=\inf\{s\ge \tau_1:{\rm MODE}^n(s)=\xi^H\}\w t,
\]
we have $\tau_2-\tau_1\le2e^n_{\rm max}\w t$.

According to the rules of \teetwo{},
server 2 must be busy throughout the interval $[\tau_2,t]$.
Thus
$\hat I^n_2(t)=\hat I^n_2(\tau_2)$. Hence by \eqref{242},
and using $\hat I_k^n[a,b]\le n^{1/2}(b-a)$ (by \eqref{004+}), we have
\begin{align*}
	\hat W^n_t&= \hat W^n_{\tau_1}+\hat F^n[\tau_1,t]+\hat L^n[\tau_1,t]
	\\
	&=  \hat W^n_{\tau_1}+\hat F^n[\tau_1,t]
	+\beta_1\hat I^n_1[\tau_1,t]+\beta_2\hat I^n_2[\tau_1,t]
	\\
	&\le K+1+2\|\hat F^n\|_t+cn^{1/2}e^n_{\rm max}+\beta_1\hat I^n_1[\tau_2,t]
\end{align*}
holds on the event $\hat W^n_t>K$. On the complementary event,
$\hat W^n_t\le K$. Use \eqref{096} and \eqref{250}
to obtain the bound
$\E[\|\hat F^n\|_t^{1+\eps_0}]\le\{\E[\|\hat F^n\|_t^2]\}^{(1+\eps_0)/2}
\le c(1+t)^{(1+\eps_0)/2}$. Use Lemma \ref{lem:1service} to bound the second moment of $n^{1/2}e^n_{\rm max}$
by $c(1+t)$. By Minkowski's inequality
and the inequality $(a+b)^{1+\eps_0}\le4a^{1+\eps_0}+4b^{1+\eps_0}$
which holds for $a,b\ge0$, $\eps_0\in(0,1)$, this yields
\[
\E[(\hat W^n_t)^{1+\eps_0}]
\le c(1+t)^{(1+\eps_0)/2}
+c\E[\one_{\{\hat W^n_t>K\}}\hat I^n_1[\tau_2,t]^{1+\eps_0}].
\]
Hence it suffices to bound the last term above by a polynomial in $t$.

To this end, let $\tau_3=\inf\{s\ge\tau_2:\hat X^n_1(s)\ge\hat\Th^n\}\w t$.
Then
\begin{align*}
	\E[\one_{\{\hat W^n_t>K\}}\hat I^n_1[\tau_2,t]^{1+\eps_0}]
	&\le 4\Del^n_1+4\Del^n_2, \quad\text{where}
	\\
	\Del^n_1&=\E[\one_{\{\hat W^n_t>K\}}\hat I^n_1[\tau_2,\tau_3]^{1+\eps_0}]
	\\
	\Del^n_2&=\E[\one_{\{\hat W^n_t>K,\tau_3<t\}}\hat I^n_1[\tau_3,t]^{1+\eps_0}].
\end{align*}
To bound $\Del^n_1$, consider the event $\{\hat W^n_t>K,\tau_3>\tau_2\}$.
On it, during the interval $[\tau_2,\tau_3]$,
one has $\hat X^n_1<\hat\Th^n$, hence
by the rules of \teetwo, server 2 prioritizes class 2, except possibly it completes
a service that started when $\hat{X}^n_1\geqslant \hat{\Th}^n$.
Moreover, $\hat W^n\ge K$ during the same interval,
by which we know that there are multiple class-2 jobs in the system,
and thus server 2 gives no service to class 1, with the only exception of service
to a job that started when $\hat{X}^n_1\geqslant \hat{\Th}^n$.
Hence the departure process associated
with activity $(1,2)$ increases by at most $1$ during this interval, that is,
$0\le D^n_{12}[\tau_2,\tau_3]\le 1$. Recalling the definition
of $\hat S^n$ in \eqref{004}, this can be expressed as
$0\le e^n_1[\tau_2,s]\le n^{-1/2}$, $s\in[\tau_2,\tau_3]$, where
\[
e^n_1(s):=n^{-1/2}D^n_{12}(s)=\hat S^n_{12}(T^n_{12}(s))+n^{-1/2}\mu^n_{12}T^n_{12}(s).
\]
Using this in \eqref{eq:balancexx} gives, for $s\in[\tau_2,\tau_3]$,
\[
\hat X^n_1(s)=\hat X^n_1(\tau_2)+\hat f^n_1[\tau_2,s]-e^n_1[\tau_2,s]
+n^{1/2}\la_1(s-\tau_2)-n^{1/2}\mu_{11}T^n_{11}[\tau_2,s],
\]
where
\[
\hat f^n_1(s)=\hat A^n_1(s)-\hat S^n_{11}(T^n_{11}(s))
+(\hat\la^n_1s-\hat\mu^n_{11}T^n_{11}(s)).
\]
Next, by \eqref{41},
\[
T^n_{11}[\tau_2,s]=(s-\tau_2)-I^n_1[\tau_2,s]-T^n_{21}[\tau_2,s].
\]
However, activity $(2,1)$ is not in use by \teetwo{}.
Denoting $c_1=\la_1-\mu_{11}$ and $g(s)=c_1s$, this yields
\begin{equation}\label{300}
	\hat X^n_1(s)=\hat X^n_1(\tau_2)+\hat f^n_1[\tau_2,s]-e^n_1[\tau_2,s]
	+n^{1/2}g[\tau_2,s]+\mu_{11}\hat I^n_1[\tau_2,s],
	\qquad s\in[\tau_2,\tau_3].
\end{equation}
Because \teetwo{} is used after the update of modes, it must be true that $c_1=\la_1-\mu_{11}>0$.

We use the following property of the Skorohod map. Let $\eta=\Gam_2[\psi]$.
Then by \eqref{e02},
$\eta(s)=\sup_{0\le \theta \le s}[\psi(\theta)^-]$.
Assume that $\psi=\psi_1+\psi_2$ where $\psi_2\ge0$. Then
\begin{align*}
	\eta(s)&=\Gam_2[\psi_1+\psi_2](s)=\sup_{\theta\in[0,s]}[(\psi_1(\theta)+\psi_2(\theta))^-]
	\\
	&\le\sup_{\theta\in[0,s]}[\psi_1(\theta)^-]
	\\
	&\le \|\psi_1\|_s.
\end{align*}
This property is used as follows. On the interval $[\tau_2,\tau_3]$,
$\hat I^n_1$ can increase only at times when $\hat X^n_1=0$; and the latter
process is nonnegative.
This shows that $\mu_{11}\hat I^n_1$ serves as the Skorohod term in \eqref{300}
on that interval, and thus by the non-negativity of $\hat X^n_1(\tau_2)$
and $c_1$, and the bound $|e^n_1[\tau_2,s]|\le n^{-1/2}$,
this shows that
\[
\mu_{11}\hat I^n_1[\tau_2,\tau_3]
\le\sup_{\theta\in[\tau_2,\tau_3]}|\hat f^n_1[\tau_2,\theta]|
+\sup_{\theta\in[\tau_2,\tau_3]}|e^n_1[\tau_2,\theta]|
\le 2\|\hat f^n_1\|_t+n^{-1/2}.
\]
Hence $\Del^n_1\le c\{\E[\|\hat f^n_1\|_t^2]\}^{(1+\eps_0)/2}+1\le c(1+t)$.

Next we bound $\Del^n_2$. To this end, consider the event
$\Om^n:=\{\hat W^n_t>K,\tau_3<t,\hat I^n_1(t)>\hat I^n_1(\tau_3)\}$.
Because by its definition, $\hat I^n_1(t)\le n^{1/2}t$, we have
\[
\Del^n_2\le\PP(\Om^n)(n^{1/2}t)^{1+\eps_0}.
\]
On the event $\Om^n$ let
\[
\tau_5=\inf\{s>\tau_3:\hat X^n_1(s)=0\},
\qquad
\tau_4=\sup\{s<\tau_5:\hat X^n_1(s)\ge\hat\Th^n\}.
\]
Then on $\Om^n$ it must hold that $\tau_3\le\tau_4<\tau_5\le t$.
The arguments which lead to \eqref{300} are valid
for the time interval $[\tau_4,\tau_5]$. As a result,
\[
\hat X^n_1[\tau_4,\tau_5]=\hat f^n_1[\tau_4,\tau_5]
-e^n_1[\tau_4,\tau_5]+n^{1/2}g[\tau_4,\tau_5]
+\mu_{11}\hat I^n_1[\tau_4,\tau_5],
\]
with $|e^n_1[\tau_4,s]|\le n^{-1/2}$ for all $s\in[\tau_4,\tau_5]$.
Now, $\hat I^n_1$ remains flat on the interval $[\tau_4,\tau_5]$, and moreover,
$\hat X^n_1(\tau_4)=\hat\Th^n-n^{-1/2}$ and $\hat X^n_1(\tau_5)=0$.
This gives
\[
\hat f^n_1[\tau_4,\tau_5]
+n^{1/2}g[\tau_4,\tau_5]\le-\hat\Th^n+2n^{-1/2}.
\]
Given any $\del>0$, using the nonnegativity of the second term on the LHS,
in the case that $\tau_5-\tau_4\le n^{-\del}$ one must have
$\hat f^n_1[\tau_4,\tau_5]\le-\hat\Th^n/2$. On the other hand,
in the case $\tau_5-\tau_4>n^{-\del}$ one must have
$\hat f^n_1[\tau_4,\tau_5]+c_1n^{1/2}n^{-\del}<0$.
As a result,
\[
\PP(\Om^n)\le p^n_1+p^n_2
:=\PP(w_t(\hat f^n_1,n^{-\del})\ge\hat\Th^n/2)+\PP(2\|\hat f^n_1\|_t>c_1n^{\frac{1}{2}-\del}).
\]
We have by \eqref{141} $\hat\Th^n=n^{-\frac{1}{2}}\lceil n^{\bar a} \rceil $,
where we recall that $\bar a<\frac{1}{2}$.
Thus by Lemma \ref{lem10}, with $\nu_1=\del$, $\nu_2=\frac{1}{2}-\bar a$,
one has, for any $\eps_1>0$,
\[
p^n_1\le c(1+t)^{\frac{\m}{2}}n^{-h_0+\eps_1}
\]
where
\[
h_0=(\frac{\m}{2}-1)\del-\m(\frac{1}{2}-\bar a),
\qquad
\del\le 1-\bar a.
\]
Next, by \eqref{250} and Chebychev's inequality,
\[
p^n_2\le c(1+t)^{\frac{\m}{2}}n^{-\frac{\m}{2}+\del \m}.
\]
As a result, we have
\[
\Del_2^n\le c(n^{\zeta_1(\del,\eps_0,\eps_1)}
+n^{\zeta_2(\del,\eps_0)})
(1+t)^{\frac{\m}{2}+1+\eps_0},
\]
where
\[
\zeta_1(\del,\eps_0,\eps_1)=
-\Big(\frac{\m}{2}-1\Big)\del+\m\Big(\frac{1}{2}-\bar a\Big)+\frac{1}{2}+\frac{\eps_0}{2}+\eps_1,
\qquad
\zeta_2(\del,\eps_0)=
-\frac{\m}{2}+\del \m+\frac{1}{2}+\frac{\eps_0}{2}.
\]
By our assumptions, we have $\m>\m_0=\frac{1}{2}(5+\sqrt{17})$
and $\bar a>\frac{1}{2}-\frac{\m^2-5\m+2}{2\m(3\m-2)}$.
Using this, a calculation shows that with the choice $\del=\m/(3\m-2)
\in(1/3,1/2)$, one has
$\zeta_1(\del,0,0)\vee\zeta_2(\del,0)<0$.
It follows that there exist $\eps_0>0$ and $\eps_1>0$
for which $\zeta_1(\del,\eps_0,\eps_1)\vee\zeta_2(\del,\eps_0)<0$.
Therefore
$\Del^n_2\le c(1+t)^{\frac{\m}{2}+1+\eps_0}$ and the proof is complete.
\qed

\subsubsection{State space collapse}\label{sec:ssc}

We now prove Proposition \ref{lem:ssc}.
Assume that either the active mode $\xi^A$ or the lower workload mode $\xi^L$ (whichever is applicable) is in canonical form. 
Let
\[\tau=\tau^n_c= \inf\left\lbrace t\geqslant 0: \, \hat{X}^n_p(t)\geqslant 2\hat\Th^n\right\rbrace.\]
This random time is used outside this proof with the notation
$\tau^n_c$, but in this proof the shorter notation $\tau$ is used.
Then
\[\mathbb{P}\left(\sup_{t\leqslant t_0}\hat{X}^n_p(t)\geqslant 2\hat\Th^n\right)\le\mathbb{P}\left(\tau\leqslant t_0\right).\]
We will prove the lemma by showing that the RHS above converges
to zero as $n\to\iy$.
On the event $\left\lbrace\tau\leqslant t_0\right\rbrace$, define
\[\sig=\sigma^n=\sup\left\lbrace t\leqslant \tau: \, \hat{X}^n_p(t)\leqslant \dfrac{3\hat\Th^n}{2}\right\rbrace.\] 
The proof relies on the fact that, under our policies, when the number of HPC jobs is above $\Th^n$, it is served at a rate that is enough to deplete the queue in all cases. 
The first step toward this goal is the following.
\begin{lemma}\label{lem:servicemax}
	There exist constants $c_1,c_2>0$ such that,
	on $\left\lbrace\tau\leqslant t_0\right\rbrace$,
	\begin{equation}
		\int_{\sigma}^{\tau}\big(\lambda_p-\sum_k\mu_{pk}\X^n_{pk}(t)\big)dt\leqslant  -c_1(\tau-\sigma)+c_2e_{\max}^n,\label{eq:servicemax}
	\end{equation}
	with $e_{\max}^n$ defined in \eqref{eq:maxservice}.	
\end{lemma}

\noi{\bf Proof.}
First, by definition of $\sigma$, and $\tau$, and the fact that jumps of $\hat{X}^n$ are of size $n^{-1/2}<\hat{\Th}^n/2$ for large $n$, we obtain
\[\inf\limits_{t\in [\sigma, \tau]}\hat{X}^n_p(t)\geqslant  \hat\Th^n.\]
We will now distinguish along the possible cases.
\begin{itemize}
	\item \textbf{Case 1(a)}: In this case, we use the \p{} policy, which is single mode. Because the active mode is in canonical form and $i_1(\xi^A)=p$, $p=2$ and server 2 prioritizes class 2. It is possible that server 2 is busy with the ``wrong'' class of job at time $\sigma$ but as soon as that job finishes, server 2 will only serve the class 2 jobs in the system:
	$$\int_{\sigma}^{\tau}\X^n_{22}(t)dt\geqslant \tau-\sigma-e_{\max}^n.$$
	Thus,
	\begin{equation*}
		\int_{\sigma}^{\tau}\big(\lambda_2-\sum_k\mu_{2k}\X^n_{2k}(t)\big)dt\leqslant (\tau-\sigma)(\lambda_2- \mu_{22})+ \mu_{22}e_{\max}^n.
	\end{equation*}
	To show that $\lambda_2- \mu_{22}<0$, note that since the active mode is in canonical form, $\frac{\lambda_1}{\alpha_1}>\beta_1$, which implies $\frac{\lambda_2}{\alpha_2}<\beta_2$ by \eqref{eq:ht}.
	\item\textbf{ Case 1(b)}: In this case, we use the \teetwo{} policy which is single mode. Between $\sigma$ and $\tau$, $\hat{X}^n_p(t)\geqslant \hat\Th^n$. Thus, both servers give priority to the HPC at all time in $[\sigma, \tau]$. It is possible that the dual activity server is occupied with the low priority class at $\sigma$ but as soon as the current job is served, the high priority class gets priority on one server and dedication by the other server. By definition of $e_{\max}^n$, we have for any $k\in \lbrace1,2\rbrace$, 
	$$\int_{\sigma}^{\tau}\X^n_{pk}(t)dt\geqslant  \tau-\sigma-e_{\max}^n\mathds{1}_{k=k_2(\xi^A)}.$$
	Thus
	\begin{equation*}
		\int_{\sigma}^{\tau}\big(\lambda_p-\sum_k\mu_{pk}\X^n_{pk}(t)\big)dt\leqslant (\tau-\sigma)(\lambda_p-\sum_k \mu_{pk})+\mu_{pk_2(\xi^A)}e_{\max}^n.
	\end{equation*}
	Finally, $\lambda_p-\sum_k \mu_{pk}<0$  by \eqref{eq:ht} and $\min_i\lambda_i>0$.
	
	\item \textbf{Case 2(a)}: In this case, we use the \p{}\p{} policy, which  only changes mode at the completion of a service at the single activity server. This is precisely the server that gives priority to the high priority class after switching mode because we are in a \textbf{SS} case. Since there are always HPC jobs to serve between $\sigma$ and $\tau$,  we claim that excluding an initial period, at any given time the HPC is being served by at least one server, as discussed in Section \ref{sec:disc}. We now discuss how long this initial period can be.
	
	The only way some service is lost is if there were no high priority jobs at some point in the system, both servers become busy with low priority jobs and $\sigma$ occurs before the service completion of those jobs. After completion of those two jobs, the HPC keeps priority on at least one server regardless of switching of the current mode. For the initial jobs, with respect to $\mathrm{MODE}^n(\sigma)$ either the service ends first at the dual activity server or the service ends first at the single activity server.
	
	In the first case, a service of the HPC job starts at the dual activity server because it has priority there until  either $\tau$ is reached or a mode switch occurs.
	In the second case the available server is currently dedicated to the LPC.  Since this corresponds to a service  completion at the single activity server, the workload is sampled and there could also be a switching of modes. 
	If the mode changes then this server becomes dual activity and the HPC has priority there. If there is no mode switch, the service of another LPC job starts.  LPC jobs are served at this server until  the mode switches. 
	
	Thus, if the   service ends  at the dual activity server before the mode switches then an HPC job will start there. If the mode switches before the   service ends  at the current dual activity server, then the current single activity server becomes dual activity and an HPC job begins service there.
	
	In either case the the HPC is served by at least one server after that time because it has priority on the dual activity server, and the dual activity server is always available when switching modes.
	In addition, the time it takes for the HPC to begin service is smaller than the service of the job present at the dual activity server at time $\sigma$, which is smaller than $e_{\max}^n$. Hence
	$$\int_{\sigma}^{\tau}\left(\mu_{p1}\X^n_{p1}(t)+\mu_{p2}\X^n_{p2}(t)\right)dt\geqslant \min_k\mu_{pk}(\tau-\sigma)-\max_k \mu_{pk}e_{\max}^n.$$ 
	This yields
	\begin{equation*}
		\int_{\sigma}^{\tau}\big(\lambda_p-\sum_k\mu_{pk}\X^n_{pk}(t)\big)dt\leqslant (\tau-\sigma)(\lambda_p- \min_k\mu_{pk})+\max_k \mu_{pk}e_{\max}^n.
	\end{equation*} 
	In addition, $\lambda_p- \min_k\mu_{pk}<0$. Putting $\xi^L$ in canonical form forces in this case $\frac{\lambda_1}{\al_1}>\beta_1$ and  $p=2$ because $i_1(\xi^L)=p$. In addition, by \eqref{92}, either
	$\frac{\lambda_1}{\alpha_1}>\beta_1 \vee \beta_2$
	or $\frac{\lambda_2}{\alpha_2}> \beta_1\vee \beta_2$.
	Thus  $\frac{\lambda_1}{\alpha_1}>\max_k\beta_k$, which implies $\frac{\lambda_2}{\al_2}<\min_k\beta_k$ by \eqref{eq:ht}. 
	
	\item \textbf{Case 2(b)}: In this case, we use the \teetwo{}\teetwo{} policy. Since the number of HPC jobs stays above $\Th^n$ throughout the period $[\sigma,\tau]$, both servers only take new jobs from the HPC regardless of the mode. Similarly to 1(b), there is at most one job served at either activity before the HPC gets served. Hence, for any $k\in \lbrace1,2\rbrace$, 
	$$\int_{\sigma}^{\tau}\X^n_{pk}(t)dt\geqslant  \tau-\sigma-e_{\max}^n.$$
	Thus
	\begin{equation*}
		\int_{\sigma}^{\tau}\big(\lambda_p-\sum_k\mu_{pk}\X^n_{pk}(t)\big)dt\leqslant (\tau-\sigma)(\lambda_p-\sum_k \mu_{pk})+\sum_k \mu_{pk}e_{\max}^n.
	\end{equation*}
	We obtain the result similarly as before because $\lambda_p-\sum_k \mu_{pk}<0$.

	\item \textbf{Cases 2(c) and 2(d)}: the HPC is single activity in one mode and dual activity in the other. Recall that this case is \textbf{CS}. This means the dual activity server stays the same regardless of switching modes. When a \teeone{} rule is applied, the dual activity server prioritizes the single activity class as long as there are more than $\Th^n$ jobs of this class.  The HPC is single activity when we apply the \teeone{} rule. Similarly, under a \teetwo{} rule, the dual activity server prioritizes the dual activity class as long as there are more than $\Th^n$ jobs of this class and the HPC is dual activity when we apply the \teetwo{} rule. Once again put $\xi^L$ in canonical form, server 2 is the dual activity server in both modes. Regardless of which rule is used, as long as the number of HPC jobs is above $\Th^n$, server 2 only takes new jobs from the HPC. As before, we have to exclude a residual service of a low priority job that started before $\sigma$. Hence
	$$\int_{\sigma}^{\tau}\X^n_{p2}(t)dt\geqslant  \tau-\sigma-e_{\max}^n.$$ 
	It remains to show that the activity processing HPC jobs throughout the period is enough to deplete them. Because of the canonical form of $\xi^L$, $\frac{\lambda_1}{\al_1}>\beta_1$ and thus by \eqref{93}, $\frac{\lambda_1}{\al_1}<\beta_2$. Moreover, $\frac{\lambda_2}{\al_2}<\beta_2$ by \eqref{eq:ht}. Whichever the HPC is, server 2 has enough capacity to deplete it. In other words, 
	$\lambda_p-\mu_{p2}<0$ and 
	\begin{equation*}
		\int_{\sigma}^{\tau}\big(\lambda_p-\sum_k\mu_{pk}\X^n_{pk}(t)\big)dt\leqslant (\tau-\sigma)(\lambda_p-\mu_{p2})+\mu_{p2}e_{\max}^n.
	\end{equation*}
\end{itemize}
\qed

Now that Lemma \ref{lem:servicemax} has been proved in all cases,
the proof of Proposition \ref{lem:ssc} does not need to differentiate between
them.

\noi{\bf Proof of Proposition \ref{lem:ssc}.}
Let $\nu_2=1/2-\bar a\leqslant 1/4$ and $\nu_1\in (\nu_2 , \,1/2 )$. Recall that 
$ \hat \Th^n=n^{-1/2}\lceil n^{\bar a} \rceil$, as defined in \eqref{141}, so that $ \hat \Th^n=n^{-1/2}\lceil n^{1/2-\nu_2} \rceil \geqslant  n^{-\nu_2}$.
Notice that, as required in Lemma \ref{lem10}, $\nu_1\leqslant \nu_2+\frac{1}{2}$. Let us introduce the event
$$\Omega_1= \left\lbrace \tau\leqslant t_0,\,  \hat{X}^n_p(\tau)-\hat{X}^n_p(\sigma)\geqslant \dfrac{\hat{\Th}^n}{4}, \inf\limits_{t\in [\sigma, \tau]}\hat{X}^n_p(t)\geqslant \hat{\Th}^n\right\rbrace.$$
For $n$ large enough one has
$\hat X^n_p(\sig)<\frac{7}{4}\hat\Th^n$.
Consequently $\{\tau\le t_0\}=\Omega_1$.
We obtain
\begin{equation}\label{eq:sep}
	\mathbb{P}\left(\tau\leqslant t_0\right)
	= \mathbb{P}\left(\Omega_1\right)\\
	= \mathbb{P}\left(\Omega_1\cap\{\tau-\sigma\leqslant n^{-\nu_1}\}\right)
	+\mathbb{P}\left(\Omega_1\cap\{\tau-\sigma> n^{-\nu_1}\}\right).
\end{equation}
Picking up on \eqref{eq:balancexxx}, by Lemma \ref{lem:servicemax}
\begin{align*}
	\hat{X}^n_p(\tau)-\hat{X}^n_p(\sigma)&=\hat{f}_p^n(\tau)-\hat{f}_p^n(\sigma)+\sqrt{n}\int_{\sigma}^{\tau}\big(\lambda_p-\sum_k\mu_{pk}\X^n_{pk}(t)\big)dt\\
	&	\leqslant \hat{f}_p^n(\tau)-\hat{f}_p^n(\sigma)-c\sqrt{n}(\tau-\sigma)+c\sqrt{n}e_{\max}^n.
\end{align*}
Notice that on the event $\left\lbrace \tau-\sigma\leqslant n^{-\nu_1},\, \tau\leqslant t_0\right\rbrace $, 
\begin{equation*}
	\hat{X}^n_p(\tau)-\hat{X}^n_p(\sigma)	\leqslant w_{t_0}(\hat{f}_p^n,n^{-\nu_1})+c\sqrt{n}e_{\max}^n.
\end{equation*}
Thus
\begin{align*} \mathbb{P}\left(\Omega_1\cap\{\tau-\sigma\leqslant n^{-\nu_1}\}\right)&\leqslant \mathbb{P}\left(w_{t_0}(\hat{f}_p^n,n^{-\nu_1}) +e_{\max}^n\sqrt{n}\geqslant \frac{n^{-\nu_2}}{2}\right)\label{eq:modul} \\
	&\leqslant  \mathbb{P}\left(w_{t_0}(\hat{f}_p^n,n^{-\nu_1})\geqslant\frac{n^{-\nu_2}}{4}\right)+ \mathbb{P}\left(ce_{\max}^n\sqrt{n}\geqslant \frac{n^{-\nu_2}}{4}\right).\end{align*}
Recall that $\nu_1\in (\nu_2,1/2)$, so $n^{-\nu_1}< n^{-\nu_2}$. By Lemma \ref{lem10} and Remark \ref{rem:stillok}, the first term converges to zero. By lemma \ref{lem:1service} and $ \nu_2=1/2-\bar a$,
\begin{equation*}
	\mathbb{P}\left(e_{\max}^n\geqslant \frac{c}{4}n^{-\nu_2-1/2}\right)=\mathbb{P}\left(e_{\max}^n\geqslant \frac{c}{4}n^{\bar a-1}\right)\to 0.
\end{equation*}

For the second term in \eqref{eq:sep}, notice that on $\left\lbrace \tau-\sigma> n^{-\nu_1},\, \tau \leqslant t_0\right\rbrace$, 
\begin{equation*}
	\hat{X}^n_p(\tau)-\hat{X}^n_p(\sigma)	\leqslant 2\lVert \hat{f}_p^n(t)\rVert _{t_0}+e_{\max}^n\sqrt{n}-c\sqrt{n}n^{-\nu_1},
\end{equation*}
and 
\begin{equation*}
	\hat{X}^n_p(\tau)-\hat{X}^n_p(\sigma)\geqslant  \dfrac{\hat{\Th}^n}{2}-n^{-1/2}.
\end{equation*}
Hence
\begin{equation*}
	\mathbb{P}\left(\Omega_1\cap\{\tau-\sigma> n^{-\nu_1}\}\right)\leqslant \mathbb{P}\left(2\lVert \hat{f}_p^n\rVert _{t_0} +e_{\max}^n\sqrt{n}\geqslant c\sqrt{n}n^{-\nu_1} +\dfrac{\hat{\Th}^n}{2}-n^{-1/2}\right).
\end{equation*}
By Lemma \ref{lem:1service}, $e_{\max}^n$ is smaller than $n^{\bar a-1}=o(n^{-1/2})$. By Remark \ref{rem:stillok}, $2\lVert \hat{f}_p^n\rVert _{t_0}$ is a tight sequence of RVs (for $t_0$ fixed). Because $\nu_1<1/2$, $\sqrt{n}n^{-\nu_1} \to +\infty$.
The claim follows.
\qed

\subsubsection{Boundary behavior}
\label{sec:boundarybeh}

The goal of this section and the following one
is to prove Proposition \ref{lem:reflection}.
The key is the follows lemma, which states, roughly speaking,
that $\hat L^n$ is approximately the boundary term
corresponding to $\hat F^n$.

Let $c_3= \frac{3}{\al_1 \wedge \al_2}$.

\begin{lemma}\label{lem:reflectionprelim}
	Fix $t_0>0$. 
	With
	\begin{equation}\label{eq:rbar}
		\bar R^n_t=\int_0^t
		\mathds{1}_{\hat{W}^n_s\geqslant c_3 \hat\Th^n }d\hat L^n_s,
	\end{equation}
	$\bar R^n_{t_0}\to0$ in probability
	as $n\to\iy$.
\end{lemma}

This lemma is proved in the next subsection.
Let us show how Proposition \ref{lem:reflection} now follows.

\noi{\bf Proof of Proposition \ref{lem:reflection}.}
In addition to $\bar R^n$ just introduced, we shall need
the following definitions:
\begin{align*}
	Z^n_t&=\max(\hat{W}^n_t-c_3\hat{\Th}^n,\, 0), \\
	\wt{R}^n_t&=\int_{0}^{t}\mathds{1}_{\hat{W}^n_s< c_3\hat{\Th}^n}d\hat L^n_s, \\
	\Del^n_t& = Z^n_t - \hat{W}^n_t.
\end{align*}
By \eqref{242} and definition of the new processes, one has
\begin{equation*}
	Z^n_t=\hat F^n_t+\Del^n_t+\bar{R}^n_t+\wt{R}^n_t.
\end{equation*}
Consider the pair $(Z^n,\wt R^n)$.
The first component is nonnegative. The second is
nonnegative, nondecreasing, and moreover,
\[
\int_{[0,\iy)}Z^n_td\wt{R}^n_t=\int_{[0,\iy)}
\max(\hat{W}^n_t-c_3\hat{\Th}^n,\, 0)\mathds{1}_{\hat{W}^n_t
	< c_3\hat\Th^n}d\hat{L}^n_t=0.
\]
Hence by Skorohod's Lemma,
$(Z^n,\wt{R}^n)=\Gam(\hat F^n+\Del^n+\bar{R}^n)$.
Hence
\begin{equation}\label{303}
	\hat W^n=\Gam_1(\hat  F^n+\Del^n+\bar R^n)-\Del^n,
	\qquad
	\hat L^n=\Gam_2(\hat  F^n+\Del^n+\bar R^n)+\bar R^n.
\end{equation}
Recall that we are considering a subsequence along which
(according to the assumptions and to Lemma \ref{lem4}),
$(\hat A^n,\hat S^n,T^n,\hat F^n)\To(A,S,T,F)$, with
$F$ as in \eqref{eq:wtilde}.
Moreover, by the definition of $Z^n$ and $\Del^n$, we have
$0\le -\Del^n_t\le c_3\hat\Th^n$, whereas
by Lemma \ref{lem:reflectionprelim}, $\bar{R}^n\to0$ in probability.
By the continuity of the map $\Gam$ we therefore conclude
from \eqref{303} that, on the same subsequence,
\[
(\hat A^n,\hat S^n,T^n,\hat F^n,\hat  W^n,\hat L^n)
\To(A,S,T,F,W,L),
\]
where $(W,L)=\Gam(F)$.
This completes the proof.
\qed

\subsubsection{Proof of Lemma \ref{lem:reflectionprelim}}
\label{sec525}

We first explain how the proof changes between the cases. 
\begin{itemize}
	\item Under the \p{} policy, both servers can process the low priority jobs. This means that idling can only occur if the low priority class has few jobs and in that case the total number of jobs is also low by Proposition \ref{lem:ssc}.
	\item Under the \teetwo{} rule, no server can incur idleness when the high priority class is above the threshold. In addition, the high priority class takes a small time to reach above $\Th^n$ and does not empty below two jobs after that time with high probability  unless the total workload is close to $0$ (Lemma \ref{lem:nonidling2}). 
	This means that neither server  will  idle  except when there are almost no jobs in the system.
	\item Under the \p{}\p{} policy, Proposition \ref{lem:ssc} is still enough to prove Lemma \ref{lem:reflectionprelim} in the same way as in the single mode case \p{}.
	\item Under the \teetwo{}\teetwo{} policy, Lemma \ref{lem:nonidling1} and \ref{lem:nonidling2} hold. Because of that, Lemma \ref{lem:reflectionprelim} holds for the same reason as in the non switching case \teetwo{}.
	\item 
	When using a \p{} rule, the high priority class could become  zero with a lot of LPC jobs in the system, and switching to the \teetwo{} rule could lead to idleness even though there are a lot of low priority jobs (approximately $\alpha_q z^*$). Thus we introduce the \teeone{} rule in place of the \p{} rule in this case. This ensures that the number of high priority jobs does not decrease too much during the corresponding period.
\end{itemize}

In some cases, some idleness can occur when the high priority class starts with too few jobs but in those cases, it takes a small time to leave such states.

In cases 1(a) and 2(a), Lemma \ref{lem:reflectionprelim} is a direct consequence of the state space collapse.  Let us introduce two random times and a lemma that we will use in cases 1(b), 2(b)--(d). Fix $t_0>0$.
Let
\[\rho^n = \inf\left\lbrace t\geqslant 0:\, \hat{X}^n_p(t)\geqslant \hat{\Th}^n\right\rbrace\w t_0, \]
\[ \tau^n_r = \inf\left\lbrace t>\rho^n:\, \hat{X}^n_p(t)=2n^{-1/2},\,  \text{ and } \hat{X}^n_{q}(t)\geqslant  \hat\Th^n\right\rbrace.\]

\begin{lemma}\label{lem:nonidling1}
	If any of the following conditions hold, we have	
	\begin{equation}\label{eq:nonidlingpart1}
		\bar R^n_{\rho^n}\to0 \text{ in probability}.
	\end{equation}	
	\begin{itemize}
		\item Case 1(b), \eqref{90} and $i_2(\xi^A)=p$.
		\item  Case 2(b), \eqref{91} and $i_2(\xi^L)=i_2(\xi^H)=p$.
		\item Cases 2(c) and 2(d), \eqref{91} and $i_1(\xi^L)=i_2(\xi^H)$.
	\end{itemize}
	
\end{lemma}
\begin{lemma}
	Under the same assumptions as the previous lemma,\label{lem:nonidling2}
	\begin{equation}\label{eq:nonidlingpart2} \mathbb{P}\left(\tau^n_r\leqslant t_0\right)\to 0.\end{equation}
\end{lemma}
\begin{remark}
	We do not prove the previous lemma in all cases because we can directly prove Lemma \ref{lem:reflectionprelim} when using a \p{} or \p{}\p{} policy in cases 1(a) and 2(a).
\end{remark}
We now use a similar reasoning as in Proposition \ref{lem:ssc} to prove Lemma \ref{lem:nonidling1}.

\noi{\bf Proof of Lemma \ref{lem:nonidling1}.}
Fix $\delta$. Let us introduce  
\[\sigma_r^n= \sup\left\lbrace t\leqslant \tau_r^n : \, \hat{X}^n_p(t)\geqslant \hat\Th^n,\, \text{ or } \hat{X}^n_{q}(t)\leqslant \frac{\hat\Th^n}{2}\right\rbrace.  \]
In cases 2(b), (c) and (d), introduce  
\[\widetilde{\tau}^n= \inf\left\lbrace t\geqslant 0: \, \hat{X}^n_q(t)\geqslant \frac{z^*\al_q}{2}\right\rbrace, \]
\[\widetilde{\sigma}^n= \sup\left\lbrace t\leqslant \widetilde{\tau}^n: \, \hat{X}^n_q(t)\leqslant \frac{z^*\al_q}{4}\right\rbrace. \]
To simplify, we write throughout this proof
$$\rho^n=\rho \text{, } \tau_r^n=\tau\text{, }\widetilde{\tau}^n=\widetilde{\tau}\text{, }\widetilde{\sigma}^n=\widetilde{\sigma} \text{ and } \sigma_r^n=\sigma.$$ 
We briefly describe the interaction between the times we just introduced.
Under the state space collapse, $\tilde{\tau}$ has to occur before the first time the current mode switches. The times $\sigma$ and $\tilde{\sigma}$ allow us to have some knowledge about the state of queue lengths, uniformly over an interval. The first step of the proof of \eqref{eq:nonidlingpart1} is establishing that
\begin{equation}
	\mathbb{P}\left(\tau^n_c\geqslant t_0,\, \widetilde{\tau}\leqslant \rho\right)\to 0\label{eq:step1}
\end{equation}
holds in cases 2(b), (c) and (d).
In those cases, under the event $\left\lbrace \tau^n_c\geqslant t_0\right\rbrace\cap \left\lbrace \widetilde{\tau}\geqslant \rho\right\rbrace$, for any $t\leqslant\rho$, 
\begin{equation}\label{eq:samode}
	\mathrm{MODE}^n(t)=\mathrm{MODE}^n(0)=\xi^L.
\end{equation}
In addition, in case 1(b), for any $t\leqslant\rho$
\begin{equation}\label{eq:samodebis}
	\mathrm{MODE}^n(t)=\xi^A,
\end{equation}
which means \eqref{eq:step1} is not needed in this case.

We now proceed with the proof of \eqref{eq:step1}. On $\left\lbrace \tau^n_c\geqslant t_0,\, \widetilde{\tau}\leqslant \rho\right\rbrace$, the following hold:
\begin{enumerate}
	\item The number of HPC job is below $\hat{\Th}^n$ and there are LPC jobs in the system during $[\widetilde{\sigma},\widetilde{\tau}]$:
	\begin{equation*}
		\sup_{t\in [\widetilde{\sigma},\widetilde{\tau}]}\hat{X}^n_p(t)<\hat{\Th}^n \text{ and } \inf_{t\in [\widetilde{\sigma},\widetilde{\tau}]}\hat{X}^n_q(t)>0.
	\end{equation*}
	\item The current mode is the same as in the initial state: for any $t \in  [\widetilde{\sigma},\widetilde{\tau}]$,
	\begin{equation*}
		\mathrm{MODE}^n(t)=\mathrm{MODE}^n(0)=\xi^L.
	\end{equation*}
	\item In addition, the number of LPC jobs must have grown during that time:
	\begin{equation*}
		\hat{X}^n_q(\widetilde{\tau})-\hat{X}^n_q(\widetilde{\sigma})\geqslant \frac{z^*\al_q}{4}.
	\end{equation*}
\end{enumerate}
1. comes from the definition of $\rho$, $\tilde{\sigma}$ and $\tilde{\tau}$. 2. comes from the fact that  the number of HPC jobs is bounded by $2\Th^n$ under $\lbrace \tau^n_c\geqslant t_0\rbrace$ and the number of LPC jobs is bounded by $\frac{z^*\al_q}{2}$ before $\tilde{\tau}$ and the workload cannot cross above $z^*$. 3. comes from the definition of $\tilde{\tau}$ and $\tilde{\sigma}$. Because of 2., only the rule in the lower workload mode is in use. When using a \teeone{} rule, both servers take new jobs from the LPC because of 1. and $\lambda_q-\sum_k\mu_{qk}<0$. When using a \teetwo{} rule, again because of 1., the dual activity server takes new jobs from the LPC. Since $\xi^L$ is in canonical form $\frac{\lambda_1}{\al_1}>\beta_1$, $k_2(\xi^L)=2$ and $q=2$ because of the \teetwo{} rule. This means that $\frac{\lambda_2}{\alpha_2}<\beta_2$ and thus $\lambda_q-\mu_{qk_2(\xi^L)}<0$.  Hence, with \eqref{eq:balancexxx}, regardless of the case there exists $c>0$ such that
$$\hat{X}_q^n(\widetilde{\tau})\leqslant\hat{X}_q^n(\widetilde{\sigma})+\hat{f}^n_q[\widetilde{\sigma},\widetilde{\tau}]-c\sqrt{n}(\widetilde{\tau}-\widetilde{\sigma})+\sum_k\mu_{qk}e^n_{\max}.$$
From this, we obtain two bounds:
\begin{align*}
	\hat{X}_q^n(\widetilde{\tau})-\hat{X}_q^n(\widetilde{\sigma})&\leqslant w_{t_0}(\hat{f}^n_q, \widetilde{\tau}-\widetilde{\sigma})-c\sqrt{n}(\widetilde{\tau}-\widetilde{\sigma})+\sum_k\mu_{qk}e^n_{\max}\\
	\hat{X}_q^n(\widetilde{\tau})-\hat{X}_q^n(\widetilde{\sigma})&\leqslant 2\lVert \hat{f}^n_q\rVert_{t_0}-c\sqrt{n}(\widetilde{\tau}-\widetilde{\sigma})+\sum_k\mu_{qk}e^n_{\max}
\end{align*}

We now prove \eqref{eq:step1}: let $r_n=\frac{\log(n+1)}{\sqrt{n}}$. Because of 3.,
\begin{align*}
	&\mathbb{P}\left(\tau^n_c\geqslant t_0,\, \widetilde{\tau}\leqslant \rho\right)\leqslant \mathbb{P}\left( \widetilde{\tau}\leqslant \rho,\,\hat{f}^n_q[\widetilde{\sigma},\widetilde{\tau}]-c\sqrt{n}(\widetilde{\tau}-\widetilde{\sigma})+\sum_k\mu_{qk}e^n_{\max}\geqslant \frac{z^*\al_q}{4}\right)\\
	&\leqslant \mathbb{P}\left( w_{t_0}(\hat{f}^n_q, \widetilde{\tau}-\widetilde{\sigma})-c\sqrt{n}(\widetilde{\tau}-\widetilde{\sigma})+\sum_k\mu_{qk}e^n_{\max}\geqslant \frac{z^*\al_q}{4}, \widetilde{\tau}-\widetilde{\sigma}\leqslant r^n\right)\\
	&\quad +\mathbb{P}\left(2\lVert \hat{f}^n_q\rVert_{t_0}-c\sqrt{n}(\widetilde{\tau}-\widetilde{\sigma})+\sum_k\mu_{qk}e^n_{\max}\geqslant \frac{z^*\al_q}{4}, \widetilde{\tau}-\widetilde{\sigma}> r^n\right)\\
	&\leqslant \mathbb{P}\left(w_{t_0}(\hat{f}^n_q, r_n)+\sum_k\mu_{qk}e^n_{\max}\geqslant \frac{z^*\al_q}{4}\right)+\mathbb{P}\left(2\lVert \hat{f}^n_q\rVert_{t_0}+\sum_k\mu_{qk}e^n_{\max}\geqslant c\log(n+1)\right). 
\end{align*}
By Lemma \ref{lem:1service}, $e^n_{\max}$ is smaller than $n^{\bar a-1}=o(1)$. By Remark \ref{rem:stillok}, $\hat{f}^n_q$ are $C$-tight and $\lVert \hat{f}^n_q\rVert_{t_0}$ are tight RVs. Both terms must then converge to zero. This concludes the proof of \eqref{eq:step1} in all relevant cases for this lemma.

\textbf{\underline{Case 2(c), \teeone{} rule:}}\\
By splitting the integral that defines $\bar R^n$,
and using \eqref{eq:rbar} and \eqref{eq:samode},
\begin{align*}
	\bar R^n_{\rho}& \leqslant \int_{0}^{\rho}\one_{\hat{W}^n_t\geqslant c_3 \hat\Th^n ,\,\tau^n_c\geqslant t_0,\,\widetilde{\tau}\geqslant \rho}d\hat{L}^n_t+\int_{0}^{\rho}\one_{\tau^n_c< t_0}+\one_{\widetilde{\tau}\leqslant \rho,\,\tau^n_c\geqslant t_0}d\hat{L}^n_t\\
	&\leqslant \int_{0}^{\rho}\one_{\hat{X}^n_q(t)\geqslant \hat{\Th}^n, \mathrm{MODE}^n(t)=\xi^L}d\hat{L}^n_t+\int_{0}^{\rho}\one_{\tau^n_c< t_0}+\one_{\widetilde{\tau}\leqslant \rho,\,\tau^n_c\geqslant t_0}d\hat{L}^n_t
\end{align*}
The first term is zero by definition of the \teeone{} rule because $q$ is the dual activity class. The second term goes to zero by Proposition \ref{lem:ssc} and \eqref{eq:step1}. This concludes the proof of \eqref{eq:nonidlingpart1} in case 2(c).

\textbf{\underline{Case 1(b), 2(b) and 2(d), \teetwo{} rule:}}\\
We now deal with all cases that use a \teetwo{} rule in lower workload at the same time thanks to \eqref{eq:samode}/\eqref{eq:samodebis}.  Recall that $p=1=i_2(\xi^L)$ and $k_1(\xi^L)=1$.  By splitting the integral, we obtain
\begin{align*}
	\bar R^n_{\rho}&=\int_{0}^{\rho}\one_{\hat{W}^n_t\geqslant c_3 \hat\Th^n ,\,\tau^n_c\geqslant t_0}d\hat{L}^n_t+\int_{0}^{\rho}\one_{\hat{W}^n_t\geqslant c_3 \hat\Th^n ,\,\tau^n_c< t_0}d\hat{L}^n_t.
\end{align*}
By \eqref{004+}, $\hat{L}^n=  \sum_k\beta_k\hat{I}_k^n$. By definition of the \teetwo{} rule, as long as there are at least $2$ customers in the system, the dual activity server cannot idle. 
With $\xi^L$ in canonical form, server $2$ is the dual activity server, so
under the event $\left\lbrace \tau^n_c\geqslant t_0\right\rbrace\cap \left\lbrace \widetilde{\tau}\geqslant \rho\right\rbrace$,
\begin{equation*}\int_{0}^{\rho}\one_{\hat{W}^n_t\geqslant 2(\alpha_1 \wedge \alpha_2)^{-1}n^{-1/2},\,  \tau^n_c\geqslant t_0,\, \widetilde{\tau}\geqslant \rho}d\hat{I}^n_2(t)=0.\end{equation*}
Since $\int_{0}^{\rho}\one_{\hat{W}^n_t\geqslant c_3 \hat\Th^n ,\,\tau^n_c< t_0}d\hat{L}^n_t\Rightarrow 0$ by Proposition \ref{lem:ssc}, we are left to deal with 
$$\int_{0}^{\rho}\one_{\hat{W}^n_t\geqslant c_3 \hat\Th^n ,\,\tau^n_c\geqslant t_0,\,\widetilde{\tau}\geqslant \rho }d\hat{L}^n_t= \beta_1\int_{0}^{\rho}\one_{\hat{W}^n_t\geqslant c_3 \hat\Th^n ,\,\tau^n_c\geqslant t_0,\,\widetilde{\tau}\geqslant \rho}d\hat{I}_1^n(t).$$
We introduce the following times:
\[\bar{\tau}^n= \inf\left\lbrace t\geqslant 0: \int_{0}^{t}\one_{\hat{W}^n_t\geqslant c_3 \hat\Th^n}d\hat{L}^n_t\geqslant \frac{\delta}{2}\right\rbrace,\]
and
\[\bar{\sigma}^n= \sup\left\lbrace t\leqslant \bar{\tau}^n: \, \hat{X}^n_2(t)=0\right\rbrace.\]
We want to show that 
$$\mathbb{P}\left(\rho\geqslant \bar{\tau}\right)\to 0.$$
On the event $\left\lbrace \rho\geqslant \bar{\tau}\right\rbrace$, we have :
\begin{itemize}
	\item First, by definition of $\bar{\sigma}$,
	$$\inf_{t\in(\bar{\sigma},\bar{\tau})}\hat{X}^n_2(t)>0.$$
	\item 
	Second, by definition of $\rho$, 
	$$\sup_{t\in [0,\bar{\tau}]}\hat{X}^n_1(t)<\hat{\Th}^n.$$
	\item Finally, by definition of $\bar{\sigma}$, $$\hat{X}^n_2(\bar{\sigma}) \leqslant n^{-1/2}.$$
	In order for server 1 to be idle at time $\bar{\tau}$ (so that $d\hat{I}_1^n(\bar{\tau})>0$), it is necessary that  $\hat{X}^n_1(\bar{\tau})=0$. This means that in order to also have $\hat{W}^n_{\bar{\tau}}\geqslant c_3 \hat\Th^n$, we need to have $ \hat{X}^n_2(\bar{\tau}) \geqslant 2 \hat\Th^n$, or
	\begin{equation}
		\label{eq:need1}
		\hat{X}^n_2(\bar{\tau})-\hat{X}^n_2(\bar{\sigma})\geqslant 2 \hat\Th^n -n^{-1/2}  \geqslant \hat\Th^n .
	\end{equation}
	
\end{itemize}
We can now give the balance equation for $\hat{X}^n_2$  between $\bar{\sigma}$ and $\bar{\tau}$ under the \teetwo{} rule. 
Under the \teetwo{} rule, server 2 prioritizes class $2$ for the whole period because the number of HPC jobs remains below $\hat{\Th}^n$. Since $\lambda_1>\mu_{11}$ we obtain $\lambda_2<\mu_{22}$ from \eqref{eq:ht}. By the same reasoning as \eqref{300} there exists $c_2>0$ such that
\begin {equation}\label{eq:finaldist2}
\hat{X}_2^n[\bar{\sigma},\bar{\tau}] \leqslant \hat{f}^n_2[ \bar{\sigma},\bar{\tau}]-c_2\sqrt{n}(\bar{\tau}-\bar{\sigma}) +\sqrt{n}e^n_{\max}\sum_k\mu_{2k}. 
\end{equation}
Combining \eqref{eq:need1} and \eqref{eq:finaldist2} we obtain (on $\left\lbrace \rho\geqslant \bar{\tau}\right\rbrace$)
$$\hat{f}^n_2[ \bar{\sigma},\bar{\tau}]-c_2\sqrt{n}(\bar{\tau}-\bar{\sigma}) +\sqrt{n}e^n_{\max}\sum_k\mu_{2k} \geqslant \hat\Th^n.$$
As in the proof of Proposition \ref{lem:ssc}, we  distinguish between two cases: $\bar{\tau}-\bar{\sigma}$ smaller or larger than  $n^{-\nu_1}$. 
On $\bar{\tau}-\bar{\sigma} \leqslant n^{-\nu_1}$,
$\hat{f}^n_2[ \bar{\sigma},\bar{\tau}] \leqslant   w_{t_0}(\hat{f}^n_2, n^{-\nu_1}) $, so that
\begin{multline*}
\mathbb{P}\left(\hat{f}^n_2[ \bar{\sigma},\bar{\tau}]-c_2\sqrt{n}(\bar{\tau}-\bar{\sigma}) +\sqrt{n}e^n_{\max}\sum_k\mu_{2k} \geqslant \hat\Th^n, \bar{\tau}-\bar{\sigma} \leqslant n^{-\nu_1} \right)   \\
\leqslant \mathbb{P}\left( w_{t_0}(\hat{f}^n_2, n^{-\nu_1}) \geqslant \hat\Th^n/2 \right) + \mathbb{P}\left( \sqrt{n}e^n_{\max}\sum_k\mu_{2k} \geqslant \hat\Th^n/2 \right).
\end{multline*}
On $\bar{\tau}-\bar{\sigma} > n^{-\nu_1}$, $\hat{f}^n_2[ \bar{\sigma},\bar{\tau}] \leqslant 2\lVert \hat{f}^n_2 \rVert_{t_0}$ 
and $c_2\sqrt{n}(\bar{\tau}-\bar{\sigma}) \geqslant c_2 n^{\frac{1}{2}-\nu_1}$, so that
\begin{multline*}
\mathbb{P}\left(\hat{f}^n_2[ \bar{\sigma},\bar{\tau}]-c_2\sqrt{n}(\bar{\tau}-\bar{\sigma}) +\sqrt{n}e^n_{\max}\sum_k\mu_{2k} \geqslant \hat\Th^n, \bar{\tau}-\bar{\sigma} > n^{-\nu_1} \right) \leqslant  \\
 \mathbb{P}\left( 2\lVert \hat{f}^n_2 \rVert_{t_0} \geqslant c_2 n^{\frac{1}{2}-\nu_1}/2 \right) + \mathbb{P}\left( \sqrt{n}e^n_{\max}\sum_k\mu_{2k} \geqslant c_2 n^{\frac{1}{2}-\nu_1}/2 \right).
\end{multline*}
To summarize, 
\begin{align*}
\mathbb{P}\left(\rho\geqslant \bar{\tau}\right)&\leqslant \mathbb{P}\left(\tau^n_c\leqslant t_0\right)+ \mathbb{P}\left(\tau^n_c\geqslant t_0,\, \widetilde{\tau}<\rho\right)+\mathbb{P}\left(\tau^n_c\geqslant t_0,\, \widetilde{\tau}\geqslant \rho,\, \bar{\tau}\leqslant \rho,\, \bar{\tau}-\bar{\sigma}\leqslant n^{-\nu_1} \right)\\
&\quad+\mathbb{P}\left(\tau^n_c\geqslant t_0,\, \widetilde{\tau}\geqslant \rho,\, \bar{\tau}\leqslant \rho,\, \bar{\tau}-\bar{\sigma}> n^{-\nu_1}\right)\\
&\leqslant \mathbb{P}\left(\tau^n_c\leqslant t_0\right)+ \mathbb{P}\left(\tau^n_c\geqslant t_0,\, \widetilde{\tau}<\rho\right)+ \mathbb{P}\left( w_{t_0}(\hat{f}^n_2, n^{-\nu_1}) \geqslant \hat\Th^n/2 \right) + \mathbb{P}\left( \sqrt{n}e^n_{\max}\sum_k\mu_{2k} \geqslant \hat\Th^n/2 \right)\\
&\quad +\mathbb{P}\left( 2\lVert \hat{f}^n_2 \rVert_{t_0} \geqslant c_2 n^{\frac{1}{2}-\nu_1}/2 \right) + \mathbb{P}\left( \sqrt{n}e^n_{\max}\sum_k\mu_{2k} \geqslant c_2 n^{\frac{1}{2}-\nu_1}/2 \right).
\end{align*}
All terms converge to zero. The first two have already been treated in Proposition \ref{lem:ssc} and \eqref{eq:step1} respectively. 
The third term converges by Lemma \ref{lem10}, the fourth and sixth by Lemma \ref{lem:1service}, and the fifth by tightness of $\hat{f}^n_2$.

This completes the proof of \eqref{eq:nonidlingpart1} and Lemma \ref{lem:nonidling1} in all of the relevant cases. \qed\vspace{0.3cm}

\noi{\bf Proof of Lemma \ref{lem:nonidling2}.}

\textbf{\underline{Case 1(b), \teetwo{} policy:}}

We begin the proof by a full analysis in the case where a \teetwo{} policy is used, and then explain how to deal with the other cases. In this case, with the active mode in canonical form, $p=1$. Recall the definition of the random time $\tau$:
\[ \tau= \inf\left\lbrace t>\rho: \, \hat{X}^n_1(t)=2n^{-1/2},\,  \text{ and } \hat{X}^n_2(t)\geqslant  \hat\Th^n\right\rbrace,\]
and   
\[\sigma= \sup\left\lbrace t\leqslant \tau^n : \, \hat{X}^n_1(t)\geqslant \hat\Th^n,\, \text{ or } \hat{X}^n_{2}(t)\leqslant \frac{\hat\Th^n}{2}\right\rbrace.\]  
For any $t\in [\sigma, \tau]$, 
\begin{itemize}
\item $\hat{X}^n_1(t) < \hat\Th^n$,
\item  and $\hat{X}^n_2(t) > \frac{\hat\Th^n}{2}$.
\end{itemize}

In this case, only server 1 processes class 1 jobs and does so at most at full rate. Indeed, only server 1 gives priority to class 1 jobs and the other server is busy with class 2 jobs present in the system. Similarly, only server 2 processes class 2 jobs and does so at full rate since there are always class 2 jobs and the number of HPC jobs is below $\Th^n$ between $\sigma$ and $\tau$. Similarly to Lemma \ref{lem:servicemax}, $e_{\max}^n$ is such that for any $s, t\in [\sigma, \tau]$, $s\leqslant t$,  including a service that could start before $\sigma$ at the wrong activity (server 2 in this case),
\[\int_{s}^{t}\big(\lambda_1-\sum_k\mu_{1k}\X^n_{1k}(s)\big)ds\geqslant (\lambda_1-\mu_{11})(t-s)-\mu_{12}e_{\max}^n,\]
and excluding a service of a HPC job by server 2 starting before $\sigma$,
\[\int_{s}^{t}\big(\lambda_2-\sum_k\mu_{2k}\X^n_{2k}(s)\big)ds\leqslant (\lambda_2-\mu_{22})(t-s)+\mu_{22}e_{\max}^n.\]
With the active mode in canonical form, we have $\lambda_1>\mu_{11}$. By \eqref{eq:ht}, this means $\lambda_2<\mu_{22}$. Thus $\lambda_1-\mu_{11}>0$ and $\lambda_2-\mu_{22}<0$.

There are two possibilities for the state of the system at time $\sigma$: we have either 
\begin{itemize}
\item $\hat{X}^n_1(\sigma)  = \hat\Th^n-n^{-1/2}$,
\item  or $\hat{X}^n_2(\sigma)\in [\frac{ \hat\Th^n}{2}+\frac{1}{2}n^{-1/2}, \frac{ \hat\Th^n}{2}+n^{-1/2}]$.
\end{itemize}
We will decompose the event $\Omega^n\coloneqq \lbrace \tau\leqslant t_0\rbrace$ in two events:
\[ \Omega^n=\Omega^n_1\cup \Omega^n_2,\]
with 
$$\Omega^n_1\coloneqq \Omega^n \cap \left\lbrace \hat{X}^n_1(\sigma) = \hat\Th^n-n^{-1/2}\right \rbrace ,$$
and 
$$\Omega^n_2\coloneqq \Omega^n \cap\left \lbrace \hat{X}^n_2(\sigma)\in [\frac{ \hat\Th^n}{2}+\frac{1}{2}n^{-1/2}, \frac{ \hat\Th^n}{2}+n^{-1/2}]\right \rbrace.$$
To lighten notation, introduce also
$$\widetilde{\Omega}^n\coloneqq \left\lbrace \tau\leqslant t_0,\, \sup\limits_{ t\in [\sigma, \tau]}\hat{X}^n_1(t) < \hat\Th^n, \inf\limits_{t\in [\sigma, \tau]}\hat{X}^n_2(t)> \dfrac{\hat\Th^n}{2}\right\rbrace.$$

Similarly to Proposition \ref{lem:ssc}, let $\nu_2=1/2 -\bar a \leqslant 1/4$ and $\nu_1\in (\nu_2 , \,1/2 )$ so that $\hat\Th^n = n^{-1/2}\lceil n^{1/2-\nu_2} \rceil \geqslant n^{-\nu_2}$. Then
\begin{align*}
\mathbb{P}\left(\Omega^n_1\right)&= \mathbb{P}\left(\hat{X}^n_1(\sigma) = \hat\Th^n-n^{-1/2},\, \hat{X}^n_1(\tau)\leqslant 2n^{-1/2},\, \widetilde{\Omega}^n\right)\\
&\leqslant \mathbb{P}\left(  \hat{X}^n_1(\tau)-\hat{X}^n_1(\sigma)\leqslant -\dfrac{\hat\Th^n}{2},  \widetilde{\Omega}^n\right)\\
&= \mathbb{P}\left( \hat{X}^n_1(\tau)-\hat{X}^n_1(\sigma)\leqslant -\dfrac{\hat\Th^n}{2},  \widetilde{\Omega}^n, \tau-\sigma\leqslant n^{-\nu_1}\right)\\
&\quad+\mathbb{P}\left(  \hat{X}^n_1(\tau)-\hat{X}^n_1(\sigma)\leqslant -\dfrac{\hat\Th^n}{2}, \widetilde{\Omega}^n, \tau-\sigma> n^{-\nu_1}\right).
\end{align*}
On $\widetilde{\Omega}^n$, we also have 
\begin{multline*}
\hat{X}^n_1(\tau)-\hat{X}^n_1(\sigma)=\hat{f}_1^n(\tau)-\hat{f}_1^n(\sigma)+\sqrt{n}\int_{\sigma}^{\tau}\left(\lambda_1-\sum_k\mu_{1k}\X^n_{1k}(t)\right)dt\\
\geqslant  \hat{f}_1^n(\tau)-\hat{f}_1^n(\sigma)+(\lambda_1-\mu_{11})\sqrt{n}(\tau-\sigma)-\mu_{12}e_{\max}^n\sqrt{n}.
\end{multline*}
On the event $\left\lbrace \widetilde{\Omega}^n,\,\tau-\sigma> n^{-\nu_1}\right\rbrace$, the reasoning is very similar to the proof of Proposition \ref{lem:ssc}:
\begin{multline}
\mathbb{P}\left(\hat{X}^n_1(\tau)-\hat{X}^n_1(\sigma)\leqslant -\dfrac{\hat\Th^n}{2}, \widetilde{\Omega}^n, \tau-\sigma> n^{-\nu_1}\right)\\
\leqslant \mathbb{P}\left(2\sup_{t\leqslant t_0}\left\lvert\hat{f}_1^n(t)\right\rvert +\mu_{12}\sqrt{n}e_{\max}^n\geqslant (\lambda_1-\mu_{11})\sqrt{n}n^{-\nu_1}/2 \right).\label{eq:decreaselarge}
\end{multline}
The last probability goes to zero by tightness of $\hat{f}_1^n$, $\nu_1<1/2$ and Lemma \ref{lem:1service}. When $\tau-\sigma\leqslant n^{-\nu_1}$ the situation is again the same as in the proof of Proposition \ref{lem:ssc}: 
\begin{multline}\label{eq:decreasesmall}
\mathbb{P}\left( \hat{X}^n_1(\tau)-\hat{X}^n_1(\sigma)\leqslant -\dfrac{\hat\Th^n}{2}, \widetilde{\Omega}^n, \,\tau-\sigma\leqslant n^{-\nu_1}\right)\\
\leqslant \mathbb{P}\left(w_{t_0}(\hat{f}^n_1, n^{-\nu_1})+\mu_{12}e_{\max}^n\sqrt{n}\geqslant \dfrac{n^{-\nu_2}}{2} \right).
\end{multline}
The right hand side also converges to zero by Lemmas \ref{lem10}, Remark \ref{rem:stillok} and Lemma \ref{lem:1service}, similarly to the proof of Proposition \ref{lem:ssc}.

The second event $\Omega^n_2$ is treated similarly:

\begin{align*}
\mathbb{P}\left(\Omega^n_2\right)&\leqslant\mathbb{P}\left(\hat{X}^n_2(\sigma)\in [\frac{ \hat\Th^n}{2}+\frac{1}{2}n^{-1/2}, \frac{ \hat\Th^n}{2}+n^{-1/2}],\, \hat{X}^n_2(\tau)\geqslant \hat\Th^n,\,  \widetilde{\Omega}^n\right)\\
&\leqslant \mathbb{P}\left( \hat{X}^n_2(\tau)-\hat{X}^n_2(\sigma)\geqslant \frac{\hat\Th^n}{3},  \widetilde{\Omega}^n\right)\\
&= \mathbb{P}\left( \hat{X}^n_2(\tau)-\hat{X}^n_2(\sigma)\geqslant \frac{\hat\Th^n}{3}, \widetilde{\Omega}^n, \tau-\sigma\leqslant n^{-\nu_1}\right)\\
&\quad+\mathbb{P}\left( \hat{X}^n_2(\tau)-\hat{X}^n_2(\sigma)\geqslant \frac{\hat\Th^n}{3},  \widetilde{\Omega}^n, \tau-\sigma> n^{-\nu_1}\right).
\end{align*}
On this event, we also have 
\begin{multline*}
\hat{X}^n_2(\tau)-\hat{X}^n_2(\sigma)=\hat{f}_2^n(\tau)-\hat{f}_2^n(\sigma)+\sqrt{n}\int_{\sigma}^{\tau}\left(\lambda_2-\sum_k\mu_{2k}\X^n_{2k}(t)\right)dt\\
\leqslant \hat{f}_2^n(\tau)-\hat{f}_2^n(\sigma)+(\lambda_2-\mu_{22})\sqrt{n}(\tau-\sigma)+\mu_{22}e_{\max}^n\sqrt{n}.
\end{multline*}
On the event $\left\lbrace \tau-\sigma> n^{-\nu_1}\right\rbrace$,  similar to the treatment of $\Omega^n_1$:
\begin{multline}
\mathbb{P}\left(\hat{X}^n_2(\tau)-\hat{X}^n_2(\sigma)\geqslant \frac{\hat\Th^n}{3},  \widetilde{\Omega}^n, \tau-\sigma> n^{-\nu_1}\right)\\
\leqslant \mathbb{P}\left(2\sup_{t\leqslant t_0}\lvert\hat{f}_2^n(t)\rvert +\mu_{22}\sqrt{n}e_{\max}^n\geqslant -(\lambda_2-\mu_{22})\sqrt{n}n^{-\nu_1}/3 \right).\label{eq:increaselarge}
\end{multline}
The last probability goes to zero by tightness of $\hat{f}_2^n$, $\nu_1<1/2$ and Lemma \ref{lem:1service}. When $\tau-\sigma\leqslant n^{-\nu_1}$ the situation is also the same as for $\Omega^n_1$: 
\begin{multline}
\mathbb{P}\left( \hat{X}^n_2(\tau)-\hat{X}^n_2(\sigma)\geqslant  \frac{\hat\Th^n}{3}, \widetilde{\Omega}^n, \tau-\sigma\leqslant n^{-\nu_1}\right)\\
\leqslant \mathbb{P}\left(w_{t_0}(\hat{f}_2^n, n^{-\nu_1})+\mu_{22}e_{\max}^n\sqrt{n}\geqslant \dfrac{n^{-\nu_2}}{3} \right).\label{eq:increasesmall}
\end{multline}
The right hand side also converges to zero by Lemma \ref{lem10} similarly to the $\Omega^n_1$ case. This concludes the proof in case 1(b).\vspace{0.5cm}

We now present the changes required to adapt the proof to the other cases described in the lemma:

\textbf{\underline{Case 2(b), \teetwo{}\teetwo{} policy:}}

This case involves switching between two \teetwo{} rules: $i_2(\xi^L)=i_2(\xi^H)=p$. Because  $\xi^L$ is in canonical form, $p=1$.  In this case, mode switching only occurs at a service completion at the single activity server that is dedicated to HPC jobs. 
There could be either type of job being served at either server immediately before time $\sigma$. We will see that after those jobs exit the system, there is at least one server processing the LPC and at most one processing the HPC. 
If the first job to finish is from the dual activity server, the service of a low priority job starts because there are fewer than $\Theta^n$ HPC jobs. This server will continue to take LPC jobs until the minimum between the next mode switching time and $\tau$. 
If the first job to finish at or after $\sigma$ is at the single activity server, either the current mode switches or the service of an HPC job starts. If there is no mode switch then this server will continue to take HPC jobs until the minimum between the next mode switching time and $\tau$. 
When the dual activity server completes its job it will, as noted before, serve LPC jobs.
If there is a mode switch then the formerly single activity server becomes dual activity, and (because there are fewer than $\Theta^n$ HPC jobs) begins service on an LPC job. When the formerly dual activity server completes its job, it becomes single activity and serves HPC. This continues until the minimum between the next mode switching time and $\tau$.

After the two jobs present at $\sigma$, there cannot be a time where both servers are occupied with HPC jobs. This is because the HPC is only processed by the server dedicated to it and there can be no residual service of an HPC job at the single activity server whenever the current mode switches. Similarly, after both initial jobs have been processed, there is always at least one server occupied with LPC jobs. This is because one server gives priority to the LPC and the service of a LPC job starts whenever the current mode switches.

Keeping the definition of $\tau$, $\sigma$, just as in the single mode case, for any $t,s\in [\sigma,\tau]$,  including/excluding a service that could start before $\sigma$ at the wrong activity, there is always at most one server processing HPC jobs and at least one server processing LPC jobs between $\sigma$ and $\tau$. Thus,
\[\int_{s}^{t}\big(\lambda_1-\sum_k\mu_{1k}\X^n_{1k}(s)\big)ds\geqslant (\lambda_1-\max_k\mu_{1k})(t-s)-\sum_k\mu_{1k}e_{\max}^n,\]
and 
\[\int_{s}^{t}\big(\lambda_2-\sum_k\mu_{2k}\X^n_{2k}(s)\big)ds\leqslant (\lambda_2-\min_k\mu_{2k})(t-s)+\sum_k\mu_{2k}e_{\max}^n.\]
With one mode in canonical form we have $\frac{\lambda_1}{\alpha_1}>\beta_1$. In addition, by Lemma \ref{rem:alg}, in case 2(b), \eqref{92} holds and thus $\frac{\lambda_1}{\alpha_1}>\max_k \beta_k$, which also means $\frac{\lambda_2}{\alpha_2}<\min_k \beta_k$ by \eqref{eq:ht}.

The rest of the proof is the same as in the single mode case, obtaining \eqref{eq:decreaselarge}, \eqref{eq:decreasesmall}, \eqref{eq:increaselarge} and \eqref{eq:increasesmall}. This concludes the proof in case 2(b).\vspace{0.5cm}

\textbf{\underline{Case 2(c), \teeone{}\teetwo{} policy:}}

With $\xi^L$ in canonical form, we have $p=2$, so   $\tau$ and $\sigma$ are defined as
\[ \tau= \inf\left\lbrace t>\rho:\, \hat{X}^n_2(t)=2n^{-1/2},\,  \text{ and } \hat{X}^n_1(t)\geqslant  \hat\Th^n\right\rbrace,\]
and   
\[\sigma= \sup\left\lbrace t\leqslant \tau^n : \, \hat{X}^n_2(t)\geqslant \hat\Th^n,\, \text{ or } \hat{X}^n_{1}(t)\leqslant \frac{\hat\Th^n}{2}\right\rbrace.\]  

In this case,  in the upper workload mode only the single activity server is allowed to begin service of the HPC between $\sigma$ and $\tau$,
while in the lower workload mode neither server can begin service of HPC between $\sigma$ and $\tau$. 
Between those times, there are always low priority class jobs to serve, and the number of HPC jobs stays below $\Th^n$. The most service the HPC can get between $\sigma$ and $\tau$ occurs if there are no switches and the current mode is always upper workload. In this mode, the single activity server (server 1) is dedicated to service of the HPC. Hence, including a service that could start before $\sigma$ at the wrong activity,
\[\int_{s}^{t}\big(\lambda_2-\sum_k\mu_{2k}\X^n_{2k}(s)\big)ds\geqslant (\lambda_2-\mu_{21})(t-s)-\mu_{21}e_{\max}^n.\]
Between $\sigma$ and $\tau$, server 2 gives priority to class 1 regardless of switches. Excluding a service of a HPC job that could have started before $\sigma$,
\[\int_{s}^{t}\big(\lambda_1-\sum_k\mu_{1k}\X^n_{1k}(s)\big)ds\leqslant (\lambda_1-\mu_{12})(t-s)-\mu_{12}e_{\max}^n.\]
Because we have one mode in canonical form, $\frac{\lambda_1}{\alpha_1}>\beta_1$. In addition, by Lemma \ref{rem:alg} in case 2(c), \eqref{93} holds, which means that $\frac{\lambda_1}{\alpha_1}<\beta_2$. Finally $\frac{\lambda_2}{\alpha_2}>\beta_1$ by the previous observation and \eqref{eq:ht}.

The rest of the proof is the same as in the single mode case, obtaining \eqref{eq:decreaselarge}, \eqref{eq:decreasesmall}, \eqref{eq:increaselarge} and \eqref{eq:increasesmall}.\vspace{0.5cm}

\textbf{\underline{Case 2(d), \teetwo\teeone{}policy:}}

This case is handled the same way as 2(c), by interchanging the roles of upper workload and lower workload mode.
\qed \vspace{0.3cm}

We now turn to the proof of the main result of this section.\vspace{0.3cm}

\noi{\bf Proof of Lemma \ref{lem:reflectionprelim}.}
We  distinguish along the same cases as in Theorem \ref{th-ao-s}.

\begin{itemize}
\item \textbf{\underline{Case 1(a), \p{} policy:}}
\end{itemize}
In this case $p=2$, and no server can be idle when $X_1 \geqslant 2$.	Thus
\begin{equation}
\int_0^{t_0}\mathds{1}_{\hat X^n_1(t)\geqslant 2n^{-1/2}}d\hat L^n_t=0.\label{q1}
\end{equation}
For any $\delta>0$, 
\begin{align*}
\mathbb{P}\left(\bar R^n_{t_0}\geqslant \delta\right)&\leqslant 	\mathbb{P}\left(\int_{0}^{t_0}(\mathds{1}_{\hat{X}^n_2(t)\geqslant 2\hat{\Th}^n }+\mathds{1}_{\hat{X}^n_1(t)\geqslant 2n^{-1/2}})d\hat L^n_t\geqslant \delta\right)\\
&=\mathbb{P}\left(\int_{0}^{t_0}\mathds{1}_{\hat{X}^n_2(t)\geqslant 2\hat\Th^n }d\hat L^n_t\geqslant \delta\right)\\
&\leqslant \mathbb{P}\left(\mathds{1}_{\sup_{t\leqslant t_0}\hat{X}^n_2(t)\geqslant2\hat\Th^n }\hat L^n(t_0)\geqslant \delta\right)\\
&\leqslant \mathbb{P}\left(\sup_{t\leqslant t_0}\hat{X}^n_2(t)\geqslant 2\hat\Th^n\right)\\
&=\mathbb{P}\left(\tau^n_c \leqslant t_0\right).
\end{align*}
By Proposition \ref{lem:ssc}, 
\begin{equation}
\label{q2}\mathbb{P}\left(\tau^n_c \leqslant t_0\right)\to 0 .
\end{equation}

\begin{itemize}
\item \textbf{\underline{Case 1(b), \teetwo{} policy:}}
\end{itemize}

In this case, we already know by Lemma \ref{lem:nonidling1} that 
$\mathbb{P}\left( \bar R^n_{\rho}\geqslant \frac{\delta}{2}\right)\to 0.$
We need to show the same thing for the time after $\rho$:
$$\mathbb{P}\left( \int_{\rho}^{t_0}\mathds{1}_{\hat{W}^n_t\geqslant c_3\hat{\Th}^n}d\hat{L}^n_t\geqslant \frac{\delta}{2}\right)\to 0.$$
First, by putting $\xi^A$ in canonical form, $p=1$. When $X^n_1(t)$ is above the threshold, both servers can serve class 1 jobs (high priority class) so almost surely,
\begin{equation}
\int_\rho^{t_0}\mathds{1}_{\hat{X}^n_1(t)\geqslant \hat\Th^n}d\hat L^n_t=0 \label{eq:idle1}.
\end{equation}
Similarly, 
\begin{equation}
\int_\rho^{t_0}\mathds{1}_{\hat{X}^n_2(t)\geqslant\label{eq:idle2} 2n^{-1/2}}d\hat{I}_2^n(t)=0,\end{equation}
and
\begin{equation}
\int_\rho^{t_0}\mathds{1}_{ 2n^{-1/2}\leqslant \hat{X}^n_1(t)\leqslant \hat\Th^n}d\hat{I}^n_1(t)=0 \label{eq:idle3}.
\end{equation}
Notice now that because of the three identities,
\begin{align*}
\mathbb{P}\left(\int_{\rho}^{t_0}\mathds{1}_{W^n_t\geqslant c_3\Th^n}d\hat L^n_t\geqslant \frac{\delta}{2}\right)&\leqslant \mathbb{P}\left(\int_{\rho}^{t_0}\left(\mathds{1}_{\hat{X}^n_1(t)\geqslant \hat\Th^n}+\mathds{1}_{\hat{X}^n_1(t)\leqslant \hat\Th^n, \hat{X}^n_2\geqslant \hat\Th^n}\right)d\hat L^n_t\geqslant \frac{\delta}{2}\right)\\
&=\mathbb{P}\left(\beta_1\int_{\rho}^{t_0}\mathds{1}_{\hat{X}^n_1(t)\leqslant \hat\Th^n, \hat{X}^n_2\geqslant \hat\Th^n}d\hat{I}^n_1(t)\geqslant \frac{\delta}{2}\right)\\
&=\mathbb{P}\left(\beta_1\int_{\rho}^{t_0}\mathds{1}_{\hat{X}^n_1(t)\leqslant 2n^{-1/2}, \hat{X}^n_2\geqslant \hat\Th^n}d\hat{I}^n_1(t)\geqslant \frac{\delta}{2}\right)
\end{align*}
By Lemma \ref{lem:nonidling2}, 
\begin{equation}\label{eq:idle4}\mathbb{P}\left(\int_{\rho}^{t_0}\mathds{1}_{\hat{X}^n_1(t)\leqslant 2n^{-1/2}, \hat{X}^n_2\geqslant \hat\Th^n}d\hat{I}^n_1(t)\geqslant \frac{\delta}{2}\right)\leqslant \mathbb{P}\left(\tau_r^n\leqslant t_0\right)\to 0  .\end{equation}

\begin{itemize}
\item \textbf{\underline{Case 2(a), \p{}\p{} policy:}}
\end{itemize}
In this case, $p=2$. The  reasoning in case 1(a) is still valid when switching between two P rules because \eqref{q1} and \eqref{q2} still hold.

\begin{itemize}
\item \textbf{\underline{Case 2(b), \teetwo{}\teetwo{} policy:}}
\end{itemize}
The result in this case has a proof very similar to case 1(b) because Lemmas \ref{lem:nonidling1} and \ref{lem:nonidling2} are still valid in this case.
We already know by Lemma \ref{lem:nonidling1} that
$$\mathbb{P}\left( \int_{0}^{\rho}\mathds{1}_{\hat{W}^n_t\geqslant c_3\hat{\Th}^n}d\hat{L}^n_t\geqslant \frac{\delta}{2}\right)\to 0.$$
We need to show the same thing for the time after $\rho$:
$$\mathbb{P}\left( \int_{\rho}^{t_0}\mathds{1}_{\hat{W}^n_t\geqslant c_3\hat{\Th}^n}d\hat{L}^n_t\geqslant \frac{\delta}{2}\right)\to 0.$$
First, by putting $\xi^L$ in canonical form $p=1$.

The idea is to split the integrals between the times $\mathrm{MODE}(t)=\xi^L$ and the times $\mathrm{MODE}(t)=\xi^H$. In this case, we still have \eqref{eq:idle1} but this time \eqref{eq:idle2} and \eqref{eq:idle3} only hold when $\mathrm{MODE}(t)=\xi^L$. In the other case, we have

\begin{equation*}
\int_\rho^{t_0}\mathds{1}_{\hat{X}^n_2(t)\geqslant 2n^{-1/2}, \mathrm{MODE}(t)=\xi^H}d\hat{I}_1^n(t)=0,\end{equation*}
\begin{equation*}
\int_\rho^{t_0}\mathds{1}_{ 2n^{-1/2}\leqslant \hat{X}^n_1(t)\leqslant \hat\Th^n,\mathrm{MODE}(t)=\xi^H}d\hat{I}^n_2(t)=0 .
\end{equation*}
We obtain
$$\int_{\rho}^{t_0}\mathds{1}_{\mathrm{MODE}(t)=\xi^L, W^n_t\geqslant  c_3\Th^n}d\hat L^n_t \leqslant \beta_1\int_{\rho}^{t_0}\mathds{1}_{\mathrm{MODE}(t)=\xi^L, \hat{X}^n_1(t)\leqslant 2n^{-1/2}, \hat{X}^n_2(t)\geqslant \hat\Th^n}d\hat{I}_1^n(t),$$
and
$$\int_{\rho}^{t_0}\mathds{1}_{\mathrm{MODE}(t)=\xi^H, W^n_t\geqslant  c_3\Th^n}d\hat L^n_t\leqslant \beta_2\int_{\rho}^{t_0}\mathds{1}_{\mathrm{MODE}(t)=\xi^H, \hat{X}^n_1(t)\leqslant 2n^{-1/2}, \hat{X}^n_2(t)\geqslant \hat\Th^n}d\hat{I}_2^n(t).$$
Finally, we obtain the result, since 
$$\mathbb{P}\left(\int_{\rho}^{t_0}\mathds{1}_{\hat{X}^n_1(t)\leqslant 2n^{-1/2}, \hat{X}^n_2(t)\geqslant \hat\Th^n}d\hat L^n_t\geqslant \frac{\delta}{2}\right)\leqslant \mathbb{P}\left(\tau_r^n \leqslant t_0\right)\to 0.$$

\begin{itemize}
\item \textbf{\underline{Case 2(c)(d), \teeone{}\teetwo{}/\teetwo{}\teeone{} policy:}}
\end{itemize}
We give the proof in case 2(d) but the reasoning is similar in case 2(c) by interchanging the role of $\xi^L$ and $\xi^H$. In this case $p=1$ with $\xi^L$ in canonical form. The idea is similar to the 2(b) case. We already know by Lemma \ref{lem:nonidling1} that 
$$\mathbb{P}\left( \int_{0}^{\rho}\mathds{1}_{\hat{W}^n_t\geqslant c_3\hat{\Th}^n}d\hat{L}^n_t\geqslant \frac{\delta}{2}\right)\to 0.$$
We need to show the same thing for the time after $\rho$:
$$\mathbb{P}\left( \int_{\rho}^{t_0}\mathds{1}_{\hat{W}^n_t\geqslant c_3\hat{\Th}^n}d\hat{L}^n_t\geqslant \frac{\delta}{2}\right)\to 0.$$

We will split the integral using $\mathds{1}_{\mathrm{MODE}(t)=\xi^L}+\mathds{1}_{\mathrm{MODE}(t)=\xi^H}$. 
We use a \teetwo{} rule when $W^n$ is small and a \teeone{} rule when $W^n$ is large but that distinction is not important here. In terms of almost sure non idling properties, we have
\begin{align*}
&\int_\rho^{t_0}\mathds{1}_{\mathrm{MODE}(t)=\xi^L,\hat{X}^n_1(t)\geqslant \hat\Th^n}d\hat L^n_t=0,
&\int_\rho^{t_0}\mathds{1}_{\mathrm{MODE}(t)=\xi^L,\hat{X}^n_2(t)\geqslant 2n^{-1/2}}d\hat{I}_2^n(t)(t)=0,\\
&\int_\rho^{t_0}\mathds{1}_{\mathrm{MODE}(t)=\xi^L, 2n^{-1/2}\leqslant \hat{X}^n_1(t)\leqslant \hat\Th^n}d\hat{I}^n_1(t)=0,
&\int_\rho^{t_0}\mathds{1}_{\mathrm{MODE}(t)=\xi^H, \hat{X}^n_2\geqslant 2n^{-1/2}}d\hat L^n_t=0.
\end{align*}
From these almost sure identities, we obtain
\begin{align*}
\mathbb{P}\left(\int_{\rho}^{t_0}\mathds{1}_{W^n_t\geqslant  c_3\Th^n}d\hat L^n_t\geqslant \frac{\delta}{2} \right)&\leqslant \mathbb{P}\left(\int_{\rho}^{t_0}\mathds{1}_{\mathrm{MODE}(t)=\xi^L, W^n_t\geqslant  c_3\Th^n}d\hat L^n_t\geqslant \frac{\delta}{4} \right)\\
&\quad+\mathbb{P}\left(\int_{\rho}^{t_0}\mathds{1}_{\mathrm{MODE}(t)=\xi^H, W^n_t\geqslant  c_3\Th^n}d\hat L^n_t\geqslant \frac{\delta}{4} \right)\\
&\leqslant  \mathbb{P}\left(\beta_1\int_{\rho}^{t_0}\mathds{1}_{\mathrm{MODE}(t)=\xi^L, \hat{X}^n_1(t)\leqslant 2n^{-1/2}, \hat{X}^n_2(t)\geqslant \hat\Th^n}d\hat{I}_1^n(t)\geqslant \frac{\delta}{4} \right)\\
&\quad +  \mathbb{P}\left(\int_{\rho}^{t_0}\mathds{1}_{\mathrm{MODE}(t)=\xi^H, \hat{X}^n_2(t)\leqslant 2n^{-1/2},\hat{X}^n_1(t)\geqslant 2\hat\Th^n}d\hat L^n_t\geqslant \frac{\delta}{4} \right)\\
&\leqslant \mathbb{P}\left(\tau_r^n\leqslant t_0\right)+\mathbb{P}\left(\tau_c^n\leqslant t_0\right)
\end{align*}
Both probabilities go to zero (by Lemma \ref{lem:nonidling2} for the first and Proposition \ref{lem:ssc} for the second) so the result is proved in this case as well.

As mentionened, the proof is the same in case 2(c): this time, the policy uses a \teeone{} rule when $W^n$ is small and \teetwo{} rule when $W^n$ is large. Keeping $\xi^L$ in canonical form, $p=2$. In addition, in terms of almost sure non idling properties, we have
\begin{align*}
&\int_\rho^{t_0}\mathds{1}_{\mathrm{MODE}(t)=\xi^H,\hat{X}^n_2(t)\geqslant \hat\Th^n}d\hat L^n_t=0,
&\int_\rho^{t_0}\mathds{1}_{\mathrm{MODE}(t)=\xi^H,\hat{X}^n_1(t)\geqslant 2n^{-1/2}}d\hat{I}_2^n(t)=0, \\
&\int_\rho^{t_0}\mathds{1}_{\mathrm{MODE}(t)=\xi^H, 2n^{-1/2}\leqslant \hat{X}^n_2(t)\leqslant \hat\Th^n}d\hat{I}^n_1(t)=0,
&\int_\rho^{t_0}\mathds{1}_{\mathrm{MODE}(t)=\xi^L, \hat{X}^n_1\geqslant 2n^{-1/2}}d\hat L^n_t=0.
\end{align*}
We can obtain the result using the same decomposition.

\subsubsection{Fast switching}\label{sec:fs}

\noi{\bf Proof of Proposition \ref{lem:correctmode}.}
Fix $t_0$ and $\eps>0$. Assume that $\xi^L$ is in canonical form.
Then $(i^l,k^l)=(2,1)$. Let
\begin{align*}
\tau_f&= \inf\left\lbrace t\geqslant 0: \, \int_{0}^{t}\mathds{1}_{ \hat{W}^n_t \leqslant  z^*-\eps}dT^n_{21}(t)>0\right\rbrace,\\
t_{\min}&= \sup\left\lbrace t\leqslant \tau_f: \, \hat{W}^n_t\geqslant z^*\right\rbrace,\\
t_{\max}&= \inf\left\lbrace t\geqslant t_{\min}: \, \hat{W}^n_t\leqslant z^*-\eps\right\rbrace,\\
\tau_1&= \inf\left\lbrace t\geqslant t_{\min}:  \mathrm{MODE}^n(t)=\xi^L\right\rbrace,
\end{align*}
where the dependence on $n$ is suppressed.
The first statement of the lemma will be proved once we show
$\mathbb{P}\left(\tau_f\leqslant t_0\right)\to 0.$
We omit the proof of the second statement, which is similar.
To this end, let
$$\kappa^n= \inf\lbrace s\geqslant 0\,\text{ s.t. } \exists t\leqslant t_0,\,  \hat{W}^n_{t-s}\geqslant z^*,\, \hat{W}^n_t\leqslant z^*-\eps\rbrace.$$
If $\sup_{t\leqslant t_0}\hat{W}^n_t\leqslant z^*$, the current mode never changes so the single activity server is dedicated to only one class for the whole period and the non-basic activity is never used. Thus, 
$$\left\lbrace \tau_f\leqslant t_0\right\rbrace \subset \left\lbrace t_{\max}\leqslant t_0\right\rbrace.$$
Note that we have used the fact that, for all of our policies,
jobs are never routed to a non basic activity of the current mode.

The remainder of the argument is based on the following fact,
which must be argued separately for each case.
This is concerned with
the difference  $\tau_1-t_{\min}$, which while
case dependent, can in all cases  be shown to satisfy
\begin{equation}
\lim\limits_{n\to +\infty}	\mathbb{P}\left(
\tau_1-t_{\min} > e_{\max}^n  \right)= 0.\label{eq:unifbound}
\end{equation}

By Proposition \ref{lem:reflection}, $\hat{W}^n$ is $C$-tight. Hence
for any constant $c>0$,
$\mathbb{P}\left( c\kappa^n\leqslant n^{\bar a-1}\right)\to 0$.
On the other hand, by Lemma \ref{lem:1service},
$\mathbb{P}\left( e_{\max}^n\geqslant n^{\bar a-1}\right)\to 0$.
Thus
\begin{equation}\label{eq:crossec}
\mathbb{P}\left(e_{\max}^n\geqslant c\kappa^n\right)\to 0.
\end{equation}
In order for $\int_{0}^{t_0}\mathds{1}_{ W^n_t \leqslant  z^*-\eps}dT^n_{21}(t)$ to become positive, $t_{\max}$ must be smaller than $t_0$ and 
$$t_{\max}-\tau_1 < e_{\max}^n.$$
It is not possible for $t_{\max}-\tau_1\geqslant e_{\max}^n$ to occur on $\left\lbrace \tau_f\leqslant t_0\right\rbrace$. If that were the case, server 1 would necessarily finish service of the job it was in the process of serving at time $\tau_1$ before the workload has time to reach $z^*-\eps$. 
When $\mathrm{MODE}^n(t)=\xi^L$ in canonical form, server 1 can only take new class 1 jobs regardless of the rule. Even if class 1 has no job in the queue, the non basic activity is not used after the possible residual job that was in service at time $\tau_1$. In addition, by definition of $t_{\min}$, this is the last time the mode switches from upper to lower workload before $\tau_f$. This would prevent $\int_{0}^{\tau_f}\one_{\hat{W}^n_t\leqslant z^*-\eps}dT^n_{21}(t)$ from becoming positive and is also the reason why $t_{\max}$ needs to be smaller than $t_0$ for $\tau_f$ to be smaller than $t_0$. 

By definition of $t_{\max}$ and $\kappa^n$, 
\begin{equation*}
\mathbb{P}\left( t_{\max}\leqslant t_0,\, t_{\max}-t_{\min}< \kappa^n\right)=0.
\end{equation*}
We next show that
\begin{equation}
\lim\limits_{n\to +\infty}	\mathbb{P}\left(t_{\max}\leqslant t_0,\, \tau_1\geqslant t_{\max}\right)= 0.\label{eq:fastswitch}
\end{equation}
We have
\begin{align*}
\mathbb{P}\left(t_{\max}\leqslant t_0,\, \tau_1\geqslant t_{\max}\right) &= \mathbb{P}\left(t_{\max}\leqslant t_0,\, \tau_1\geqslant t_{\max}, t_{\max}-t_{\min} \geqslant \kappa^n\right) \\
&\leqslant  \mathbb{P}\left( \tau_1\geqslant t_{\min} + \kappa^n\right)\\
&\leqslant  \mathbb{P}\left( \tau_1\geqslant t_{\min} + \kappa^n, e_{\max}^n \leqslant  \kappa^n\right) + \mathbb{P}\left( \tau_1\geqslant t_{\min} + \kappa^n , e_{\max}^n >  \kappa^n\right)\\
&\leqslant  \mathbb{P}\left( \tau_1\geqslant t_{\min} + e_{\max}^n \right) + \mathbb{P}\left( e_{\max}^n >  \kappa^n\right).
\end{align*}
Both terms converge to zero, the first by \eqref{eq:unifbound}, and the second by \eqref{eq:crossec}.

We can now prove the lemma based on
\eqref{eq:fastswitch}, \eqref{eq:unifbound} and \eqref{eq:crossec}. We have
\begin{align*}
\mathbb{P}\left(\tau_f\leqslant t_0\right)&=	\mathbb{P}\left(\tau_f\leqslant t_0,\, t_{\max}\leqslant t_0\right)\\
&=\mathbb{P}\left(\tau_f\leqslant t_0,\, t_{\max}\leqslant t_0,\,\tau_1\geqslant t_{\max}\right)  +  \mathbb{P}\left(\tau_f\leqslant t_0,\, t_{\max}\leqslant t_0,\, \tau_1\leqslant t_{\max},\,  t_{\max}-\tau_1 < e_{\max}^n\right)\\
&\leqslant \mathbb{P}\left(t_{\max}\leqslant t_0, \tau_1\geqslant t_{\max}\right)+ \mathbb{P}\left(t_{\max}\leqslant t_0,\,  t_{\max}-t_{\min}\leqslant 2e_{\max}^n\right) + \mathbb{P}\left(\tau_1-t_{\min} > e_{\max}^n  \right)\\
&\leqslant \mathbb{P}\left(t_{\max}\leqslant t_0, \tau_1\geqslant t_{\max}\right)+\mathbb{P}\left(\kappa^n \leqslant 2e_{\max}^n\right)+\mathbb{P}\left( t_{\max}\leqslant t_0,\, t_{\max}-t_{\min}< \kappa^n\right)\\
&+ \mathbb{P}\left(\tau_1-t_{\min} > e_{\max}^n  \right).
\end{align*}
The first term goes to zero by \eqref{eq:fastswitch}, the second by \eqref{eq:crossec}, the third is zero, and the fourth goes to zero by \eqref{eq:unifbound}. This completes the proof.

It remains to prove \eqref{eq:unifbound}. As noted above, the proof differs by case.

\textbf{Case 2(a):} 
The current mode changes to $\xi^L$ after the first service completion of a job at the single activity server  for $\xi^H$ (which is server 2) at or after $t_{\min}$.
Note that $t_{\min}$ must correspond to a service completion. If $t_{\min}$ corresponds to a service completion at server 2, the $\tau_1 = t_{\min}$. If $t_{\min}$ corresponds to a service completion at server 1, then the mode will not cange, and $\tau_1$ will correspond to the next service completion at server 2.
This will occur before $t_{\min}+ e_{\max}^n$ if server 2 is busy (serving class 1) at $t_{\min}$.
Note that, on $\left\lbrace \tau^n_c\geqslant t_0\right\rbrace$, $\hat{X}_2^n(t_{\min}) < 2 \hat\Th^n$, so that 
$$X_1^n(t_{\min}) \geqslant \alpha_1 (W^n_{t_{\min}}-2 \alpha_2^{-1} \Th^n)
\geqslant  \alpha_1 ( n^{1/2} z^* -1 -2 \alpha_2^{-1} \Th^n) \geqslant 2,$$
so server 2 is busy at  $t_{\min}$. Thus
\begin{equation*}
\mathbb{P}\left(\tau_1-t_{\min}> e_{\max}^n\right)\leqslant \mathbb{P}\left(\tau^n_c\leqslant t_0\right),
\end{equation*}
which converges to zero by Proposition \ref{lem:ssc}. This proves \eqref{eq:unifbound}.

\textbf{	Case 2(b):} The current mode changes after the first service completion of a job at the single activity server (which is server 2) at or after $t_{\min}$. Server 2 is dedicated to HPC jobs.  
Under $\left\lbrace \tau^n_r\geqslant t_0, \rho^n \leqslant t_{\min} \right\rbrace $, there are at least 2 HPC jobs in the system at time $t_{\min}$ so the single activity server cannot be idling at that time.   Thus
$$
\mathbb{P}\left(\tau_1-t_{\min}> e_{\max}^n\right) \leqslant \mathbb{P}\left(\tau^n_r\leqslant t_0\right)+ \mathbb{P}\left(\rho^n > t_{\min}\right).
$$
By Lemma \ref{lem:nonidling2}, $\mathbb{P}\left(\tau^n_r\leqslant t_0\right)$ converges to zero.
In addition,
\begin{align*}
\mathbb{P}\left(\rho^n > t_{\min}\right) &=\mathbb{P}\left(\rho^n > t_{\min}, t_{\min}>\tilde{\tau}\right) +\mathbb{P}\left(\rho^n > t_{\min}, t_{\min} \leqslant \tilde{\tau}\right)\\
&\leqslant \mathbb{P}\left(\rho^n > \tilde{\tau} \right) + \mathbb{P}\left( t_{\min} \leqslant \tilde{\tau}\right)\\
&\leqslant \mathbb{P}\left(\rho^n > \tilde{\tau}, \tau^n_c \geqslant t_0 \right) + \mathbb{P}\left(\rho^n > \tilde{\tau}, \tau^n_c < t_0 \right) +
\mathbb{P}\left(\hat{X}_2^n(t_{\min}) < \alpha_2 z^*/2 \right)  \\
&\leqslant \mathbb{P}\left(\rho^n > \tilde{\tau}, \tau^n_c \geqslant t_0 \right) +2\mathbb{P}\left( \tau^n_c < t_0 \right).
\end{align*}
where $\tilde{\tau}$ is defined in the proof of Lemma \ref{lem:nonidling1}. 
Both terms converge to zero, the first by Lemma \ref{lem:nonidling2}, and the second by Proposition \ref{lem:ssc}.
This proves \eqref{eq:unifbound}.

\textbf{	Cases 2(c) and 2(d):} The current mode changes after the first service completion or arrival of a low or high priority job after $t_{\min}$ so under $\left\lbrace t_{\max} \leqslant  t_0 \right\rbrace$,
$$\tau_1-t_{\min}\leqslant e_{\max}^n.$$
Thus \eqref{eq:unifbound} follows from Lemma \ref{lem:1service}.
\qed

\appendix
\section{Solution of the HJB equation and free boundary point}
\beginsec
\manualnames{A}
\label{app:a}
In this appendix we present the expression found in
\cite[Section 5.3]{Sheng78} for the
solution to the HJB equation. It includes an equation
that uniquely characterizes the free boundary point (or switching point) $z^*$
in the dual mode case.
Recall that it is  assumed in \cite{Sheng78}, without loss of generality,
that $b_1\geqslant b_2$. 
We assume further, for simplicity, (and, again, without loss of generality) that if $b_1=b_2$ then $\sigma_1 \geqslant \sigma_2$.
Note that, with these indexing assumptions,  $m=2$ in  whichever of the complementary conditions \eqref{90} or  \eqref{91} that holds.
Throughout this section, denote
the unique classical solution to \eqref{14}--\eqref{14+}
by $u(x)$, $x\in\R_+$, and let $x$ serve as the initial condition for the WCP,
which elsewhere in this paper is denoted by $z$. Let
\begin{align*}
\beta=\dfrac{b_1+\sqrt{b_1^2+2\gamma \sigma_1^2}}{\sigma^2_1},
\qquad
\rho=\dfrac{b_2+\sqrt{b_2^2+2\gamma \sigma_2^2}}{\sigma^2_2},
\qquad
\nu=\dfrac{b_1-\sqrt{b_1^2+2\gamma\sigma^2_1}}{\sigma_1^2}.
\end{align*}

\begin{theorem}[{\cite[Section 5.3]{Sheng78}}]
\label{prop:optcostbcp}
Under \eqref{90},
$u(x)=\frac{x}{\gamma}+\frac{b_2}{\gamma^2}+\frac{1}{\gamma \rho}e^{-\rho x}$,
$x\in\R_+$.
\end{theorem}

Next, consider condition \eqref{91}. Because $b_1$ and $b_2$ are distinct
in this case, we have $b_1>b_2$.
The policy \eqref{p01} from \S\ref{sec4} corresponding to switching
at $z$ is given in the present notation
by $\bar\xi_z(x)=\xi^{*,1}\one_{x\le z}+\xi^{*,2}\one_{x>z}$.
Let $\frS^{(2)}_z$ be the admissible control system from Lemma \ref{lem3}.2,
with a generic switching point $z$ in place of the specific $z^*$.
Let the corresponding expression $J_{\rm WCP}(x,\frS^{(2)}_z)$,
which is nothing but the cost associated with the switching policy $\bar\xi_z$,
be denoted by $J(x,\bar\xi^z)$. Following is an expression
for this cost.
For $z>0$, let
\begin{align*}
{A}(z)&=\dfrac{\nu \beta (e^{(\nu-\beta )z}+e^{-(\nu-\beta )z})}{\nu \beta (e^{(\nu-\beta )z}+e^{-(\nu-\beta )z})+\rho e^{-\rho z}(e^{-\nu z}-e^{-\beta z})(\nu e^{\nu z}-\beta e^{\beta z})},\\
{B}(z)&=\dfrac{\rho e^{-\rho z}(e^{-\nu z}-e^{-\beta z})(\nu e^{\nu z}-\beta e^{\beta z})}{\nu \beta (e^{(\nu-\beta )z}+e^{-(\nu-\beta )z})+\rho e^{-\rho z}(e^{-\nu z}-e^{-\beta z})(\nu e^{\nu z}-\beta e^{\beta z})},\\
{C}(z)&=\dfrac{(\nu-\beta)(e^{-\nu z}-e^{-\beta z})}{\nu \beta (e^{(\nu-\beta )z}+e^{-(\nu-\beta )z})+\rho e^{-\rho z}(e^{-\nu z}-e^{-\beta z})(\nu e^{\nu z}-\beta e^{\beta z})},\\
{D}(z)&={A}(z)+{B}(z)+\rho e^{-\rho z}{C}(z),\\
{E}(z)&=\rho e^{-\rho z}{C}(z),\\
{F}(z)&=\dfrac{\beta-\nu}{\nu e^{-\beta z}-\beta e^{-\nu z}}.
\end{align*}
($F(\cdot)$ is not to be confused with the process $F$ defined in the body of the paper).
Then, for $x\leqslant z$,
\begin{multline}\label{eq:costlow}
J(x,\bar\xi^z)=\dfrac{x}{\gamma}+\dfrac{b_1\left[(e^{-\nu z}-e^{-\beta z})+({A}(z)-1)(e^{-\nu x}-e^{-\beta x})+({D}(z)-1)(e^{-\nu z-\beta x}-e^{-\beta z-\nu x})\right]}{\gamma^2(e^{-\nu z}-e^{-\beta z})}\\
+\dfrac{b_2\left[{B}(z)(e^{-\nu x}-e^{-\beta x})-{E}(z)(e^{-\nu z-\beta x}-e^{-\beta z-\nu x})\right]}{\gamma^2(e^{-\nu z}-e^{-\beta z})}\\
+\dfrac{{C}(z)(e^{-\nu x}-e^{-\beta x})-F(z)(e^{-\nu z-\beta x}-e^{-\beta z-\nu x})}{\gamma (e^{-\nu z}-e^{-\beta z})},
\end{multline}
and for $x\geqslant z$,
\begin{equation}\label{eq:costhigh}
J(x,\bar\xi^z)=\dfrac{x}{\gamma}+\dfrac{b_1{A}(z)}{\gamma^2e^{\rho x}}+\dfrac{b_2}{\gamma^2}(1-e^{-\rho x}+e^{-\rho x}{B}(z))+\dfrac{{C}(z)-1}{\gamma e^{\rho x}}.
\end{equation}
It is here where the principle of smooth fit is applied. For the cost to be
$C^2$ (in $x$), it must satisfy $J''(z-,\bar\xi^z)=J''(z+,\bar\xi^z)$.
Using the expressions \eqref{eq:costlow} and \eqref{eq:costhigh}, this condition
can be translated to the following equation
\begin{multline}\label{59}
0=\left(\dfrac{b_1-b_2}{\gamma^2}\right)\left[\dfrac{({A}(z)-1)(\beta e^{-\beta z}-\nu e^{-\nu z})}{e^{-\nu z}-e^{-\beta z}}-\rho^2 e^{-\rho z}{A}(z)+\dfrac{(\nu-\beta)({D}(z)-1) e^{-(\nu+\beta)z}}{e^{-\nu z}-e^{-\beta z}}\right]\\
+\frac{1}{\gamma}\left[\dfrac{{C}(z)(\beta e^{-\beta z}-\nu e^{-\nu z})}{e^{-\nu z}-e^{-\beta z}}-\rho^2 e^{-\rho z}{C}(z)+\dfrac{(\beta-\nu)F(z) e^{-(\nu+\beta)z}}{e^{-\nu z}-e^{-\beta z}}\right].
\end{multline}
\begin{theorem}[{\cite[Section 5.3]{Sheng78}}]
Let \eqref{91} hold and assume $b_1-b_2\ne\gamma$. Then \eqref{59}
has a unique solution $z^*\in(0,\iy)$. Moreover,
$u(x)=J(x,\bar\xi^{z^*})$, $x\in\R_+$.
\end{theorem}

\section{Symmetry conditions}
\beginsec
\manualnames{B}
\label{app:b}

The following result is related to Remark \ref{prop:noswitch}.
\begin{lemma}\label{lem:sameb}
1. If either \eqref{sy1} or \eqref{sy2} holds then
$b_1=b_2$. In particular, \eqref{90} holds.
\\
2. If \eqref{sy3} holds then $\sig_1=\sig_2$.
In particular, \eqref{90} holds.
\end{lemma}

\noi{\bf Proof.} 1.
Recall the expressions for $b_1$ and $b_2$,
$$
b_m=b(\xi^{*,m})=\sum_i\frac{\hat\la_i-\sum_k\hat\mu_{ik}\xi^{*,m}_{ik}}{\al_i}.
$$
The difference between $b_1$ and $b_2$ is thus the difference between $\gamma_m:=\sum_{i,k}\frac{\hat\mu_{ik}\xi^{*,m}_{ik}}{\al_i}$.
For $\xi^{*,1}$, we distinguish these cases: either $\frac{\lambda_1}{\alpha_1}> \beta_2$ or $\frac{\lambda_1}{\alpha_1}< \beta_2$. 	For $\xi^{*,2}$, we distinguish these cases: either $\frac{\lambda_1}{\alpha_1}< \beta_1$ or $\frac{\lambda_1}{\alpha_1}> \beta_1$. We will see that in each of the four cases $b_1=b_2$.
\begin{align*}
\frac{\lambda_1}{\alpha_1}> \beta_2:\quad \gamma_1&=\dfrac{1}{\alpha_1}\left[\hat{\mu}_{12}+\hat{\mu}_{11}(\dfrac{\lambda_1}{\alpha_1\beta_1}-\dfrac{\beta_2}{\beta_1})\right]+\dfrac{\hat{\mu}_{21}}{\alpha_2}(1-\dfrac{\lambda_1}{\alpha_1\beta_1}+\dfrac{\beta_2}{\beta_1}),\\
\frac{\lambda_1}{\alpha_1}< \beta_2:\quad \gamma_1&=\dfrac{\hat{\mu}_{12}\lambda_1}{\alpha_1^2\beta_2}+\dfrac{1}{\alpha_2}\left[\hat{\mu}_{21}+\hat{\mu}_{22}(1-\dfrac{\lambda_1}{\alpha_1\beta_2})\right],\\
\frac{\lambda_1}{\alpha_1}< \beta_1:\quad \gamma_2&=\dfrac{\hat{\mu}_{11}\lambda_1}{\alpha_1^2\beta_1}+\dfrac{1}{\alpha_2}\left[\hat{\mu}_{22}+\hat{\mu}_{21}(1-\dfrac{\lambda_1}{\alpha_1\beta_1})\right],\\
\frac{\lambda_1}{\alpha_1}> \beta_1:\quad \gamma_2&=\dfrac{1}{\alpha_1}\left[\hat{\mu}_{11}+\hat{\mu}_{12}(\dfrac{\lambda_1}{\alpha_1\beta_2}-\dfrac{\beta_1}{\beta_2})\right]+\dfrac{\hat{\mu}_{22}}{\alpha_2}(1-\dfrac{\lambda_1}{\alpha_1\beta_2}+\dfrac{\beta_1}{\beta_2}).\\
\end{align*}
We now take the difference for each pair:
\begin{align*}
\frac{\lambda_1}{\alpha_1}> \beta_2\, \&\,\frac{\lambda_1}{\alpha_1}< \beta_1 :
\quad b_1-b_2&=\dfrac{\hat{\mu}_{12}-\hat{\mu}_{11}\frac{\beta_2}{\beta_1}}{\alpha_1}-\dfrac{\hat{\mu}_{22}-\hat{\mu}_{21}\frac{\beta_2}{\beta_1}}{\alpha_2},\\
\frac{\lambda_1}{\alpha_1}> \beta_2\,\&\,\frac{\lambda_1}{\alpha_1}> \beta_1:
\quad b_1-b_2&=(1-\dfrac{\lambda_1}{\alpha_1\beta_1}+\dfrac{\beta_2}{\beta_1})\left[\dfrac{1}{\alpha_1}(\hat{\mu}_{12}\dfrac{\beta_1}{\beta_2}-\hat{\mu}_{11})+\dfrac{1}{\alpha_2}(\hat{\mu}_{21}-\hat{\mu}_{22}\dfrac{\beta_1}{\beta_2})\right],\\
\frac{\lambda_1}{\alpha_1}< \beta_2\, \&\, \frac{\lambda_1}{\alpha_1}< \beta_1 :
\quad b_1-b_2&=\dfrac{\lambda_1}{\alpha_1^2}\left[\dfrac{\hat{\mu}_{12}}{\beta_2}-\dfrac{\hat{\mu}_{11}}{\beta_1}\right]+\dfrac{\lambda_1}{\alpha_1\alpha_2}\left[\dfrac{\hat{\mu}_{21}}{\beta_1}-\dfrac{\hat{\mu}_{22}}{\beta_2}\right],\\
\frac{\lambda_1}{\alpha_1}< \beta_2\, \&\, \frac{\lambda_1}{\alpha_1}> \beta_1:
\quad b_1-b_2&=\dfrac{\hat{\mu}_{12}\frac{\beta_1}{\beta_2}-\hat{\mu}_{11}}{\alpha_1}+\dfrac{\hat{\mu}_{21}-\hat{\mu}_{22}\frac{\beta_1}{\beta_2}}{\alpha_2}.
\end{align*}
Now, if
$\frac{\hat{\mu}_{i1}}{\beta_1}=\frac{\hat{\mu}_{i,2}}{\beta_2}$, $i=1,2,$
we get 
$$\hat{\mu}_{11}\frac{\beta_2}{\beta_1}-\hat{\mu}_{12}=\hat{\mu}_{22}-\hat{\mu}_{21}\frac{\beta_2}{\beta_1}=0,$$
and consequently $b_1=b_2$ as claimed. If
$\frac{\hat{\mu}_{1k}}{\alpha_1}=\frac{\hat{\mu}_{2k}}{\alpha_2}$,
$k=1,2,$
it is not hard to see that again the expressions can be rewritten
with a different factorization to get $b_1=b_2$.

2. Under \eqref{sy3}, denote $C_i=C_{S_{i1}}=C_{S_{i2}}$.
Then for $\xi \in \slp$,
\begin{align*}
\sigma(\xi)^2&=\sum_{i}\dfrac{\sigma^2_{A,i}+\sum_k \sigma^2_{S_{ik}}\xi_{ik}}{\alpha_i^2}\\
&=\sum_{i}\dfrac{\sigma^2_{A,i}+\sum_k C^2_{i}\mu_{ik}\xi_{ik}}{\alpha_i^2}\\
&=\sum_i\dfrac{\sigma^2_{A,i}+C^2_{i}\lambda_i}{\alpha_i^2},
\end{align*}
where the last equality follows from \eqref{01}. Hence $\sig_1=\sig_2$
as claimed.
\qed

\paragraph{Acknowledgement.}
We are very grateful to Cristina Costantini for providing us with
the proof of Lemma \ref{lem3} Part 3.
RA is supported by ISF grant 1035/20.

\footnotesize

\bibliographystyle{is-abbrv}

\bibliography{main}

\end{document}